\documentclass[11pt,twoside]{article}
\usepackage{amsmath,amssymb,euscript,theorem,epsfig,psfrag}

\numberwithin{equation}{section}
\newtheorem{theorem}{Theorem}[section]
\newtheorem{corollary}[theorem]{Corollary}
\newtheorem{lemma}[theorem]{Lemma}
\newtheorem{proposition}[theorem]{Proposition}
\newtheorem{summary}[theorem]{Summary}
{\theorembodyfont{\rmfamily} \newtheorem{rem}[theorem]{Remark}}
\newenvironment{proof}{{\bf Proof}.\ }{ \hfill $\square$}

\newcommand{\bke}[1]{\left ( #1 \right )}
\newcommand{\bkt}[1]{\left [ #1 \right ]}
\newcommand{\bket}[1]{\left \{ #1 \right \}}
\newcommand{\norm}[1]{\left \| #1 \right \|}

\newcommand{\e}{\varepsilon}

\renewcommand{\L}{{\cal L}}
\newcommand{\la}{\lambda}
\newcommand{\al}{\alpha}
\newcommand{\pd}{\partial}

\newcommand{\R}{\mathbb{R}}
\newcommand{\C}{\mathbb{C}}
\renewcommand{\Re}{\mathop{\mathrm{Re}}}
\renewcommand{\Im}{\mathop{\mathrm{Im}}}

\newcommand{\donothing}[1]{}

\newcommand{\de}{\delta}
\newcommand{\si}{\sigma}
\newcommand{\Lp}{L_{+}}
\newcommand{\Lm}{L_{-}}
\newcommand{\bP}{\bar{P}}
\newcommand{\ra}{\rightarrow}
\newcommand{\mat}[1]{\begin{bmatrix} #1 \end{bmatrix}}
\newcommand{\vect}[1]{\left[ \begin{smallmatrix} #1 \end{smallmatrix} \right]}
\newcommand{\pmax}{p_{\mathrm{max}}}
\newcommand{\be}{\beta}
\newcommand{\ph}{\varphi}
\newcommand{\ka}{\kappa}
\newcommand{\sech}{{\,\mathrm{sech}}}
\renewcommand{\th}{\theta}
\newcommand{\pispc} {\delta_{p_c}^p}

\def\pcirc{{\text{$\bigcirc$ \hspace{-0.41cm} {\tiny\rm p}~}}}
\def\ppcirc{{\text{$\bigcirc$ \hspace{-0.43cm} {\tiny$\gamma$}~}}}
\def\lb{[\hspace{-0.05cm}[}
\def\rb{]\hspace{-0.05cm}]}
\def\diag{{\rm diag}}

\begin{document}

\title{
Spectra of Linearized Operators for NLS Solitary Waves}

\author{Shu-Ming Chang\thanks
{Department of Mathematics, National Tsing Hua University, Hsinchu, Taiwan
(schang@mx.nthu.edu.tw)}, \quad
Stephen Gustafson\thanks
{Department of Mathematics, University of British Columbia,
Vancouver, BC V6T1Z2, Canada
(gustaf@math.ubc.ca)}, \quad
Kenji Nakanishi\thanks {Department of Mathematics, Kyoto University,
Kyoto 606-8502, Japan (n-kenji@math.kyoto-u.ac.jp)}, \quad
Tai-Peng Tsai\thanks
{Department of Mathematics, University of British Columbia,
Vancouver, BC V6T1Z2, Canada
(ttsai@math.ubc.ca)}
}

\date{}

\maketitle

{\bf Abstract}. Nonlinear Schr\"odinger (NLS) equations with focusing power
nonlinearities have solitary wave solutions.
The spectra of the linearized operators around these solitary waves
are intimately connected to stability properties of the solitary waves,
and to the long-time dynamics of solutions of (NLS).
We study these spectra in detail, both analytically
and numerically.
\medskip


\medskip


\medskip

\section{Introduction}\label{sec1}
Consider the nonlinear Schr\"odinger equation (NLS) with focusing
power nonlinearity,
\begin{equation}\label{eq1-1}
  i \pd_t {\psi} = -\Delta \psi - |\psi|^{p-1}\psi,
\end{equation}
where $\psi(t,x):\R \times \R^n \to \C$ and $1<p<\infty$.
Such equations arise in many physical settings,
including nonlinear optics, water waves, and quantum physics.
Mathematically, nonlinear Schr\"odinger equations with various
nonlinearities are studied as basic models of nonlinear dispersive
phenomena. In this paper, we stick to the case of a pure power
nonlinearity for the sake of simplicity.

For a certain range of the power $p$ (see below), the NLS \eqref{eq1-1}
has special solutions, of the form $\psi(t,x) = Q(x)\, e^{it}$. These are
called {\it solitary waves}.
The aim of this paper is to study the spectra of the
{\it linearized operators} which arise when~\eqref{eq1-1} is linearized
around solitary waves. The main motivation for this study is
that properties of these spectra are intimately related to the problem of the
stability (orbital and asymptotic) of these solitary waves,
and to the long-time dynamics of solutions of NLS.

Let us begin by recalling some well-known facts about \eqref{eq1-1}.
Standard references include \cite{Ca, St, SS}.  Many basic results
on the linearized operators we study here were proved by
Weinstein \cite{W83,W85}.  The Cauchy (initial value) problem for
equation \eqref{eq1-1} is locally (in time) well-posed in
$H^1(\R^n)$ if $1<p<\pmax$, where
\begin{equation*}
\pmax := 1+4/(n-2) \quad \text{ if } n \ge 3; \quad
\pmax := \infty \quad \text{ if } n=1,2.
\end{equation*}
Moreover, if $1<p<p_c$ where
\begin{equation*}
  p_c := 1+ 4/n,
\end{equation*}
the problem is globally well-posed. For $p\ge p_c$, there exist
solutions whose $H^1$-norms go to $\infty$ ({\it blow up}) in finite
time. In this paper, the cases $p<p_c$, $p=p_c$ and $p>p_c$ are called
{\it sub-critical}, {\it critical}, and {\it super-critical}, respectively.

The set of all solutions of \eqref{eq1-1} is invariant under the
symmetries of translation, rotation, phase, Galilean transform and
scaling: if $\psi(t,x)$ is a solution, then so is
\begin{equation*}
  \widetilde \psi(t,x) := \la^{2/(p-1)} \psi \bke{ \la^2 t, \ \la R x - \la^2
  tv -x_0} \, \exp \bket{ i\bkt{ \frac{\la Rx\cdot v}2 - \frac{ \la^2 t
  v^2}4 + \gamma_0}}
\end{equation*}
for any constant $x_0, v\in \R^n$, $\la >0$, $\gamma_0\in \R$ and $R
\in O(n)$.  When $p=p_c$, there is an additional symmetry
called the ``pseudo-conformal transform'' (see \cite[p.35]{SS}).

We are interested here in solutions of \eqref{eq1-1}
of the form
\begin{equation}\label{eq1-3}
  \psi(t,x) = Q(x)\, e^{it}
\end{equation}
where $Q(x)$ must therefore satisfy the nonlinear elliptic equation
\begin{equation} \label{Q.eq}
  -\Delta Q- |Q|^{p-1}Q = -Q.
\end{equation}
Any such solution generates a family of solutions by the
above-mentioned symmetries, called {\it solitary waves}.
Solitary waves are special examples of {\it nonlinear bound states},
which, roughly speaking, are solutions that are spatially
localized for all time. More precisely, one could define nonlinear
bound states to be solutions $\psi(t,x)$ which are
{\it non-dispersive} in the sense that
\begin{equation*}
\sup_{t \in \R} \inf_{ x_0 \in \R^n} \norm{ |x| \psi(t,x-x_0)
}_{L^2_x(\R^n)} < \infty.
\end{equation*}

Testing \eqref{Q.eq} with $\bar Q$ and $x.\nabla \bar Q$ and taking
real parts, one arrives at the {\it Pohozaev identity} (\cite{Poho})
\begin{equation} \label{pohozaev}
  \frac{1}{2}\int |Q|^2 = b \frac{1}{p+1} \int |Q|^{p+1}, \quad
  \frac{1}{2} \int |\nabla Q|^2  = a \frac{1}{p+1} \int |Q|^{p+1}
\end{equation}
where
\[
  a = \frac {n(p-1)}4, \quad b = \frac {n+2-(n-2)p}4.
\]
The coefficients $a$ and $b$ must be positive, and hence a necessary
condition for existence of non-trivial solutions is $p \in (1, \pmax)$.

For $p \in (1,\pmax)$, and for all space dimensions, there exists
at least one non-trivial radial solution
$Q(x) = Q(|x|)$ of \eqref{Q.eq} (existence goes back to \cite{Poho}).
This solution, called a {\it nonlinear ground state}, is smooth,
decreases monotonically as a function of $|x|$, decays exponentially
at infinity, and can be taken to be positive: $Q(x) > 0$. It is the
unique positive solution. (See~\cite{SS} for references for
the various existence and uniqueness results for various
nonlinearities.) The
ground state can be obtained as the minimizer of several different variational
problems. One such result we shall briefly use later
is that, for all $n\ge 1$ and $p \in (1,\pmax)$,
the ground state minimizes the
Gagliardo-Nirenberg quotient
\begin{equation}
\label{varprob}
  J[u] := \frac{\bke{\int |\nabla u|^2 }^a \bke{\int u^2 }^b}{\int
  u^{p+1} }
\end{equation}
among nonzero $H^1(\R^n)$ radial functions (Weinstein \cite{W83}).

For $n=1$, the ground state
is the unique $H^1(\R)$-solution of \eqref{Q.eq} up to translation and
phase \cite[p.259, Theorem 8.1.6]{Ca}.
For $n \geq 2$, this is not the case: there are countably
infinitely many radial solutions (still real-valued),
denoted in this paper by $Q_{0,k,p}(x)$, $k=0,1,2,3,\ldots$,
each with exactly $k$ positive zeros as a function
of $|x|$ (Strauss \cite{Str77}; see also \cite{BL}).
In this notation, $Q_{0,0,p}$ is the ground state.

There are also non-radial (and complex-valued) solutions, for example
those suggested by P. L. Lions \cite{Lions} with non-zero angular momenta,
\begin{equation*}
n=2, \quad Q = \phi(r) \,e^{i m \th}, \quad \text{ in polar
coordinates } r, \th;
\end{equation*}
\begin{equation*}
n=3, \quad Q = \phi(r,x_3)\, e^{i m \th}, \quad \text{ in cylindrical
coordinates } r, \th, x_3,
\end{equation*}
and similarly defined for $n \ge 4$.  When $n=2$, some of these
solutions are denoted here by $Q_{m,k,p}$ with $p \in (1,\pmax)$ and
$k=0,1,2,\ldots$ denoting their numbers of positive zeros.  
See Section~\ref{sec4} for more details.

We will refer to all the solitary waves generated by $Q_{0,0,p}$ as
{\it nonlinear ground states}, and all others as {\it nonlinear
excited states}. We are not aware of a complete characterization of
all solutions of \eqref{Q.eq}, or of \eqref{eq1-1}. For example,
the uniqueness of $Q_{m,k,p}$ with $m,k \ge 1$ is
apparently open.  Also, we do not know if there are ``breather''
solutions, analogous to those of the generalized KdV equations.  In
this paper we will mainly study radial solutions (and in particular
the ground state), but we will also briefly consider non-radial
solutions numerically in Section~\ref{sec4}.

To study the stability of a solitary wave solution~\eqref{eq1-3},
one considers solutions of (NLS) of the form
\begin{equation} \label{eq1-5}
  \psi(t,x) = [Q(x) + h(t,x)]\, e^{it} .
\end{equation}
For simplicity, let $Q=Q_{0,0,p}$ be the ground state for the
remainder of this introduction (see Section~\ref{sec4} for the
general case). The perturbation $h(t,x)$ satisfies an equation
\begin{equation} \label{eq1-6}
\pd_t h = \L h + (\text{nonlinear terms})
\end{equation}
where $\L$ is the {\it linearized operator} around $Q$:
\begin{equation} \label{eq1-9}
  \L h = -i \bket{(-\Delta +1 - Q^{p-1})h  - \tfrac{p-1}2\, Q^{p-1}(h+\bar h)}.
\end{equation}
It is convenient to write $\L$ as a matrix operator acting on
$\mat{\Re h \\ \Im h}$,
\begin{equation} \label{eq1-7}
  \L = \mat{ 0 & \Lm \\ -\Lp & 0}
\end{equation}
where
\begin{equation} \label{eq1-11}
\Lp = -\Delta + 1 - pQ^{p-1}, \qquad  \Lm = -\Delta + 1 - Q^{p-1}.
\end{equation}

Clearly the operators $\Lm$ and $\Lp$ play a central role
in the stability theory. They are self-adjoint Schr\"odinger operators
with continuous spectrum $[1,\infty)$, and with finitely many
eigenvalues below $1$. In fact, when $Q$ is the ground state,
it is easy to see that
$\Lm$ is a nonnegative operator, while $\Lp$ has exactly one negative
eigenvalue (these facts follow from Lemma~\ref{TH2-1} below).

Because of its connection to the stability problem,
the object of interest to us in this paper is the spectrum
of the non-self-adjoint operator $\L$. The simplest properties
of this spectrum are
\begin{enumerate}
\item for all $p \in (1,\pmax)$, $0$ is an
eigenvalue of $\L$.
\item the set $\Sigma_c := \bket{ ir: r \in \R, |r|\ge 1}$
is the continuous spectrum of $\L$.
\end{enumerate}
(See the next section for the first statement. The second
is easily checked.)

It is well-known that the exponent $p=p_c$ is critical for
stability of the ground state solitary wave
(as well as for blow-up of solutions). For $p < p_c$ the ground state
is {\it orbitally stable}, while for $p \geq p_c$ it is
unstable (see \cite{W86,GSS}). These facts have immediate spectral
counterparts: for $p \in (1,p_c]$, all
eigenvalues of $\L$ are purely imaginary, while for $p \in (p_c,\pmax)$,
$\L$ has at least one eigenvalue with positive real part.

The goal of this paper is to get a more detailed understanding of the
spectrum of $\L$, using both analytical and numerical techniques.  See
\cite{Gril,ComPel,CPV,Dem-Sch} for related work.  This finer
information is essential for understanding the long-time dynamics of
solutions of (NLS), for example: (i) to prove asymptotic (rather than
simply orbital) stability, one often assumes either the linearized
operator $\L$ has no nonzero eigenvalue, or its nonzero eigenvalues
are $\pm ri$ with $0.5<r<1$. These assumptions need to be verified;
(ii) to determine the rate of relaxation to stable
solitary waves when there is a unique pair of nonzero eigenvalues $\pm
ri$ with $0<r<1$.  Heuristic arguments suggest that $[1/r]$, the
smallest integer no longer that $1/r$, may decide the rate; (iii) to
construct stable manifolds of unstable solitary waves, one needs to
know if there are eigenvalues which are not purely imaginary, and to
find their locations. These are highly active areas of current
research, see e.g.~\cite{GNT,KZ,Schlag04} and the references therein.

Interesting questions with direct relevance to these stability-type
problems include:

\begin{itemize}
\item[(i)]
Can one determine (or estimate) the number and locations
of the eigenvalues of $\L$ lying on the segment between $0$ and $i$?

\item[(ii)]
Can $\pm i$, the thresholds of the continuous spectrum $\Sigma_c$,
be eigenvalues or resonances?

\item[(iii)]
Can eigenvalues be embedded inside the continuous spectrum?

\item[(iv)] Can the linearized operator have eigenvalues with non-zero
real {\it and} imaginary parts (this is already known not to happen
for the ground state -- see the next section -- and so we pose this
question with excited states in mind).

\item[(v)]
Are there bifurcations, as $p$ varies, of pairs of purely
imaginary eigenvalues into pairs of eigenvalues with non-zero
real part (a stability/instability transition)?
\end{itemize}

\medskip

The detailed discussion of the numerical methods is postponed to the
Appendix. Roughly speaking, we first compute the nonlinear ground
state by iteration and renormalization, and then compute the spectra
of various suitably discretized linear operators.

Let us now summarize the main results and observations of this paper.

\begin{enumerate}
\item
{\bf Numerics for spectra.} When $Q$ is the ground state, we compute
numerically the spectra of $\L$, $L_+$ and $L_-$ as functions of $p$,
see Figures~\ref{fig:cL1d}--\ref{fig:cL2d-m2}. In these figures, the
horizontal axis is the logarithm of $p-1$.  The vertical axis is:
Solid lines are purely imaginary eigenvalues of $\L$ (without $i$) for
$ \in (1,p_c)$; dashed lines are real eigenvalues of $\L$; dotted line
are eigenvalues of $L_+$; dashdot lines are eigenvalues of $L_-$.  We
have ignored imaginary eigenvalues of the discretized operators with
modulus greater than one, which correspond to the continuous spectra
of the original operators.  Figure~\ref{fig:cL1d} is the
one-dimensional case.  Figures~\ref{fig:cL2d} and \ref{fig:cL3d} are
the spectra of these operators restricted to radial functions, for
space dimensions $n=2$ and $3$.  Figures~\ref{fig:cL2d-m1} and
\ref{fig:cL2d-m2} are for $n=2$ and are the spectra restricted to
functions of the form $\phi(r)e^{i\th}$ and $\phi(r)e^{i2\th}$,
respectively. These pictures shed some light on questions (i), (iv),
and (v) above, and to a certain extent on question (ii).

Figures~\ref{fig:d2m1}--\ref{fig:bifm2} are concerned with the spectra
of excited states, see discussion below.

\begin{figure}
\centering \psfrag{x}{${\rm log_{10}}(p-1)$}\psfrag{t}{$p$}
\epsfig{file=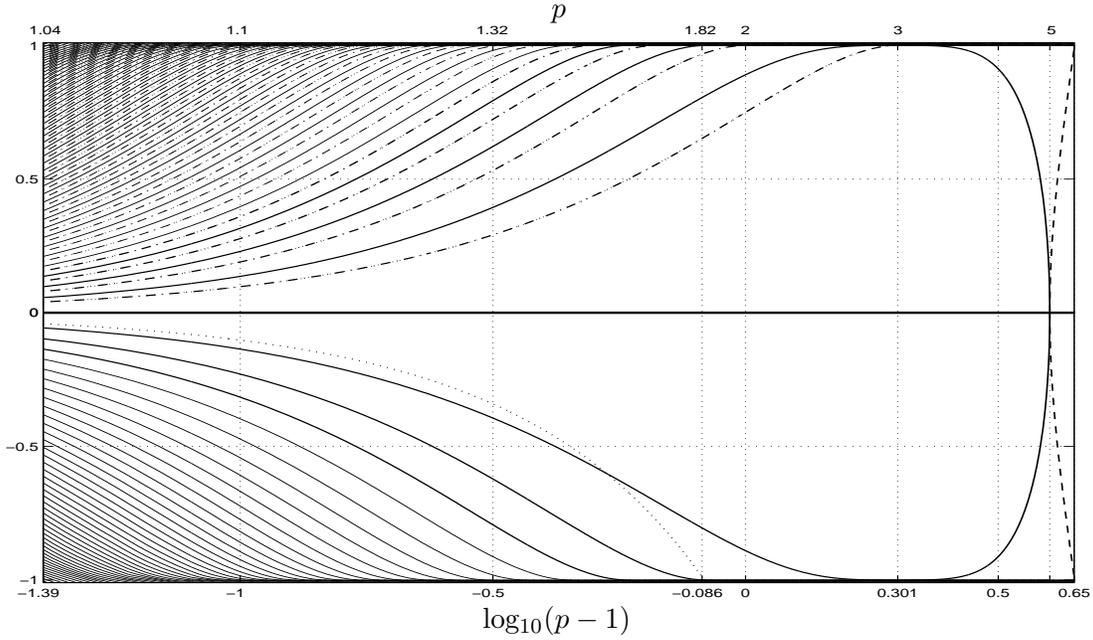,height=14.5cm,width=8cm, angle=270}
\caption{ Spectra of $\L$, $L_{+}$ and $L_{-}$ for $n=1$ with
logarithmic axis for the values of $p-1$. (solid line: purely
imaginary eigenvalues of $\L$; dashed line: real eigenvalues of
$\L$; dotted line: eigenvalues of $L_+$; dashdot line: eigenvalues
of $L_-$)} \label{fig:cL1d}
\end{figure}
%
\begin{figure}
\centering \psfrag{x}{${\rm log_{10}}(p-1)$}\psfrag{t}{$p$}
\epsfig{file=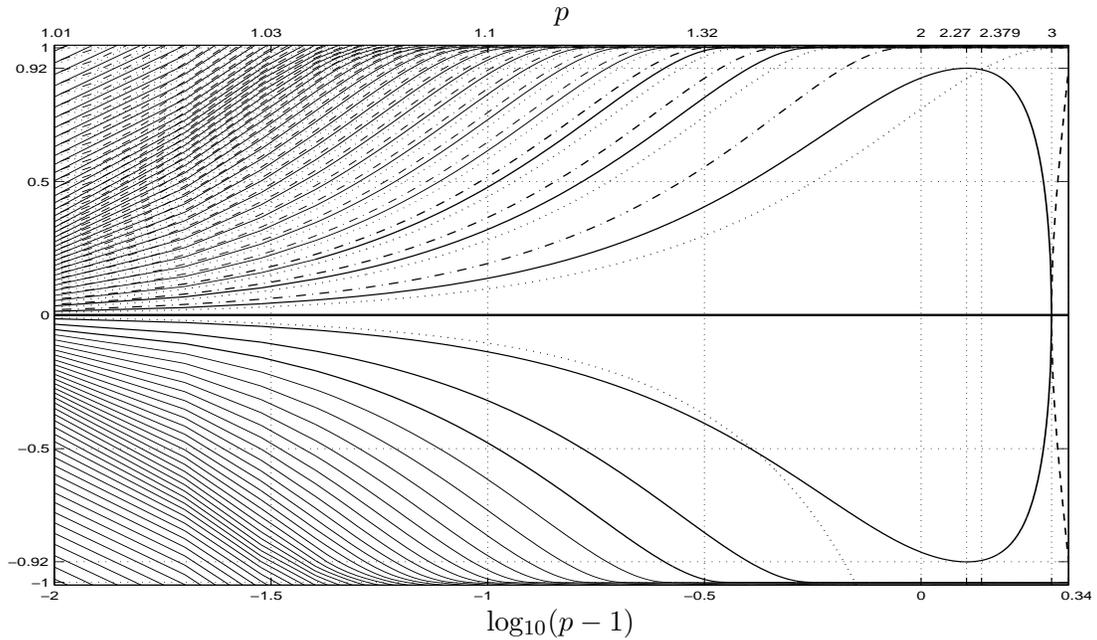,height=14.5cm,width=8cm, angle=270}
\caption{ Spectra of $\L$, $L_{+}$ and $L_{-}$ restricted to radial
functions in the two-dimensional space, with
logarithmic axis for the values of $p-1$. (solid line: purely
imaginary eigenvalues of $\L$; dashed line: real eigenvalues of
$\L$; dotted line: eigenvalues of $L_+$; dashdot line: eigenvalues
of $L_-$) } \label{fig:cL2d}
\end{figure}
%
\begin{figure}
\centering \psfrag{x}{${\rm
log_{10}}(p-1)$}\psfrag{t}{$p$}\psfrag{y}{$\frac{7}{3}$}
\epsfig{file=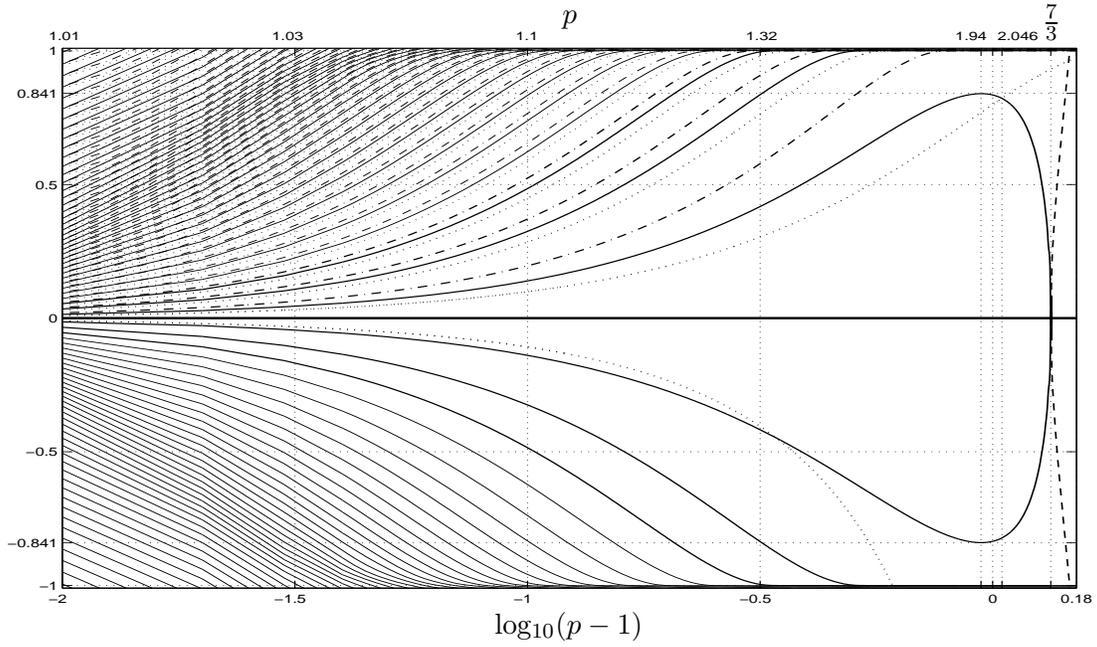,height=14.5cm,width=8cm, angle=270}
\caption{ Spectra of $\L$, $L_{+}$ and $L_{-}$ restricted to radial
functions in the three-dimensional space. } \label{fig:cL3d}
\end{figure}
%
%
\begin{figure}
\centering \psfrag{x}{${\rm log_{10}}(p-1)$}\psfrag{t}{$p$}
\epsfig{file=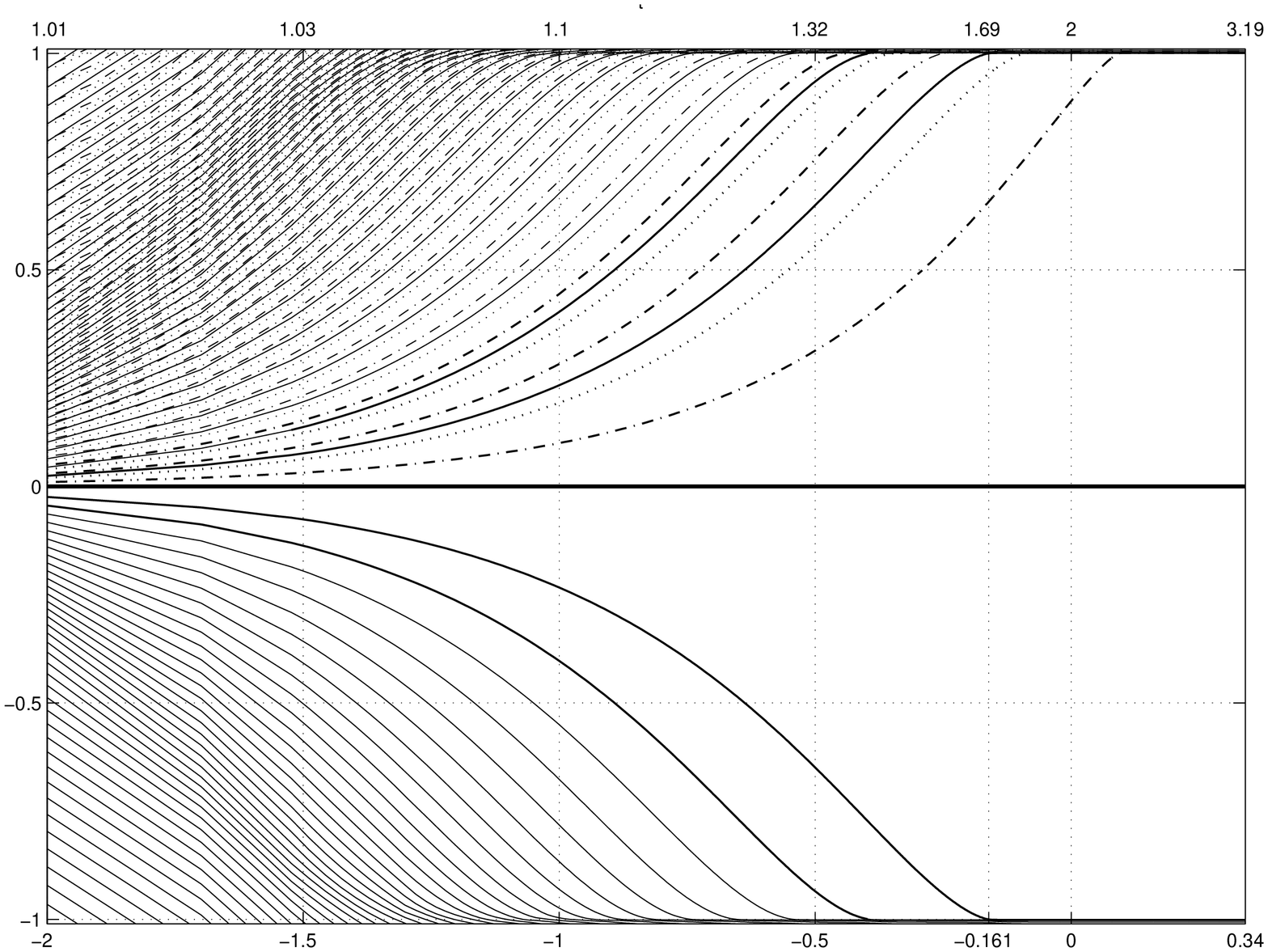,height=14.5cm,width=8cm, angle=270}
\caption{Spectra of $\L$, $L_{+}$ and $L_{-}$ restricted to
functions of the form $\phi(r) e^{i \th}$ in the two-dimensional
space.} \label{fig:cL2d-m1}
\end{figure}
%
\begin{figure}
\centering \psfrag{x}{${\rm log_{10}}(p-1)$}\psfrag{t}{$p$}
\epsfig{file=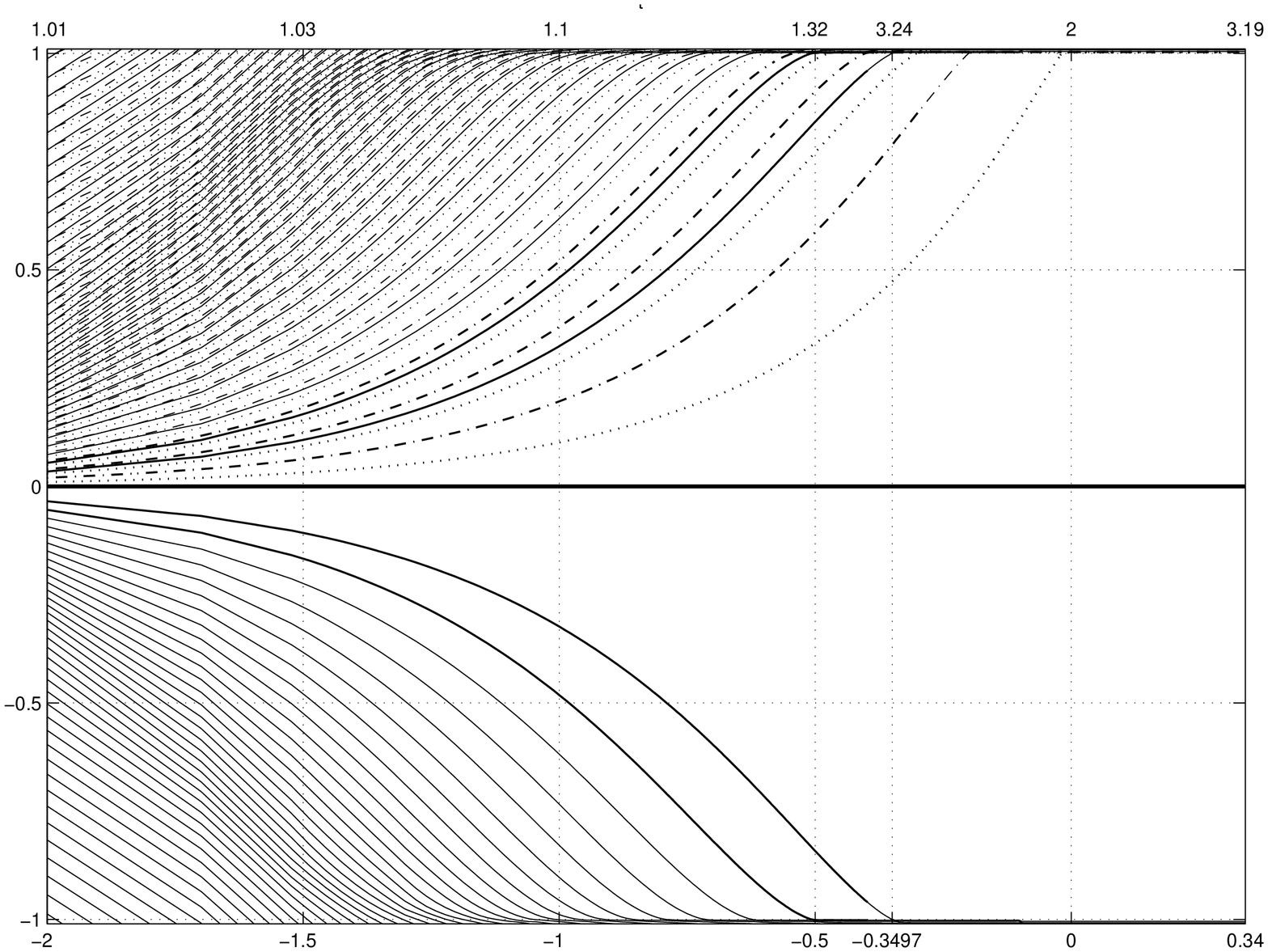,height=14.5cm,width=8cm, angle=270}
\caption{Spectra of $\L$, $L_{+}$ and $L_{-}$ restricted to
functions of the form $\phi(r) e^{i 2\th}$ in the two-dimensional
space.} \label{fig:cL2d-m2}
\end{figure}
%
%

\item
{\bf One-dimensional phenomena.} The case $n=1$ is the easiest case
to handle analytically. In Section~\ref{sec3}, we undertake a
detailed study of the one-dimensional problem, giving rigorous
proofs of a number of phenomena observed in Figure~\ref{fig:cL1d}.
One simple such phenomenon is the (actually
classical) fact that the eigenvalues of $\Lp$ and $\Lm$ exactly
coincide, with the exception of the first, negative, eigenvalue of
$\Lp$ (note that this appears to be a strictly one-dimensional
phenomenon: the eigenvalues of $L_+$ and $L_-$ are different for
$n\ge 2$, as Figures~\ref{fig:cL2d}--\ref{fig:cL2d-m2} indicate).
In fact, we are able to prove sufficiently precise upper and lower
bounds on the eigenvalues of $\L$ (lying outside the continuous
spectrum) to determine their number, and estimate their positions,
as functions of $p$ (see Theorem~\ref{th3-7}). We use two basic
techniques: an embedding of $\Lp$ and $\Lm$ into a hierarchy of
related operators, and a novel variational problem for the
eigenvalues, in terms of a 4-th order self-adjoint differential
operator (see Theorem~\ref{th3-5}). In this way, we get a fairly
complete answer to question (i) above for $n=1$.

\item
{\bf Variational characterization of eigenvalues}. We present
self-adjoint variational formulations of the eigenvalue problem for
$\L$ in any dimension (see Summary~\ref{TH2-5A}), including the
novel $n=1$ formulation mentioned above. In principle, these provide
a means of counting/estimating the eigenvalues of $\L$ (and hence
addressing question (i) above in higher dimensions), though we only
obtain detailed such information for $n=1$.

\item
{\bf Bifurcation at $p=p_c$.} In each of
Figures~\ref{fig:cL1d}--\ref{fig:cL3d}, a pair of purely imaginary
eigenvalues for $p < p_c$ appears to collide at $0$ at $p=p_c$, and
become a pair of real eigenvalues for $p > p_c$. This is exactly the
stability/instability transition for the ground state. We rigorously
verify this picture, determining analytically the spectrum near 0
for $p$ near $p_c$, and making concrete a bifurcation picture
suggested by M.~I.~Weinstein (personal communication): see
Theorem~\ref{th2-2}. This gives a partial answer to question (v)
above. It is worth pointing out that for $n=1$, the imaginary part
of the (purely imaginary) eigenvalue bifurcating for $p < p_c$ is
always larger than the third eigenvalue of $\Lp$ (the first is
negative and the second is zero) -- this is proved analytically in
Theorem~\ref{th3-7}. For $n \ge 2$, however, they intersect at $p
\approx 2.379$ for 2D and $ p \approx 2.046$ for 3D, 
(see Figures~\ref{fig:cL1d}--\ref{fig:cL3d}).


\item
{\bf Interlacing property.}
A numerical observation: in all the figures,
the adjacent eigenvalues of $\L$ seem each to bound an
eigenvalue of $\Lp$ and one of $\Lm$ (at least for $p$ small enough).
We are able to establish this ``interlacing'' property analytically
in dimension one (see Theorem~\ref{th3-7}).

\item
{\bf Threshold resonance}. An interesting fact observed numerically
(Figure~\ref{fig:cL1d}) is that, in the 1D case,
as $p \to 3$, one eigenvalue curve converges to $\pm i$, the
threshold of the continuous spectrum. One might suspect that, at
$p=3$, $\pm i$ corresponds to a resonance or embedded eigenvalue. It
is indeed a resonance: we find an explicit non spatially-decaying
``eigenfunction'', and show numerically in Section~\ref{sec3-7} that
the corresponding eigenfunctions converges, as $p \to 3$, to this
function. This observation addresses question (ii) above for $n=1$.

\item
{\bf Excited states}. In Section~\ref{sec4} we consider the spectra of
linearized operators around excited states with non-zero angular
momenta. We observe that, in addition to the bifurcation mentioned
above at $p=p_c$, there are complex eigenvalues which are neither real
nor purely imaginary (addressing question (iv) above; see
Figures~\ref{fig:d2m1}--\ref{fig:bifm2}), symmetric with respect to
both real and purely imaginary axes.  These complex eigenvalues also
come from bifurcation: as $p$ decreases, a quadruple of complex
eigenvalues will collide into the imaginary axis away from $0$, and
then split to 4 purely imaginary eigenvalues. It seems that all
eigenvalues lie on the imaginary axis for $p\in (1,p_*)$ for some
$p_*$ close to $1$. In other words, {\it numerically these excited
states are spectrally stable for $p$ close to $1$}. It is possible
that the numerical error increases enormously as $p \to 1_+$ due to
the artificial boundary condition, since the spectrum is approaching
to the continuous one for $p=1$. This has to be verified analytically
in the future. Even if they are indeed spectrally stable, it is not
clear if they are {\it nonlinearly} stable.

\end{enumerate}

It is worth mentioning some important questions we
{\it cannot} answer:
\begin{enumerate}

\item We are so far unable to give precise rigorous estimates on the
number and positions of the eigenvalues of $\L$ for $n \geq 2$
(question (i) above).

\item We cannot exclude the possible existence of embedded
eigenvalues (question (iii) above).

\item We do not know a nice variational formulation for
eigenvalues of $\L$ when $Q$ is an {\it excited} state
(this problem is also linked to question (i) above).

\item We do not have a complete characterization
of solitary waves, or more generally of nonlinear bound states.

\end{enumerate}

We end this introduction by describing
some related numerical work. Buslaev-Grikurov
\cite{BG, Grik} study the linearized operators for solitary waves of
the following 1D NLS with $p<q$,
\begin{equation*}
 i \psi_t + \psi_{xx} + |\psi|^{p}\psi - \al |\psi|^q \psi = 0.
\end{equation*}
They draw the bifurcation picture for eigenvalues near zero when the
parameter $\al>0$ is near a critical value, with the frequency of
the solitary wave fixed. This picture is similar to Weinstein's
picture which we study in Section~\ref{sec2-3}.

Demanet and Schlag \cite{Dem-Sch} consider the same linearization as
us and study the super-critical case $n=3$ and $p\le 3$ near $3$. In
this case, it is numerically shown that both $\Lp$ and $\Lm$ have no
eigenvalues in $(0,1]$ and no resonance at $1$, a condition which
implies (see \cite{Schlag04}) that $\L$ has no purely imaginary
eigenvalues in $[-i,0) \cup (0,i]$ and no resonance at $\pm i$.

We outline the rest of the paper: in Section~\ref{sec2} we consider
general results for all dimensions for ground states. In
Section~\ref{sec3} we consider one dimensional theory.  In
Section~\ref{sec4} we discuss the spectra for excited states with
angular momenta. In the Appendix we discuss the numerical methods.

{\it Notation}: For an operator $A$, $N(A)=\bket{\phi\in L^2| \ A
\phi=0}$ denotes the nullspace of $A$. $N_g(A)=\cup_{k=1}^\infty
N(A^k)$ denotes the generalized nullspace of $A$. The $L^2$-inner
product in $\R^n$ is $(f,g)=\int_{\R^n} \bar f g \, dx$.


\section{Revisiting the general theory for ground states}\label{sec2}

In this section we review mostly well-known results which are valid
for all dimensions, for the ground state $Q(x)=Q_{0,0,p}(x)$, and give
new proofs of some statements.

We begin by recalling some well-known results for
the linearized operator $\L$ defined by \eqref{eq1-9}.
As is well known for
linearized Hamiltonian system (and can be checked directly), if $\la$
is an eigenvalue, then so are $-\la$ and $\pm \bar \la$. Hence if
$\la\not = 0$ is real or purely imaginary, it comes in a pair. If it
is complex with nonzero real and imaginary parts, it comes in a
quadruple.  It follows from nonlinear stability and instability
results \cite{W86,GSS} that all eigenvalues are purely imaginary if $p
\in (1,p_c)$, and that there is at least one eigenvalue with positive
real part when $p \in (p_c,\pmax)$.  It is also known (see
e.g.~\cite{CPV}) that the set of isolated and embedded eigenvalues is
finite, and the dimensions of the corresponding generalized
eigenspaces are finite.

\subsection{$\Lp$, $\Lm$, and the generalized nullspace of
$\L$}\label{sec2-1}

Here we recall the makeup of the generalized nullspace
$N_g(\L)$ of $\L$. Easy computations give
\begin{equation}
\label{eq2-1}
  \Lp Q_1 = -2Q, \quad \Lm Q = 0, \quad \text{where }
Q_1 := (\tfrac 2{p-1}  + x \cdot \nabla )Q,
\end{equation}
and
\begin{equation}
\label{eq2-2}
  \Lm xQ = -2\nabla Q, \quad \Lp \nabla Q = 0.
\end{equation}
In the critical case $p = p_c$, we also have
\begin{equation}
\label{eq2-3}
  \Lm (|x|^2 Q) = -4 Q_1, \quad  \Lp \rho = |x|^2 Q
\end{equation}
for some radial function $\rho(x)$ (for which
we do not know an explicit formula in terms of $Q$).
Denote
\begin{equation} \label{pispc.def}
  \pispc = \left\{ \begin{array}{cc}
  1 &  p = p_c \\
  0 &  p \not= p_c.
  \end{array} \right.
\end{equation}
For $1 < p < p_{max}$, the generalized nullspace of $\L$ is
given by  (see \cite{W85})
\begin{equation}
\label{NgcalL}
  N_g(\L) = \mathrm{span} \bket{ \mat{0\\ Q}, \mat{0\\
  xQ}, \pispc \mat{0\\ |x|^2 Q},
  \mat{\nabla Q \\ 0}, \mat{Q_1\\ 0}, \pispc \mat{\rho \\ 0}}.
\end{equation}
In particular
\[
  \dim N_g(\L) = 2n + 2 + 2 \pispc.
\]
The fact that the vectors on the r.h.s of~\eqref{NgcalL} lie
in $N_g(\L)$ follows immediately from
the computations~\eqref{eq2-1}-\eqref{eq2-3}.
That these vectors span $N_g(\L)$ is established
in \cite{W85}, Theorems B.2 and B.3, which rely on
the non-degeneracy of the kernel of $\Lp$:
\begin{lemma}
\label{TH2-2A}
For all $n \ge 1$ and $p\in (1,\pmax)$,
\begin{equation*}
  N(\Lp) = \mathrm{span}\bket{\nabla Q}
\end{equation*}
\end{lemma}
This lemma is proved in \cite{W85} for certain $n$ and $p$
($n=1$ and $1<p<\infty$, or $n=3$ and $1<p \le 3$),
and is completely proved later by a general result of \cite{Kw}.
We present here a direct proof of this lemma, without referring to \cite{Kw},
relying only on oscillation properties of Sturm-Liouville ODE
eigenvalue problems.
A similar argument (which in the present case, however,
applies only for $p \leq 3$) appears in~\cite{FGJS}, Appendix C.
For completeness, we also
include some arguments of \cite{W85}.

\medskip

{\bf A new proof.} \quad
We begin with the cases $n \ge 2$.
Since the potential in $\Lp$ is radial, any solution of $\Lp v=0$ can
be decomposed as
$v = \sum_{k \ge 0} \sum_{{\bf j} \in \Sigma_k}
v_{k,{\bf j}}(r) Y_{k,{\bf j}}(\hat x)$, where $r=|x|$,
$\hat x = \frac x{r}$ is the spherical variable,
and $Y_{k,{\bf j}}$ denote spherical harmonics:
$- \Delta_{S^{n-1}} Y_{k,{\bf j}} = \la_k Y_{k,{\bf j}}$
(a secondary multi-index ${\bf j}$, appropriate to the dimension,
runs over a finite set $\Sigma_k$ for each $k$).
Then $\Lp v = 0$ can be written as
$A_k v_{k,{\bf j}} = 0$, where, for $k=0,1,2,3,\ldots$,
\[
  A_0 = - \pd _r^2 - \frac {n-1}r \pd_r + 1 - p Q^{p-1}(r), \quad
  A_k = A_0 + \la_k r^{-2} ,\quad \la_k = k(k+n-2).
\]

Case 1: $k=1$. Note $\nabla Q = Q'(r) \hat x$. Since $A_1 Q' = 0$ and
$Q'(r) < 0$ (monotonicity of the ground state) for $r \in (0,\infty)$,
$Q'(r)$ is the unique ground state of $A_1$ (up to a factor),
and so $A_1 \ge 0$, $A_1 |_{\{Q'\}^{\perp}} > 0$.

Case 2: $k \ge 2$. Since $A_k = A_1 + (\la_k - \la_1) r^{-2}$ and
$\la_k > \la_1$, we have $A_k >0$, and hence $A_k v_k=0$ has no
nonzero $L^2$-solution.

Case 3: $k=0$. Note that the first eigenvalue of $A_0$ is negative
because $(Q,A_0 Q) = (Q,-(p-1)Q^p) <0$. The second eigenvalue is
non-negative due to \eqref{TH2-1-eq2} and the minimax
principle. Hence, if there is a nonzero solution of $A_0 v_0=0$,
then $0$ is the second eigenvalue.  By Sturm-Liouville theory, $v_0(r)$
can be taken to have only one positive zero,
which we denote by $r_0>0$.
By \eqref{eq2-1}, $A_0 Q = -(p-1) Q^p$ and $A_0 Q_1 = -2Q$. Hence
$(Q^p,v_0)=0= (Q,v_0)$. Let $\al = (Q(r_0))^{p-1}$.  Since $Q'(r)<0$
for $r>0$, the function $Q^p-\al Q = Q(Q^{p-1}- \al)$ is positive for
$r<r_0$ and negative for $r>r_0$. Thus $v_0( Q^p-\al Q)$ does not
change sign, contradicting $(v_0, Q^p-\al Q) = 0$.
Combining all these cases gives Lemma \ref{TH2-2A} for $n \geq 2$.

Finally, consider $n=1$. Suppose $\Lp v = 0$. Since $\Lp$ preserves
oddness and evenness, we may assume $v$ is either odd or even.
If it is odd, it vanishes at the origin, and so by linear ODE
uniqueness, $v$ is a multiple of $Q'$.
So suppose $v$ is even. As in Case 3 above, since $\Lp$ has precisely
one negative eigenvalue, and has $Q'$ in its kernel,
$v(x)$ can be taken to have two zeros,
at $x = \pm x_0$, $x_0 \not= 0$. The argument of Case 3 above then
applies on $[0,\infty)$ to yield a contradiction.
{ \hfill $\square$}

\medskip

We complete this section by summarizing some positivity estimates
for the operators $\Lp$ and $\Lm$. These estimates are
closely related to the stability/instability of the ground state.

\begin{lemma} \label{TH2-1}
\begin{equation} \label{TH2-1-eq0}
  \Lm \geq 0, \quad \Lm|_{\bket{Q}^\perp} > 0 \quad (1 < p < \pmax)
\end{equation}
\begin{equation} \label{TH2-1-eq2}
  (Q, \Lp Q) < 0, \quad
  \Lp |_{\bket{Q^p}^\perp} \ge 0 \quad (1 < p < \pmax),
\end{equation}
\begin{equation} \label{TH2-1-eq1}
  \Lp |_{\bket{Q}^\perp} \ge 0 \quad (1< p \le p_c)
\end{equation}
\begin{equation} \label{TH2-1-eq3}
  \Lp|_{ \bket{Q,xQ}^\perp } > 0, \quad \Lm |_{\bket{ Q_1}^\perp} > 0 \quad
  (1< p < p_c)
\end{equation}
\begin{equation} \label{TH2-1-eq4}
  \Lp|_{\bket{Q,xQ,|x|^2Q}^\perp } > 0, \quad
  \Lm |_{\bket{ Q_1,\rho}^\perp } > 0 \quad (p = p_c).
\end{equation}
\end{lemma}
\begin{proof}
Most estimates here are proved in \cite{W85} except the second part of
\eqref{TH2-1-eq2} when $p >p_c$. It can be proved for $p
\in (1,\pmax)$ by modifying the proof of \cite[Prop.~2.7]{W85} for
\eqref{TH2-1-eq1} as follows.  (It is probably also well-known but we
do not know a reference.)
%

Recall the ground state $Q$ is obtained by the
minimization problem~\eqref{varprob}. If a minimizer $Q(x)$ is
rescaled so that
\[
  \frac {\int |\nabla Q|^2}{2a} = \frac {\int Q^2}{2b} =
  \frac {\int Q^{p+1}}{p+1} = \text{constant } k > 0,
\]
i.e., \eqref{pohozaev} is satisfied, then $Q(x)$ satisfies
\eqref{Q.eq}. The minimization inequality
$\frac {d^2}{d\e^2}\big|_{\e =0} J[Q+\e \eta] \ge 0$
for all real functions $\eta$, is equivalent to
\begin{equation} \label{eq2-1-eq5}
  k (\eta,\,\Lp\eta) \ge \frac 1a (\int \eta \Delta Q )^2 +
\frac 1b (\int Q \eta)^2 - (\int  Q^p \eta)^2 .
\end{equation}
Thus $(\eta,\,\Lp\eta) \ge 0$ if $\eta \perp Q^p$.  Note that, if
$\eta \perp Q$, by \eqref{Q.eq} the right side of \eqref{eq2-1-eq5}
is positive if $a\le 1$, i.e. $p \le p_c$. In this way, we recover
\eqref{TH2-1-eq1}.
\end{proof}

\subsection{Variational formulations of the eigenvalue problem for $\L$}
\label{var-Rn}

In this subsection we summarize various variational formulations for
eigenvalues of $\L$.  The generalized nullspace is given by
\eqref{NgcalL}.  Suppose $\la\not = 0$ is a (complex) eigenvalue of
$\L$ with corresponding eigenfunction $\vect{u \\ w} \in L^2$,
\begin{equation} \label{eq2-6}
   \mat{ 0 & \Lm \\ -\Lp & 0} \mat{u \\ w} = \la  \mat{u \\ w}.
\end{equation}
The functions $u$ and $w$ satisfy
\begin{equation}\label{eq2-7A}
  \Lp u =  -\la w, \quad \Lm w =  \la u.
\end{equation}
Therefore
\begin{equation} \label{eq2-8}
 \Lm \Lp u = \mu u,  \quad \mu = -\la^2.
\end{equation}
Since $(\mu u, Q)= (\Lm \Lp u,Q)=(\Lp u,\Lm Q)=0$ and $\mu \not =0$,
we have $u \perp Q$.

Denote by $\Pi$ the $L^2$-orthogonal projection onto $Q^\perp$.  We
can write $\Lp u = \Pi \Lp u + \al Q$.  Eq.~\eqref{eq2-8} implies $\Lm
\Pi \Lp u =\mu u$ and hence, using $u \perp Q$
and~\eqref{TH2-1-eq0}, $\Pi \Lp u =\Lm^{-1} \mu u$. Thus
\begin{equation} \label{eq2-9}
(u,Q)=0, \quad \Lp u = \mu \Lm^{-1} u + \al Q.
\end{equation}
Since \eqref{eq2-8} is also implied by Eq.~\eqref{eq2-9}, these two
equations are equivalent.

If $Q(x)$ is a general solution of \eqref{Q.eq}, $\mu = -\lambda^2$ may not be
real. However, it must be real for the nonlinear ground state $Q=Q_{0,0,p}$.
This fact is already known (see \cite{RSSarxiv}).
We will give a different proof.

\medskip

{\bf Lemma.}
For $Q=Q_{0,0,p}$, every eigenvalue $\mu$ of \eqref{eq2-8} is real.

\medskip

{\bf A new proof.}\quad
Multiply \eqref{eq2-7A} by $\bar{u}$ and $\bar{w}$ respectively and
integrate. Then we get
\begin{equation} \label{eq2-11}
 (u, \Lp u) = -\lambda (u,w), \quad
 (w, \Lm w) = \lambda (w,u) = \lambda \overline{(u,w)}.
\end{equation}
Taking the product, we get
\begin{equation*}
 (u,\Lp u)(w,\Lm w) = -\lambda^2 |(u,w)|^2 = \mu |(u,w)|^2.
\end{equation*}
If $\mu \not= 0$, $w$ is not a multiple of $Q$, and
so by~\eqref{TH2-1-eq0},  $(w,\Lm w) > 0$. Hence
$(u,w)\not = 0$ by \eqref{eq2-11}.
Thus
\begin{equation*}
  \mu = \frac {(u,\Lp u)(w,\Lm w)}{ |(u,w)|^2} \in \R.
\end{equation*}
{ \hfill $\square$}

This argument does not work when $Q$ is an excited state, since
$(u,w)$ may be zero (see e.g.~\cite[Eq.(2.63)]{TY3}).  The fact $\mu
\in \R$ implies that eigenvalues $\la$ of $\L$ are either real or
purely imaginary.  Thus $\L$ has no complex eigenvalues with nonzero
real and imaginary parts. This is not the case for excited states
(see Section~\ref{sec4}, also~\cite{TY3}).

The proof of reality of $\mu$ in \cite{RSSarxiv}
uses the following formulation.  For the
nonlinear ground state $Q$, $\Lm$ is nonnegative and the operator
$\Lm^{1/2}$ is defined on $L^2$ and invertible on $Q^\perp$.  A
nonzero $\mu\in \C$ is an eigenvalue of \eqref{eq2-8} if and only if
it is also an eigenvalue of the following problem:
\begin{equation} \label{eq2-10}
\Lm^{1/2} \Lp \Lm^{1/2} g = \mu g,
\end{equation}
with $ g= \Lm^{-1/2} u$. The operator $\Lm^{1/2} \Lp \Lm^{1/2}$
already appeared in \cite{TY1}. Since it can be realized as a
self-adjoint operator, $\mu$ must be real.

Furthermore, the eigenvalues of $\Lm^{1/2} \Lp \Lm^{1/2}$ can be
counted using the minimax principle.  Note that $Q$ is an eigenfunction
with eigenvalue $0$. For easy comparison with other formulations, we
formulate the principle on $Q^\perp$. Let
\begin{equation}  \label{varA}
\mu_j := \inf_{ g \perp Q,
g_k, k=1,\ldots,j-1} \frac{(g,\Lm^{1/2}
\Lp \Lm^{1/2} g)}{ (g,g)}, \qquad (j=1,2,3,\ldots)
\end{equation}
with a suitably normalized minimizer denoted by $g_j$
(if it exists -- the definition terminates once a minimizer
fails to exist).
The corresponding definition for \eqref{eq2-9} is
\begin{equation} \label{varB}
\mu_j := \inf _{u \perp Q,\, (u ,\, \Lm^{-1} u_k)=0,\,k=1,\ldots ,j-1}
\, \frac{(u,L_+ u)}{ (u,\Lm^{-1} u)} , \qquad
(j=1,2,3,\ldots)
\end{equation}
with a suitably normalized minimizer denoted by $u_j$ (if it
exists). In fact, the minimizer $u_j$ satisfies
\begin{equation} \label{eq2-19}
\Lp u_{j} =  \mu_j \Lm^{-1}u_{j} + \al_j Q + \beta_1 \Lm^{-1} u_1 +
\cdots + \beta_{j-1} \Lm^{-1} u_{j-1},
\end{equation}
for some Lagrange multipliers $ \beta_1,\ldots \beta_{j-1}$.
Testing~\eqref{eq2-19} with $u_k$ with $k<j$, we get $ (u_k, \beta_k
\Lm^{-1}u_k)= (u_k,\Lp u_j) = (\Lp u_k , u_j)=0$ by \eqref{eq2-19}
for $u_k$ and the orthogonality conditions.  Thus $\beta_k=0$ and $
\Lp u_{j} =  \mu_j \Lm^{-1}u_{j} + \al_j Q $ and hence $u_j$
satisfies \eqref{eq2-9}.


\begin{lemma} \label{TH2-3}
The eigenvalues of \eqref{varA} and \eqref{varB} are the same, and
\begin{align*}
\text{if } 1<p<p_c  &: \quad \mu_1 = \cdots = \mu_n = 0,  \quad \mu_{n+1}>0.
\\
\text{if } p=p_c  &: \quad \mu_1 = \cdots = \mu_{n+1} = 0,  \quad \mu_{n+2}>0.
\\
\text{if } p_c<p<\pmax  &: \quad \mu_1 <0, \quad
\mu_2= \cdots = \mu_{n+1} = 0,  \quad \mu_{n+2}>0.
\end{align*}
The $0$-eigenspaces are span\,$\Lm^{-1/2}\{\nabla Q,\pispc Q_1\}$
for \eqref{varA} and span$\{\nabla Q,\pispc Q_1\}$ for \eqref{varB},
where $\pispc$ is defined in \eqref{pispc.def}.
\end{lemma}

\begin{proof}
The eigenvalues of \eqref{varA} and \eqref{varB} are seen to be the
same by taking $g = \Lm^{-1/2}u$ up to a factor.  By
estimate~\eqref{TH2-1-eq1}, $\mu_1 \ge 0$ for $p \in (1,p_c]$. For
$p \in (p_c,\pmax)$, using \eqref{pohozaev}, $\Pi Q_1 = Q_1 - \frac
{(Q_1,Q)}{(Q,Q)}Q$, and elementary computations (such as~\eqref{Q1Q}
below), one finds
\[
  (\Pi Q_1, \Lp \Pi Q_1) = \frac {n^2 (p-1)}{4}(p_c -p)
  \frac{1}{p+1} \int Q^{p+1}
\]
which is negative for $p>p_c$. Thus $\mu_1 <0$. By estimate
\eqref{TH2-1-eq2}, $\mu_2 \ge 0$ for $p \in (1,\pmax)$.

It is clear that $u = \frac {\pd }{\pd x_j} Q$, $j=1,\ldots,n$,
provides $n$ $0$-eigenfunctions. For $p=p_c$, another
$0$-eigenfunction is $u = Q_1$ since $Q_1 \perp Q$
(see again~\eqref{Q1Q} below), $\Lm^{-1}\nabla Q = - \frac 12 xQ$,
and $(Q_1,\Lp Q_1)=0$.
It remains to show that
$\mu_{n+1}>0$ for $p \in (1,p_c)$ and $\mu_{n+2}>0$ for $p \in
[p_c,\pmax)$.  If $\mu_{n+2}=0$ for $p \in (p_c,\pmax)$, the argument
after \eqref{eq2-19} shows the existence of a function $u_{n+2} \not =
0$ satisfying
\[
\Lp u_{n+2} = \al Q \quad \text{for some } \al\in \R, \quad u_{n+2}
\perp Q,\, \Lm^{-1} u_1,\, \Lm^{-1} \nabla Q= - \tfrac 12 xQ.
\]
By Lemma~\ref{TH2-2A}, $u_{n+2} + \frac \al 2 Q_1 = c \cdot\nabla Q$
for some $c\in \R^d$.  The orthogonality conditions imply
$u_{n+2}=0$. The cases $p \in (1,p_c]$ are proved similarly.
\end{proof}

\begin{rem}

The formulation~\eqref{varB} for $\mu_1$ has been used for the
stability problem, see e.g. \cite[p.73, (4.1.9)]{SS}, which can be
used to prove that $\mu_1<0$ if and only if $p \in (p_c,\pmax)$ by a
different argument. The later fact also follows from \cite{W85,GSS}
indirectly.

\end{rem}

We summarize our previous discussion in the following theorem.

\begin{summary}  \label{TH2-5A}
Let $Q(x)$ be the unique positive radial ground state solution of
\eqref{Q.eq}, and let $\L$, $\Lp$ and $\Lm$ be as in \eqref{eq1-9} and
\eqref{eq1-11}.  The eigenvalue problems~\eqref{eq2-8}, \eqref{eq2-9},
and \eqref{eq2-10} for $\mu \not =0$ are equivalent, and the
eigenvalues $\mu$ must be real. These eigenvalues can be counted by
either \eqref{varA} or \eqref{varB}.  $\mu_1 < 0$ if and only if $p
\in (p_c,\pmax)$. Furthermore, all eigenvalues of $\L$ are purely
imaginary except for an additional real pair when $p \in (p_c,
\pmax)$.
\end{summary}

The last statement follows from the relation $\mu = - \la^2$ in \eqref{eq2-8}.

\subsection{Spectrum near $0$ for $p$ near $p_c$} \label{sec2-3}

We now consider eigenvalues of $\L$ near $0$ when $p$ is near $p_c$.
It was suggested by M.I. Weinstein that as $p$ approaches $p_c$ from
below, a pair of purely imaginary eigenvalues will collide at the
origin, and split into a pair of real eigenvalues for $p >
p_c$.  In the following theorem and corollary we prove this
picture rigorously and identify the leading terms of the eigenvalues
and eigenfunctions.

Note that Comech-Pelinovsky \cite{ComPel} considers a different problem
where the equation is fixed and the varying parameter is frequence
$\omega$ rather than exponent $p$ of the nonlinearity.  That problem
has only $U(1)$ symmetry and no translation, but its situation is
similar to ours since we consider radial functions only in our
proof. It seems one can adapt their approach to give an alternative
proof.  They use an abstract projection (Riesz projection) onto the
discrete spectrum to reduce the problem to a 4x4 matrix problem (and
exploit the complex structure), while we are more direct. We thank
the referee for pointing out \cite{ComPel} to us.

\begin{theorem} \label{th2-2}
There are small constants $\mu_*>0$ and $\e_* >0$ so that for every
$p \in (p_c - \e_*, p_c + \e_*)$, there is a solution of
\begin{equation} \label{eq2-14}
  \Lp \Lm w = \mu w
\end{equation}
of the form
\begin{equation*}
  w = w_0 + (p-p_c)^2 g, \quad
  w_0 = Q + a(p-p_c) |x|^2 Q, \quad
  g \perp Q,
\end{equation*}
\begin{equation*}
  \mu = 8a(p-p_c) + (p-p_c)^2 \eta,
  \quad
  a = a(p) = \frac{n (Q_1, Q^p)}{4(Q_1, x^2Q)}<0,
\end{equation*}
with $\norm{g}_{L^2}$, $|\eta|$, $|a|$ and $1/|a|$ uniformly bounded
in $p$. Moreover, for $p \not= p_c$,
this is the unique solution of \eqref{eq2-14} with
$0<|\mu|\le \mu_*$.
\end{theorem}

\begin{proof}
Set $\e := p-p_c$.
Computations yield
\begin{equation} \label{Q1Q}
  (Q_1, Q )
  = \bke{ \frac 2{p-1} - \frac n2} (Q, Q )
  = -\frac{\e n}{2(p-1)} (Q, Q ),
\end{equation}
\begin{equation} \label{Q1Qp}
  (Q_1,Q^p  )= -\frac{1}{p-1} (\Lp Q, Q_1 )
  = -\frac{1}{p-1} (Q, \Lp Q_1 )
  = \frac{2}{p-1} (Q, Q ),
\end{equation}
and
\begin{equation} \label{Q1x2Q}
  (Q_1, |x|^2 Q )= \left( \frac{2}{p-1} - \frac{n+2}{2} \right)
  (Q,  |x|^2 Q )
  = -( 1 + \frac{\e n}{2(p-1)}) (Q,  |x|^2 Q ).
\end{equation}
Since by~\eqref{eq2-14} with $\mu \not= 0$,
\[
  (Q_1,w) = \mu^{-1}(Q_1, \Lp \Lm w) = \mu^{-1}(\Lm \Lp Q_1,  w)  = 0,
\]
we require the leading term $(Q_1,w_0)=0$, which decides the value of
$a$ using \eqref{Q1Q} and \eqref{Q1x2Q}. Thus we also need
$(Q_1,g)=0$. That $a<0$ (at least for $\e$ sufficiently small) follows
from \eqref{Q1Qp} and \eqref{Q1x2Q}.

Using the computations
\begin{equation}\label{L-x2Q}
  \Lm  |x|^2 Q = [\Lm, |x|^2]Q
  = -4 x \cdot \nabla Q
  - 2nQ = -4Q_1 - \frac {2n}{p-1} \e Q
\end{equation}
and
\begin{equation*}
  \Lp Q = [\Lm - (p-1)Q^{p-1}] Q = -(p-1)Q^p,
\end{equation*}
we find
\begin{equation*}
  \Lp \Lm w_0 = a\e\Lp[-4Q_1-\frac {2n}{p-1} \e Q]
  = a \e[8Q + 2n\e Q^p].
\end{equation*}
Thus $\mu = 8a\e + o(\e)$ and we need to solve
\begin{equation*}
  0 = [\Lp \Lm - 8a\e - \e^2 \eta][w_0+\e^2 g]
\end{equation*}
which yields our main equation
for $g$ and $\eta$:
\begin{equation}\label{eq:main}
  \Lp \Lm g =
  8a^2( |x|^2Q) - 2an(Q^p) + \eta w_0
  +(8a\e+\e^2\eta)g.
\end{equation}
Recall that on radial functions
(we will only work on radial functions here)
\begin{equation*}
  \ker[(\Lp \Lm)^*] = \ker[ \Lm \Lp ] = {\rm span} \{Q_1 \}.
\end{equation*}
Let $P$ denote the $L^2$-orthogonal projection
onto $Q_1$ and $\bP := {\bf 1} - P$.
It is necessary that
\begin{equation*}
  P [8a^2( |x|^2Q) - 2an(Q^p) + \eta w_0
  +(8a\e+\e^2\eta)g] = 0
\end{equation*}
for \eqref{eq:main} to be solvable.
This solvability condition holds
since $(Q_1,g) = (Q_1,w_0) = 0$, and,
using the relations \eqref{Q1x2Q} and
\eqref{Q1Qp}, $ (Q_1, 8a^2( |x|^2Q) - 2an(Q^p) )= 0$.

Consider the restriction (on radial functions)
\[
T = \Lp \Lm :\ [\ker \Lm]^{\perp} = Q^{\perp}
\longrightarrow Ran(\bP) = Q_1^{\perp}.
\]
Its inverse
$T^{-1}=(\Lm)^{-1} (\Lp)^{-1}$ is bounded because
$(\Lp)^{-1} : Q_1^\perp \ra Q^\perp$ and
$(\Lm)^{-1} : Q^\perp \ra Q^\perp$
are bounded. So our strategy is to solve~(\ref{eq:main}) as
\begin{equation}\label{eq:main1a}
  g = T^{-1} \bP [8a^2( |x|^2Q) - 2an(Q^p) + \eta w_0
  +(8a\e+\e^2\eta)g]
\end{equation}
by a contraction mapping argument, with $\eta$
chosen so that $ (Q_1,g)=0$. Specifically, we
define a sequence $g_0=0$, $\eta_0=0$, and
\begin{align*}
g_{k+1} &= \bP T^{-1} \bP [8a^2( |x|^2Q) - 2an(Q^p) + \eta_k w_0
  +(8a\e+\e^2\eta_k)g_k],
\\
\eta_{k+1} &= - \frac{1}{(Q_1,T^{-1} w_0)}\, (Q_1,T^{-1} \bP
  [8a^2( |x|^2Q) - 2an(Q^p)  +(8a\e+\e^2\eta_k)g_k]).
\end{align*}
We need to check $(Q_1,T^{-1} w_0)$ is of order one.  Since $w_0 = Q +
O(\e)$ and $\Lp Q_1 = -2Q$, we have $(\Lp)^{-1} w_0 = -\frac 12 \Pi Q_1
+ O(\e)$ where $\Pi$ denotes the orthogonal projection onto $Q^\perp$. Thus,
using \eqref{L-x2Q} and \eqref{Q1Q},
\begin{align*}
 (Q_1,T^{-1} w_0) &=  -\frac{1}{2}(Q_1, (\Lm)^{-1} \Pi Q_1 )+ O(\e)
\\
& = \frac{1}{8} (Q_1,  \Pi |x|^2 Q )+ O(\e)
 =\frac{1}{8} (Q_1, |x|^2 Q )+ O(\e) ,
\end{align*}
which is of order one because of \eqref{Q1x2Q}.
One may then check that $N_{k} :=
\norm{g_{k+1}-g_{k}}_{L^2} + \e^{1/2}
|\eta_{k+1}-\eta_{k}|$ satisfies $N_{k+1} \le C
\e^{1/2} N_k$, and hence $(g_k,\eta_k)$ is indeed
a Cauchy sequence.

Finally, the uniqueness follows from the invariance of the total dimension of
generalized eigenspaces near $0$ under perturbations.
\end{proof}

\begin{rem}
To understand heuristically the leading terms in $w$ and $\mu$,
consider the following analogy. Let $A_\e = \mat{ 0 & 1 \\ 0 & \e}$,
which corresponds to $\Lp \Lm$. One has $A_\e \vect{1\\ 0}=\vect{ 0 \\
0}$, $A_\e \vect{0 \\ 1} = \vect{ 1 \\ \e}$ and $A_\e \vect{1 \\ \e} =
\e \vect{1 \\ \e}$. The vectors $\vect{1\\ 0}$, $\vect{0 \\ 1}$ and
$\vect{1 \\ \e}$ correspond to $Q$, $|x|^2Q$ and $w$, respectively.
\end{rem}

The theorem yields an eigenvalue $\mu$ with the same sign as $p_c-p$.
Since the eigenvalues of $\L$ are given by $\la = \pm \sqrt {-\mu}$,
we have the following corollary.
\begin{corollary}
With notations as in Theorem~\ref{th2-2}, $\L$ has a pair of
eigenvalues $\la = \pm \sqrt {-\mu} = \pm \sqrt{8|a|(p-p_c) -
(p-p_c)^2 \eta}$ with corresponding eigenvectors $\vect{u \\ w}$
solving \eqref{eq2-6} and
\[
u = \la^{-1} L_- w = \mp \sqrt {2 |a| (p-p_c)} \,Q_1 +
O((p-p_c)^{3/2}).
\]

When $p \in (p_c - \e_*,\ p_c)$
(stable case), $\la$ and $u$ are purely imaginary.

When $p \in (p_c ,
\ p_c + \e_*) $ (unstable case), $\la$ and $u$ are real.
\end{corollary}

In deriving the leading term of $u$ we have used \eqref{L-x2Q}.  We
solved for $w$ before $u$ simply because $w$ is larger than $u$.


\newcommand{\p}{\partial}
\newcommand{\ti}{\tilde}
\newcommand{\Del}[1]{}
\newcommand{\EQ}[1]{\begin{equation}#1\end{equation}}
\newcommand{\EQAL}[1]{\begin{equation}\begin{split} #1 \end{split}\end{equation}}
\newcommand{\CA}[1]{\begin{equation}\begin{cases} #1 \end{cases}\end{equation}}
\newcommand{\Z}{\mathbb{Z}}

\section{One dimensional theory}\label{sec3}
In this section we focus on the one dimensional theory. For $n=1$,
the ground state $Q(x)$ has an explicit formula for all $p\in (1,\infty)$,
\begin{equation} \label{eq3-1}
Q(x) = c_p \cosh^{-\be} (x /\be),
\quad  c_p := (\frac{p+1}2)^{\frac 1{p-1}} ,
\quad \be := \frac 2{p-1}.
\end{equation}
The function $Q(x)$ satisfies \eqref{Q.eq} and is the unique
$H^1(\R)$-solution of \eqref{Q.eq} up to translation and phase
\cite[p.259, Theorem 8.1.6]{Ca}.

\subsection{Eigenfunctions of $\Lp$ and $\Lm$}\label{sec3-1}
We first consider eigenvalues and eigenfunctions of $\Lp$ and $\Lm$.
For $n=1$,
\begin{equation} \label{1DLpLm.def}
\Lp = - \pd_{xx}+1- p Q^{p-1}, \qquad \Lm = - \pd_{xx}+1-  Q^{p-1}.
\end{equation}
By \eqref{eq3-1}, these operators are both of the form
\begin{equation*}
  -\pd_{xx}+1- C \sech^2(x/\beta).
\end{equation*}
Such operators have essential spectrum $[1,\infty)$, and
finitely many eigenvalues below $1$.
A lot of information about such operators is available in the
classical book \cite{Ti}, p. 103:
\begin{itemize}
\item all eigenvalues are simple, and can be computed explicitly, as
zeros and poles of an explicit meromorphic function;
\item all eigenfunctions can be expressed in terms of the
hypergeometric function.
\end{itemize}

We begin by presenting another way to derive the eigenvalues, as well as
different formulas for the eigenfunctions. We will not prove
right here that this set contains all of the
eigenvalues/eigenfunctions. This fact is a consequence of
the more general Theorem~\ref{thm:hier}, proved later
(and see also~\cite{Ti}).

Define
\begin{equation} \label{lam.def}
\begin{split}
&\la_m := 1- k_m^2, \quad\quad
k_m := \frac {p+1}2 - \frac {m(p-1)}2, \\
&p_m := \frac{m+1}{m-1} \;\; \text{ for } m>1, \;\;\ p_1 = \infty.
\end{split}
\end{equation}
The following theorem agrees with the numerical
observation Figure~\ref{fig:cL1d}.

\begin{theorem} \label{th3-1}
For $n=1$ and $1<p<\infty$, let $Q(x)$ be defined by \eqref{eq3-1},
$\Lp$ and $\Lm$ be defined by \eqref{1DLpLm.def}, and $\la_m,k_m,p_m$
be defined by \eqref{lam.def}.  Suppose for $M \in \Z^+$,
\begin{equation} \label{eq3-4}
  p_{M+1} \le p < p_{M}.
\end{equation}
Then the operator $\Lp$ has eigenvalues
$\la_m$, $ 0 \le m \le M$,
with eigenfunctions of the form
\[
\ph_{2\ell} = \sum_{j=0}^{\ell} c^{2\ell}_{2j} Q^{k_{2j}}, \quad
\ph_{2\ell-1} = \sum_{j=1}^\ell c^{2\ell-1}_{2j-1} (Q^{k_{2j-1}})_x ,
\]
and the operator $\Lm$ has eigenvalues $\la_m$, $ 1\le m\le M$, with
eigenfunctions of the form
\[
\psi_{2\ell-1}= \sum_{j=1}^\ell d^{2\ell-1}_{2j-1} Q^{k_{2j-1}} ,
\qquad
\psi_{2\ell}= \sum_{j=1}^{\ell} d^{2\ell}_{2j} (Q^{k_{2j}})_x.
\]
In particular, all eigenvalues of $\Lm$ are eigenvalues of $\Lp$, and
$\Lp$ always has one more eigenvalue ($\la_0<0$) than $\Lm$.  
\end{theorem}

\begin{proof} It can be proved by induction, using
\[
Q^{p-1} = \frac {p+1}2 \, \cosh^{-2} (x /\be), \quad
Q_x = - Q \tanh (x/\be),\quad
Q_x ^2 = Q^2 (1- \frac 2{p+1}Q^{p-1}),
\]
and
\begin{align}
Q^{-k} L_+ Q^k  &= \frac {(k+p)\,(2k-p-1)} {p+1} \,  Q^{p-1} +
(1 - k^2), \label{eq3-5}
\\
[(Q^k)_x]^{-1} L_+ (Q^k)_x &= \frac{(k-1)(2k+3p-1)}{p+1} \,Q^{p-1} +
(1 - k^2),
\nonumber
\\
Q^{-k} L_- Q^k &= \frac { (k-1)(2k+ p+1)}{p+1}\,  Q^{p-1} + (1- k^2), 
\label{eq3-8}
\\
[(Q^k)_x]^{-1} L_- (Q^k)_x &= \frac {(k+p)(2k+ p-3)}{p+1}\,  Q^{p-1} +
(1 - k^2).
\nonumber
\end{align}
The coefficients of $Q^{p-1}$ vanish when $k=\frac {p+1}2, 1, 1,
\frac{3-p}2$, respectively. It is why the highest power of $Q$
is $Q^{\frac {p+1}2}$ in $\ph_{2\ell}$, $Q_x$ in $\phi_{2\ell-1}$, $Q$ in
$\psi_{2\ell-1}$, and $(Q^{\frac {3-p}2})_x$ in $\psi_{2\ell}$.
\end{proof}

\subsection{Connection between $L_+$ and $\Lm$ and their factorizations}
\label{sec3-2}

In light of Theorem \ref{th3-1}, it is natural to ask {\it why} all
eigenvalues of $\Lm$ are also eigenvalues of $\Lp$. Is there a simple
connection between their eigenfunctions?  In this section we prove
this is indeed so.

We first look for an operator $U$ of the form
\begin{equation*}
U = \pd_x + R(x), \quad\quad (\text{ so } U^* = -\pd_x + R(x)),
\end{equation*}
such that
\begin{equation} \label{eq3-9}
 L_- U  = U  L_+ , \quad\quad (\text{ so } U^* L_- = L_+ U^*).
\end{equation}
It turns out that there is a unique choice of $R(x)$:
\begin{equation*}
R(x) = - \frac {p+1}2 \frac {Q_x}Q =  \frac {p+1}2
\tanh(\frac {(p-1)x}{2}).
\end{equation*}
In fact, with this choice of $R(x)$,
\begin{equation}\label{U.def}
U =\ph_0  \pd_x \ph_0^{-1}, \quad\quad
(\text{ so } U^* = -\ph_0 ^{-1} \pd_x \ph_0),
\end{equation}
where $\ph_0 = Q^{\frac{p+1}2}$ is the ground state of $L_+$, and is
considered here as a multiplication operator:
$U f= \ph_0 \pd_x (\ph_0 ^{-1} f)$.

Suppose now $\psi$ is an eigenfunction of $L_-$ with eigenvalue $\la$:
$L_- \psi = \la \psi$. By \eqref{eq3-9},
\begin{equation*}
0 = U^* (L_- - \la) \psi = (L_+ - \la) U^*  \psi.
\end{equation*}
Thus $U^* \psi$ is an eigenfunction of $L_+$ with same eigenvalue
$\la$ (provided $U^* \psi \in L^2$). Therefore, the map
\begin{equation*}
\psi \mapsto U^*  \psi
\end{equation*}
sends an eigenfunction of $L_- $ to an eigenfunction of $L_+$ with
same eigenvalue. This map is not onto because $U^* $ is not
invertible.  Specifically, the ground state $\ph_0$ is not in the
range. In fact, $U \ph _0 =\ph_0 \pd_x \ph_0 ^{-1}\ph_0 =0$. If $\ph_0
= U^* \psi$, then $(\ph_0,\ph_0)= (\ph_0,U^* \psi) = (U \ph_0 ,
\psi)=0$, a contradiction.  We summarize our finding as the
following proposition.

\begin{proposition}
Under the same assumptions and notation as Theorem~\ref{th3-1}, the
eigenfunctions $\ph_m$ and $\psi_m$ of $\Lp$ and $\Lm$ satisfy
\begin{equation*}
\ph_m = U^* \psi_m, \qquad (m=1,\ldots,M),
\end{equation*}
up to constant factors.
Note that $U^*$ sends even functions to odd functions and vice versa.
\end{proposition}

\begin{proof}
We only need to verify that $U^* \psi_m \in L^2$. This is the case since
$U^* =- \pd_x + \frac {p+1}2 \tanh(x/\beta)$, $\psi_m(x)$ are sums of
powers of $Q$ and $Q_x$, and that $\tanh(x/\beta)$, $Q_x/Q$, and
$Q_{xx}/Q_x$ are bounded.
\end{proof}

Analogous to the definition of $U$, we define
\begin{equation} \label{S.def}
S := Q \pd_x Q^{-1} = \pd_x - \frac {Q_x}Q, \quad\quad
(\text{ so } S^* = -Q^{-1} \pd_x Q).
\end{equation}
Clearly $SQ =0$. Recall that $\la_0$ is the first eigenvalue of $\Lp$
with eigenfunction $\ph_0$. Hence $L_+-\la_0$ is a nonnegative
operator. In fact we have the following factorizations.

\begin{lemma} \label{th3-3}
Let $U$ and $S$ be defined by \eqref{U.def} and \eqref{S.def},
respectively. One has
\begin{equation} \label{eq3-44}
 L_+ - \la_0 = U^* U , \qquad L_- - \la_0 = U  U^* .
\end{equation}
\begin{equation}\label{eq3-14}
\Lm = S^* S, \qquad SS^* = -\pd_x ^2 + 1 + \frac {p-3}{p+1} Q^{p-1}.
\end{equation}
Moreover, $SS^*  = \Lm + \frac {2(p-1)}{p+1} Q^{p-1}>0$.
\end{lemma}

The formula $\Lm = S^* S$ was known, see e.g.~\cite[p.73,
(4.1.8)]{SS}. It is an example of the Darboux
transformations, see e.g.~\cite{Matveev-Salle}.  Factorization of
Schr\"odinger operators into first-order operators has been known
since the times of Darboux (1840s).

\subsection{Hierarchy of Operators}\label{sec3-3}
In this subsection we generalize Theorem~\ref{th3-1} and
Lemma~\ref{th3-3} to a family of operators containing $\Lp$ and
$\Lm$. As a reminder, we have \EQAL{
 &Q''/Q = 1 - Q^{p-1},\quad (Q'/Q)^2 = 1-\frac{2}{p+1}Q^{p-1},\\
 &(Q'/Q)'=Q''/Q-(Q'/Q)^2 = -\frac{p-1}{p+1}Q^{p-1}.}
Let $S(a):=Q^a\p_x Q^{-a}$. We have
\EQAL{
 &S(a)=\p_x - aQ'/Q,\ S(a)^*=-\p_x -aQ'/Q,\\
 &S(a)^* S(a)=-\p_x^2 + a^2 - a\left\{a+\frac{p-1}{2}\right\}\frac{2}{p+1} Q^{p-1}.}
Define the following hierarchy of operators:
\EQAL{ \label{hiera}
 &S_j := S(k_j),\quad \text{ where recall }
 \quad k_j = 1 - (j-1) \frac{p-1}{2},\\
 &L_j := S_{j-1} S_{j-1}^* + \la_{j-1} = S_j^* S_j + \la_j, \quad
 \text{ where recall } \quad \la_j=1-k_j^2.}
Then we have
\EQAL{ \label{hiera2}
 &S_0 = U,\ S_1 = S, \dots \\
 &L_0 = L_+,\ L_1 = L_-,\ L_2 = SS^*,\dots \\
 &S_j L_j = L_{j+1} S_j, \quad L_j S_j^* = S_j^* L_{j+1}.}
More explicitly,
\EQAL{ \label{Lj}
  L_j = -\p_x^2 + 1 - k_{j-1} k_j \frac{2}{p+1} Q^{p-1}.}
Note that $j$ here can be any real number.

Recall the definition $p_j := 1+2/(j-1)$ for $j>1$, and
set $p_j = \infty$ for $j \leq 1$.
Then $p_j$ is a monotone decreasing function of $j$,
$k_j>0$ for $p<p_j$, $k_j=0$ for $p=p_j$ and
$k_j<0$ for $p>p_j$. Let \EQ{
 \la_j' :=
 \begin{cases}
  \la_j &(1 < p \leq p_j),\\
  1 &(p_j<p<p_{j-1}),\\
  \la_{j-1} &(p_{j-1} \leq p).
 \end{cases}}
By the second identity of \eqref{hiera}, and
\eqref{Lj} together with the fact $k_{j-1}k_j<0$ for
$p_j<p<p_{j-1}$, we have the lower bound
\EQ{ \label{low bd}
 L_j \ge \la_j'.}
In fact, this estimate is sharp: for $p \in (1,p_j) \cup (p_{j-1},\infty)$,
the ground state is obvious from the second identity of \eqref{hiera}:
\CA{
 L_j Q_j = \la_j Q_j, &(1<p<p_j),\\
 L_j Q_{j-1}^* = \la_{j-1} Q_{j-1}^*, &(p_{j-1}<p),}
where we denote
\EQ{
 Q_j := Q^{k_j}, \quad Q_j^* := Q^{-k_j}.}
For $p \in [p_j, p_{j-1}]$, there is no ground state.
Thus we have completely determined the ground
state of $L_j$ for all $p>1$. The complete
spectrum, together with explicit eigenfunctions,
are derived using the third identity of \eqref{hiera2} as follows.
\begin{theorem}
\label{thm:hier}
For any $j\in\R$ and $p>1$, the point spectrum of
$L_j$ consists of simple eigenvalues \EQAL{
\label{pspec}
 spec_{p}(L_j) = &\{\la_k \mid p<p_k,\ k \in \{j,j+1,j+2,\dots\}  \}\\
 & \cup \{\la_k \mid p>p_k,\ k \in \{ j-1,j-2,j-3,\dots \}  \},}
and the eigenfunction for the eigenvalue $\la_k$
is given uniquely up to constant multiple by
\CA{ \label{egf}
 S_j^* \cdots S_{k-1}^* Q_k & (k \in \{j,j+1,\dots\} ),\\
 S_{j-1} \cdots S_{k+1} Q_k^* & (k \in \{j-1,j-2,\dots\} ),
}
each of which is a linear combination of
\CA{ \label{lc}
 Q_j, Q_{j+2},\dots Q_k & (k \in \{j,j+2,\dots\}),\\
 Q_{j+1}R, Q_{j+3}R,\dots Q_k R & (k \in \{j+1,j+3,\dots\}),\\
 Q_{j-1}^*, Q_{j-3}^*,\dots Q_k^* & (k \in \{j-1,j-3,\dots\}),\\
 Q_{j-2}^*R, Q_{j-4}^*R,\dots Q_k^* R & (k \in \{j-2,j-4,\dots\})}
where $R := Q'/Q$.
\end{theorem}
\begin{proof}
The ground states have been determined. The third
identity of \eqref{hiera2} implies that
\eqref{egf} belong to the eigenspace of $L_j$
with eigenvalue $\la_k$. Moreover, each function is nonzero
because $S_k^*$ is injective for $p<p_k$ and so
is $S_k$ for $p_k<p$. Since $S_j$ annihilates only the
ground state $Q_j$ for $p<p_j$ and $S_{j-1}^*$
annihilates only the ground state $Q_{j-1}^*$ for
$p>p_j$, all the excited states of $L_j$ for
$p<p_j$ are mapped injectively by $S_j$ to
bound states of $L_{j+1}$, and for $p>p_j$ by
$S_{j-1}^*$ to those of $L_{j-1}$. Hence we have
\eqref{pspec} and all the eigenvalues are simple
because the ground states are so. \eqref{lc}
follows from the fact that $S_j$ and $S_j^*$ act
on $Q^a$ like $C(a,j)R$, while $S_j S_{j-1}$ and
$S_{j-1}^* S_j^*$ act on $Q^a$ like $C_1(a,j) +
C_2(a,j)Q^{p-1}$.
\end{proof}

\subsection{Mirror conjugate identity}\label{sec3-4}
The following remarkable identity has application to estimating
eigenvalues of $\L$ (see Section~\ref{sec:est}): \EQ{ \label{twist}
 S_j (L_{j-1}-\la_j) S_j^* = S_j^* (L_{j+2}-\la_j) S_j.}
To prove this, start with the formula
\begin{multline}
 (\p_x + R)(\p_x^2 + V)(\p_x -R)\\
 = \p_x^4 + (-3R'-R^2+V)\p_x^2 + (-3R'-R^2+V)'\p_x \\
 - R'''-(VR)'-RR''-R^2V,
\end{multline}
which implies that
$(\p_x+R)(\p_x^2+V_+)(\p_x-R)=(\p_x-R)(\p_x^2+V_-)(\p_x+R)$
is equivalent to \EQ{ \label{eq V}
 V_\pm = -R''/R \pm 3 R' - R^2 + C/R.}
Now set $R := aQ'/Q$. Plugging the following identities
\EQAL{
 &R^2 = a^2(1-\frac{2}{p+1}Q^{p-1}),\ R'= -a\frac{p-1}{p+1}Q^{p-1},\\
 &R''/R = -\frac{(p-1)^2}{p+1}Q^{p-1}.}
into \eqref{eq V}, we get, for $C=0$,
\EQ{
 V_\pm = - a^2 + \frac{2}{p+1}(a\pm(p-1))(a\pm(p-1)/2)Q^{p-1} .}
Hence for $a=k_j$ we have
\EQ{
 V_\pm = -k_j^2 + \frac{2}{p+1} k_{j\pm 2} k_{j\pm 1},}
which gives the desired identity \eqref{twist}.
The above proof also shows that $L_{j-1}$ and
$L_{j+2}$ are the unique choice for the identity
to hold with $S_j$ (modulo a constant multiple of
$Q/Q_x$, which is singular).

\subsection{Variational formulations for eigenvalues of $\L$}\label{sec3-5}
We considered two variational formulations for nonzero eigenvalues
of $\L$ in general dimensions in Section~\ref{var-Rn}.  Here we
present a new variational formulation for 1-D. Define the
selfadjoint operator
\begin{equation}
  H := S L_+ S^*.
\end{equation}
This is a fourth-order differential operator,
with essential spectrum $[1,\infty)$.
By a direct check, we have
\begin{equation*}
  H Q = S L_+ S^* Q = SL_+ (-2 Q_x) = 0.
\end{equation*}
Thus $Q$ is an eigenfunction with eigenvalue $0$.  Since
$(Q,S^*f)=(SQ,f)=0$ for any $f$, we have
\begin{equation}\label{eq3-33A}
\text{Range } S^* \perp Q.
\end{equation}
In particular, since $\Lp|_{Q^\perp}$ is nonnegative for $p\le 5$ by
Lemma \ref{TH2-1}, so is $H$.

\begin{lemma} \label{th3-5A} The null space of $H$ is
\[
N(H) = \text{span} \bket{ Q , \pispc xQ },
\]
where, recall, $\pispc$ is $0$ if $p \not= p_c$, and $1$ if $p=p_c$.
\end{lemma}

{\it Remark.} Note that dim $N(H) = 1 + \pispc$ which is different
from dim $N(\Lm^{1/2} \Lp \Lm^{1/2}) = 2 + \pispc$. We will show below
that $H$ and $\Lm^{1/2} \Lp \Lm^{1/2}$ have the same
{\em nonzero} eigenvalues.

\smallskip

\begin{proof}
If $Hf =0$, then $\Lp S^* f = -2c Q$ and $S^* f = c Q_1 + d Q_x$ for
some $c,d \in \R$ by Lemma \ref{TH2-2A}.  We have $Q_1 \perp Q$ iff $p
=p_c =5$. Thus, if $p \not = 5$, $c=0$ by \eqref{eq3-33A}, and
$S^* (f + \frac d2 Q)=0$. We conclude $f  = - \frac d2 Q$.

When $p = 5$, we have $S^* xQ = - Q^{-1} \pd_x (Q xQ) = -2Q_1 $. Thus
$S^* (f + \frac c2 Q_1 + \frac d2 Q)=0$ and $f = -\frac c2 Q_1 - \frac
d2 Q$.
\end{proof}

\smallskip

Define eigenvalues of $H$ as follows:
\begin{equation} \label{varC}
\tilde \mu_j := \inf _{f \perp f_k, k <j} \, \frac {(f, Hf)}{(f,f)},
\qquad (j=1,2,3,\ldots)
\end{equation}
with a suitably normalized minimizer denoted by $f_j$, if it exists.
By standard variational arguments,
if $\tilde \mu_j < 1$, then a minimizer $f_j$ exists.
By convention, if $\mu_k$ is the first of the $\mu_j$'s to
hit $1$ (and so $f_k$ may not be defined), we set
$\mu_j := 1$ for all $j > k$.

We can expand Summary \ref{TH2-5A} to the following.

\begin{theorem}[Equivalence] \label{th3-5}
Let $n=1$. Let $\mu_j$ be defined as in Summary \ref{TH2-5A} and
$\tilde \mu_j$ be defined by \eqref{varC}. Then $\mu_j = \tilde
\mu_j$.  When $\mu_j \not = 0$ and $\mu_j<1$, the eigenfunctions of
\eqref{varB} and \eqref{varC} can be chosen to satisfy
\[
  u_j = S^* f_j, \quad f_j = \frac 1{\mu_j} \,S \Lp u_j.
\]
\end{theorem}

\begin{proof}
First we establish the equivalence of nonzero eigenvalues.  Suppose $f
= f_j$ is an eigenfunction of \eqref{varC} with eigenvalue $\tilde \mu \not
=0$, then $S\Lp S^* f = \tilde\mu f$. Let $u := S^* f \not= 0$
and apply $S^*$ on both sides.
By $\Lm = S^* S$ we get $\Lm \Lp u = \tilde \mu u$.  Thus $u$ is
an eigenfunction satisfying \eqref{eq2-8} with $\mu=\tilde\mu$.
On the other hand, suppose
$u$ satisfies $\Lm \Lp u = \mu u$ with $\mu \not = 0$.  Applying
$S\Lp$ on both sides and using $\Lm = S^* S$, we get $S\Lp S^* S \Lp u
= \mu S \Lp u$, i.e., $H f = \mu f$ for $f = \mu^{-1}S \Lp u$.

Now use Lemmas \ref{TH2-3} and \ref{th3-5A}. If $p \in (1,5)$, then
$\mu_1 =\tilde \mu_1= 0$, corresponding to $Q_x$ and $Q$, and $\mu_2 =
\tilde \mu_2 >0$. If $p=5$, then $\mu_1 = \mu_2 =\tilde\mu_1 =
\tilde\mu_2 = 0$, corresponding to $Q_x, Q_1$, and $Q,xQ$, and
$\mu_3 =\tilde \mu_3=1$.  If $p \in (5,\infty)$, then $\mu_1 = \tilde
\mu_1 <0$, $\mu_2 = \tilde \mu_2 =0$, corresponding to $Q_x$ and $Q$,
and $\mu_3=\mu_3 =1$.  We have shown $\tilde \mu_j = \mu_j$.
\end{proof}

\medskip

In the following we will make no distinction between $\mu_j$ and
$\tilde \mu_j$.  By the minimax principle, \eqref{varC} has the following
equivalent formulations:
\begin{equation}\label{varC2}
\mu_j  = \inf_{\dim M= j}\, \sup_{f \in M}  \frac {(f, Hf)}{(f,f)}
= \sup_{\dim M= j-1} \, \inf_{f \perp M}  \frac {(f, Hf)}{(f,f)}.
\end{equation}
Here $M$ runs over all linear subspaces of $L^2(\R)$ with the
specified dimension.

\subsection{Estimates of eigenvalues of $\L$}
\label{sec:est}

In this subsection we prove lower and upper bounds for eigenvalues
of $\L$, confirming some aspects of the numerical computations shown
in Figure~\ref{fig:cL1d}.  Recall that, by Lemma~\ref{TH2-3}, the
first positive $\mu_j$ is $\mu_2$ for $p \in (1,p_c)$ and $\mu_3$
for $p \in [p_c,\pmax)$.  The first theorem concerns upper bounds
for $\mu_1$ and $\mu_2$.

\medskip

\begin{figure}[ht]
\centering
\psfrag{x}{$\log_{10}(p-1)$}\psfrag{y}{$\mu_j(p)$}\psfrag{z}{$p$}
\psfrag{a}{$\mu_1$}\psfrag{b}{$\mu_2$}\psfrag{c}{$\mu_3$}
\psfrag{d}{$\mu_4$}\psfrag{e}{$\mu_5$}
\epsfig{file=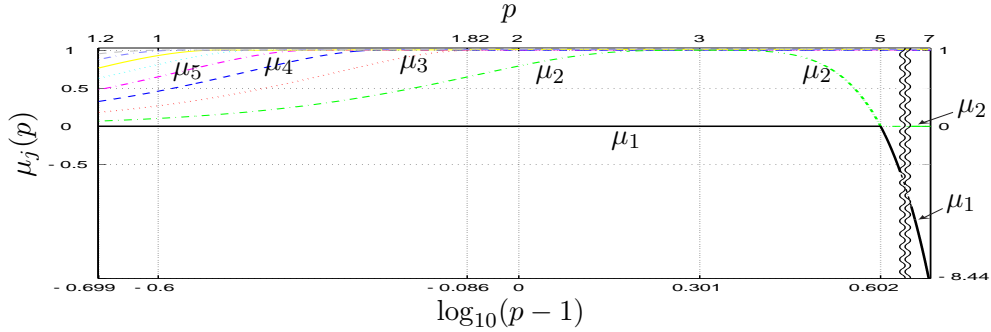,height=15cm,width=5cm,angle=270}
\vspace*{-1cm} \caption{$p$ vs. $\mu_j$.} \label{fig:mu_1D}
\end{figure}

\begin{theorem} \label{th3-6}
Suppose $n=1$ and $1<p<\infty$.

(a) If $p \not = 3$, then $\mu_2 \le C_p$ for some explicitly
computable $C_p < 1$. In particular $f_2$ exists.

(b) $\mu_1<0$ if and only if $p>5$. For any $C>0$,
we have $\mu_1(p) \le - C p^3$ for $p$
sufficiently large.

\end{theorem}

\begin{proof}
For part (a), we already know $\mu_2 = 0$ for $p \geq 5$. Assume $p \in
(1,5)$.  Consider test functions of the form $f = SQ^k$ with $k>0$.
$f$ is odd and hence $f\perp Q$, the $0$-eigenfunction of $H$.  Since
$H=S\Lp S^*$ and $S^*S=\Lm$, we have
\begin{equation*}
\mu_2 \le
\frac {(f, Hf)}{(f,f)} = \frac {(L_-Q^k, L_+ L_-Q^k)}{(Q^k,L_-Q^k)} .
\end{equation*}
By formulas \eqref{eq3-5} and \eqref{eq3-8},
\begin{equation*}
L_- Q^k = a Q^{k+p-1} + b Q^k, \quad a  = \frac 1{p+1} (k-1)(2k+p+1)
, \quad b = 1-k^2.
\end{equation*}
\begin{equation*}
L_+ Q^{k+p-1} = \sigma  Q^{k+2p-2} + d Q^{k+p-1},
\end{equation*}
\begin{equation*}
 \sigma 
= \frac 1{p+1} (k+2p-1)(2k+p-3)
, \quad d = 1-(k+p-1)^2.
\end{equation*}
\begin{equation*}
L_+ Q^k = c Q^{k+p-1} + b Q^k, \quad c  = \frac 1{p+1} (k+p)(2k-p-1).
\end{equation*}
Thus
\begin{equation} \label{fHf/ff}
\frac {(f, Hf)}{(f,f)} = \frac {a^2\sigma  J_3 + a(ad + bc+ b\sigma )J_2 + b
(ad+ab + bc) J_1 + b^3 J_0} {a J_1 + b J_0}
\end{equation}
where
\begin{equation*}
J_m = \int_\R Q^{2k+ m(p-1)}(x)\,dx, \quad (m=0,1,2,3),
\end{equation*}
which are always positive.  If $k \to 0^+$, then $J_m$ converges to
$\int_\R Q^{m(p-1)}\,dx$ for $m > 0$, and $J_0 = O(k^{-1})$.
The above quotient can be written as
\[
\eqref{fHf/ff} = b^2 + \frac{ J} {a J_1 + b J_0}
\]
where
\[
J = a^2\sigma  J_3 + a(ad + bc+ b\sigma )J_2 + b
(ad + bc) J_1.
\]
Note that $J_m|_{k=0} = (\frac{p+1}2)^m \, \frac2{p-1}\, \int_\R
\sech^{2m}( y)\, dy$ with $\int_\R \sech^{2m}( y)\, dy= 2,\frac 43,
\frac {16}{15} $ for $m=1,2,3$, respectively.  Also, as $k \to 0^+$, $
a \to -1$, $b\to 1$, $ c \to -p$, $\sigma  \to \frac {(2p-1)(p-3)}{p+1}$,
and $d \to 1 - (p-1)^2$. Direct calculation shows
\[
\lim_{ k \to 0^+} J =- \frac {2}{15(p-1)} \, (p+1)^2(p-3)^2.
\]
Also note $b^2 < 1$ for $k>0$.  Thus, if $1<p<\infty$ and $p \not =
3$, then $J<0$ and the quotient \eqref{fHf/ff} is less than $1$ for
$k$ sufficiently small. (If $p=3$, the sign of $J$ is unclear and
\eqref{fHf/ff} may not be less than $1$.)  This proves $\mu_2 < 1$ and
provides an upper bound less than $1$ for $\mu_2$. It also implies the
existence of $f_2$. This establishes statement
(a).

For statement (b), the fact that $\mu_1 <0$ if and only if $p>5$ is
part of Lemma~\ref{TH2-3}.  We now consider the behavior of $\mu_1$
for $p$ large.  Fix $k>1$ to be chosen later. As $p \to \infty$,
\[
J_m = (\frac {p+1}2)^{\frac {2k}{p-1} + m}
\cdot \frac 2{p-1} \cdot \int_\R \bke{\sech x}^{\frac {4k}{p-1} + 2m}\, dx
\sim C_m p^{m-1},
\]
with $C_m = 2^{1-m} \int_\R (\sech x)^{2m} dx= 2,\frac 23,
\frac {4}{15} $ for $m=1,2,3$, respectively, and
\[
a \sim k-1, \quad
b = 1- k^2, \quad
c \sim -p, \quad
\sigma \sim 2p, \quad
d \sim - p^2.
\]
Thus, by \eqref{fHf/ff},
\begin{equation*}
 \frac {(f, Hf)}{(f,f)} \sim  \frac
{a\sigma  J_3 + a d J_2}{J_1 }
\sim \frac{1-k}{15}\,p^3
\quad \text{ as } p \to \infty .
\end{equation*}
By choosing $k>1$ sufficiently large, we have shown
that for any $C$, $\mu_1 \le -C p^3$ for $p$ sufficiently large.
\end{proof}

\medskip

The next theorem bounds eigenvalues of $\L$ by eigenvalues of
$\Lp$ and $\Lm$. Recall $p_j$ and $\la_j(p)$ are
defined in \eqref{eq3-4} and \eqref{lam.def}.

\begin{theorem}[Interlacing of eigenvalues] \label{th3-7}
Fix $k \ge 1 $ and $p \in [p_{k+2} ,\ p_{k+1})$ where, recall,
$p_j = \frac{j+1}{j-1}$.  Let $\la_j(p) = 1- \frac 14[(p+1)-j(p-1)]^2$ be as in
\eqref{lam.def} and so $\la_{k+1} < 1 \le \la_{k+2}$.  For the
eigenvalues $\mu_j$ defined by \eqref{varC}, we have
\begin{equation} \label{eq:th3-7}
\la_{j+1}^2(p) <  \mu_{j+1}(p) < \la_{j+2}^2 (p),\quad (1\le j<k ); \qquad
\la_{k+1}^2(p) <  \mu_{k+1}(p) \le 1.
\end{equation}
In particular, there are $K$ simple eigenvalues $\mu_2, \ldots, \mu_{K+1}$
in $(0,1)$ where $K=k$ if $\mu_{k+1}<1$ and $K=k-1$ if $\mu_{k+1}=1$.
Moreover,  $K$ is always 1 when $k=1$. Finally,
\[
 \mu_2 \ge
 \begin{cases}
  \la_2 \la_3 &(1<p \leq 2),\\
  \la_2 &(2<p<5),
 \end{cases}
\qquad
 \mu_3 \ge
 \begin{cases}
  \la_3 \la_4 & (1 < p \leq 5/3),\\
  \la_3 & (5/3 < p \leq 2),\\
  1 & (2<p<\infty),
 \end{cases}
\]
\[
\mu_1 \ge -\frac{1}{16}(p-1)^3(p-5) \quad (5 \leq p < \infty).
\]
\end{theorem}

\begin{figure}[ht]
\centering
\psfrag{x}{$\log_{10}(p-1)$}\psfrag{y}{$f_j(p)$}\psfrag{z}{$p$}
\psfrag{a}{$p_2$}\psfrag{b}{$p_3$}\psfrag{c}{$p_4$}
\psfrag{d}{$p_5$}\psfrag{e}{$p_6$}
\psfrag{g}{$f_1$}\psfrag{h}{$f_2$}\psfrag{i}{$f_3$}\psfrag{j}{$f_4$}\psfrag{k}{$f_5$}
\epsfig{file=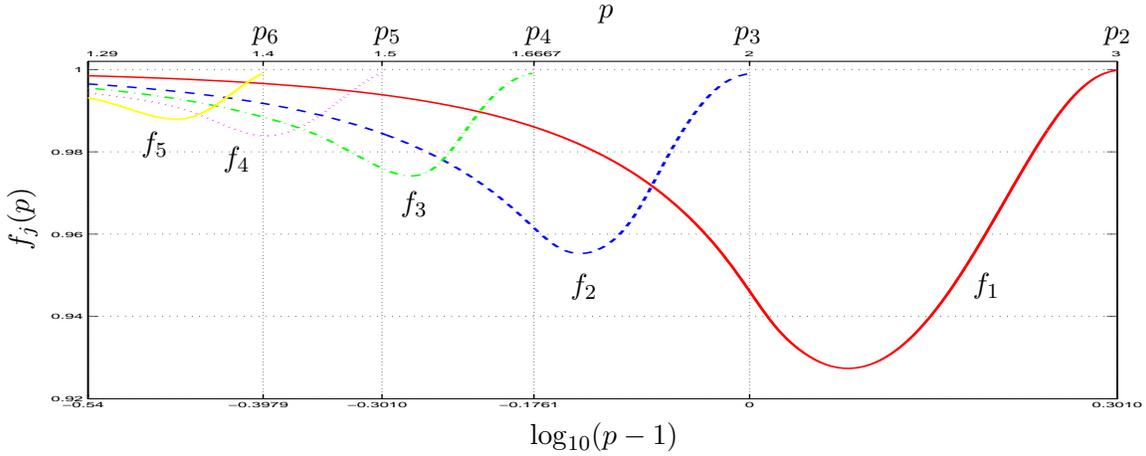,height=15cm,width=5cm,
angle=270}\vspace*{0.3cm}
 \caption{$p$ vs. $f_j$ for $j=1,\ldots,5$, where
$f_j(p)=\frac{\lambda_{j+1} \lambda_{j+2}}{\mu_{j+1}}$ for $1<p<p_{j+2}$
and $f_j(p)=\frac{\lambda_{j+1} }{\mu_{j+1}}$ for $p_{j+2}\le p<p_{j+1}$.}

\label{fig:lam_mu}
\end{figure}

\begin{rem}
In view of the above lower bounds for $\mu_2$ and $\mu_3$, we
conjecture that
\begin{equation}
\mu_{j+1} \ge \la_{j+1}  \la_{j+2} \quad (1<p< p_{j+2}); \qquad
\mu_{j+1} \ge \la_{j+1}   \quad (p_{j+2} \le p< p_{j+1}).
\end{equation}
This is further confirmed numerically for $j=3,4,5$ (see
Figure~\ref{fig:lam_mu}). Note that $\lim_{p \to p_{j+1}-}
\frac{\la_{j+1}}{\mu_{j+1}} =1$ because both $\la_{j+1}$ and
$\mu_{j+1}$ converge to $1$. It also seems that $\frac{\la_{j+1}
\la_{j+2}}{\mu_{j+1}}$ has a limit as $p \to 1+$, but it is not
clear although we have \eqref{eq:th3-7} and $\la_j =
(j-1)(p-1)+O((p-1)^2)$ as $ p \to 1+$.
\end{rem}

\begin{proof}
We first prove the upper bound: For $j<k$, use the test functions
\begin{equation*}
S\psi_2, \ S\psi_3, \ \ldots, \ S \psi_{j+2}
\end{equation*}
(we cannot use $S\psi_1$ since it is zero). Recall $\Lm \psi_m = \la_m
\psi_m$. Let $a=(a_2,\ldots, a_{j+2})$ vary over $\C^{j+1}-
\bket{0}$.  By equivalent definition \eqref{varC2}, $H=S \Lp S^*$,
$\Lm = S^* S$, and the orthogonality between the $\psi_m$'s, we have
\begin{align*}
\mu_{j+1} & \le \sup _a \frac {
(\sum_{m}  a_m S \psi_m, H \sum_{\ell}  a_\ell S \psi_\ell)}
{(\sum_{m}  a_m S \psi_m,  \sum_{\ell}  a_\ell S \psi_\ell)}
= \sup _a \frac {
(\sum_{m}  a_m  \psi_m, \Lm \Lp \Lm \sum_{\ell}  a_\ell  \psi_\ell)}
{(\sum_{m}  a_m  \psi_m,  \Lm \sum_{\ell}  a_\ell  \psi_\ell)}
\\
& \le \sup _a \frac {
(\sum_{m}  a_m  \psi_m, \Lm \Lm \Lm \sum_{\ell}  a_\ell  \psi_\ell)}
{(\sum_{m}  a_m  \psi_m,  \Lm \sum_{\ell}  a_\ell  \psi_\ell)}
= \sup _a \frac {\sum_{m} |a_m|^2 \la_m^3} {\sum_{m} |a_m|^2 \la_m}
\\
& \le\max_{m=2,\ldots, j+2} \la_m^2 = \la_{j+2}^2.
\end{align*}
Since $\mu_{j+1} \le \la_{j+2}^2<1$, it is attained
at some function, for which the second inequality above
cannot be replaced by an equality sign. Thus
$\mu_{j+1} < \la_{j+2}^2$.

\medskip

For the lower bound of eigenvalues, we use only
the special case $j=1$ of \eqref{twist}: \EQ{
 H = S L_+ S^* = S L_0 S^* = S^* L_3 S.}
In particular, we have for $1<p<3$,
\EQ{
 H \ge S^* L_2 S = S^* S S^* S = L_1^2 = L_-^2,}
which implies that
\EQ{
 \la_{j+1}^2 \le \mu_{j+1} \quad (1<p<3)}
(and again, equality is impossible).

For the second eigenvalue $\mu_2$, we can get a
more precise estimate by using \eqref{low bd} for
$L_3\ge \la_3'$ together with \EQ{
 L_1|_{Q^\perp} \ge \la_2',}
which follows from $spec(L_1)$. Combining these
estimates, we have for any $f\perp Q$ and $p<5$,
\EQ{
 (Hf,f) \ge \la_3'(Sf,Sf) \ge \la_3' \la_2' (f,f),}
which implies that $\mu_2\ge \la_3'\la_2'$, i.e.,
\EQ{ \label{better bd}
 \mu_2 \ge
 \begin{cases}
  \la_2 \la_3 &(1<p \leq 2),\\
  \la_2 &(2<p<5)
 \end{cases}}
For $p>3$, we have $L_3 \ge \la_2=-(p-1)(p-5)/4$ and
\EQ{
 L_3 - L_2 \ge -(p-1)(p-3)/2=:-a.}
Hence for any $t\in[0,1]$, we have
\EQ{
 L_3 \ge t L_2 - a t + (1-t)\la_2.}
and so for $b>0$, we have
\EQAL{
 (Hf,f)+b(f,f) &\ge (S^* (tL_2-at+(1-t)\la_2) Sf,f)+b(f,f)\\
 &=t\|L_1 f\|^2 -(at-(1-t)\la_2)(L_1 f,f) + b\|f\|^2,}
which is nonnegative if
\EQ{
 b \ge (at-(1-t)\la_2)^2/(4t),}
whose infimum is attained at
$t =-\la_2/(a + \la_2) = (p-5)/(p-1) \in (0,1)$ for $p > 5$.
Plugging this back in, we obtain the lower bound
\EQ{
 \mu_1 \ge \lambda_2(a+\lambda_2)
 = -\frac{1}{16}(p-1)^3(p-5) \quad (p>5).}

We have a similar bound on $\mu_3$ by using the even-odd decomposition
$L^2(\R)=L^2_{ev}(\R)\oplus L^2_{od}(\R)$.  Let $\psi_j$, $\xi_j$
be the eigenfunction of $L_1$ and $L_3$ such that
\EQ{
 L_1 \psi_j = \la_j \psi_j,\ L_3 \xi_j=\la_j \xi_j.}
$\psi_j$ starts from $j=1$ and $\xi_j$ starts with $j=3$. They are
even for odd $j$ and odd for even $j$. For any even function $f\perp
Q=\psi_1$, $Sf$ is odd and so we have $f\perp \psi_1=Q, \psi_2$ and
$Sf\perp\xi_3$. Hence by $spec(L_3)$ and $spec(L_1)$, we have
\EQ{
 (Hf,f)=(L_3 Sf,Sf)\ge \ti\la_4(Sf,Sf) = \ti\la_4(L_1 f,f)
\ge \ti\la_4\ti\la_3(f,f),}
where we denote
\EQ{
 \ti\la_j:=
 \begin{cases}
  \la_j & (1<p<p_j),\\
  1 & (p_j<p).
 \end{cases}}
Thus the second eigenvalue of $H$ on $L^2_{ev}$
is $\ge\ti\la_4\ti\la_3$. Next for any odd
function $f\perp\psi_2$, we have
$f\perp\psi_1,\psi_2,\psi_3$. Hence we have
\EQ{
 (Hf,f) = (L_3 Sf,Sf) \ge \la_3'(Sf,Sf) = \la_3'(L_1 f,f)
\ge \la_3'\ti\la_4(f,f).}
Similarly, every odd function $f\perp S^*\xi_3$
satisfies $f\perp\psi_1$ and
$Sf\perp\xi_3,\xi_4$, so
\EQ{
 (Hf,f) \ge \ti\la_5\ti\la_2(f,f).
}
Hence the second eigenfunction on $L^2_{od}$ is
$\ge\max(\ti\la_4\la_3',\ti\la_5\ti\la_2)\ge\ti\la_4\ti\la_3$. Therefore
we have $\mu_3 \ge \ti\la_3\ti\la_4$, i.e.,
\EQ{
 \mu_3 \ge
 \begin{cases}
  \la_3 \la_4 & (1<p<5/3),\\
  \la_3 & (5/3<p<2),\\
  1 & (2<p).
 \end{cases}}
This argument, however, does not yield any useful estimates for the
higher $\mu_j$.
\end{proof}

\subsection{Resonance for $p=3$} \label{sec3-7}
In the theory of dispersive estimates for the linear Schr\"odinger
evolution, it is important to know whether or not the endpoints of the
continuous spectrum of the linear operator are eigenvalues or
{\it resonances}.  For our $\L$, the endpoints are $\la = \pm i$.
Resonance here refers to a function $\phi$ which satisfies the eigenvalue
problem locally in space with eigenvalue $i$ or $-i$, but which does
not belong to $L^2(\R^n)$. For dimension $n=1$, one requires $\phi \in
L^\infty(\R)$.
(Note for comparison's sake that in one dimension,
the operator $-d^2/dx^2$ has a
resonance -- corresponding to the constant function --
at the endpoint $0$ of its continuous spectrum.)

Before we made the numerical calculation, we did not expect to see
any resonance. However, from Figure~\ref{fig:cL1d}, one sees that
$\ka = \sqrt{ \mu_2}$ converges to $1$ as $p \to 3$. What does the
point $\ka=1$ at $p=3$ correspond to? A natural conjecture is that
it is a resonance or an eigenvalue, since the $p=3$ case is well-known
to be completely integrable and special phenomena may occur.

This is indeed the case since we have the following solution to the
eigenvalue problem~\eqref{eq2-6} when $p=3$,
\begin{equation}
  \phi = \begin{bmatrix} 1- Q^2 \\ i \end{bmatrix}, \quad \la = i.
\end{equation}
It is clear that $\phi \in L^\infty(\R)$ but $\phi \not \in L^q(\R)$
for any $q< \infty$.

Let $u_p(x)$ denote the real-valued (and suitably normalized) solution
of \eqref{eq2-8} corresponding to $\mu=\mu_2$. It is the first
component of the eigenfunction of \eqref{eq2-6}.  A natural question
is: does $u_p(x)$ converge in some sense to $u_3(x) := 1-Q^2(x)$ as $p \to
3$?  Since $u_p -u_3$ is not in $L^q(\R)$ for all $q\in [1,\infty)$,
it seems natural to measure the convergence in the following weighted
norm,
\[
  \|f\|_{\rm w} := \int_\R {\rm w}(f)^2(x)~{dx},
\]
where a weighting operator ${\rm w}$ is defined by
${\rm w}(f)(x) := f(x)\frac{1}{\sqrt{1+x^2}}$. This
de-emphasizes the value of $u_p-u_3$ for $x$ large, and so it should
converge to $0$ as $p$ goes to $3$. This is
confirmed numerically as follows.

\begin{figure}[ht]
\centering
\psfrag{z}{$p=2.8$}\psfrag{x}{}\psfrag{y}{$u_{2.8}(x)$}\psfrag{a}{(a)}
\epsfig{file=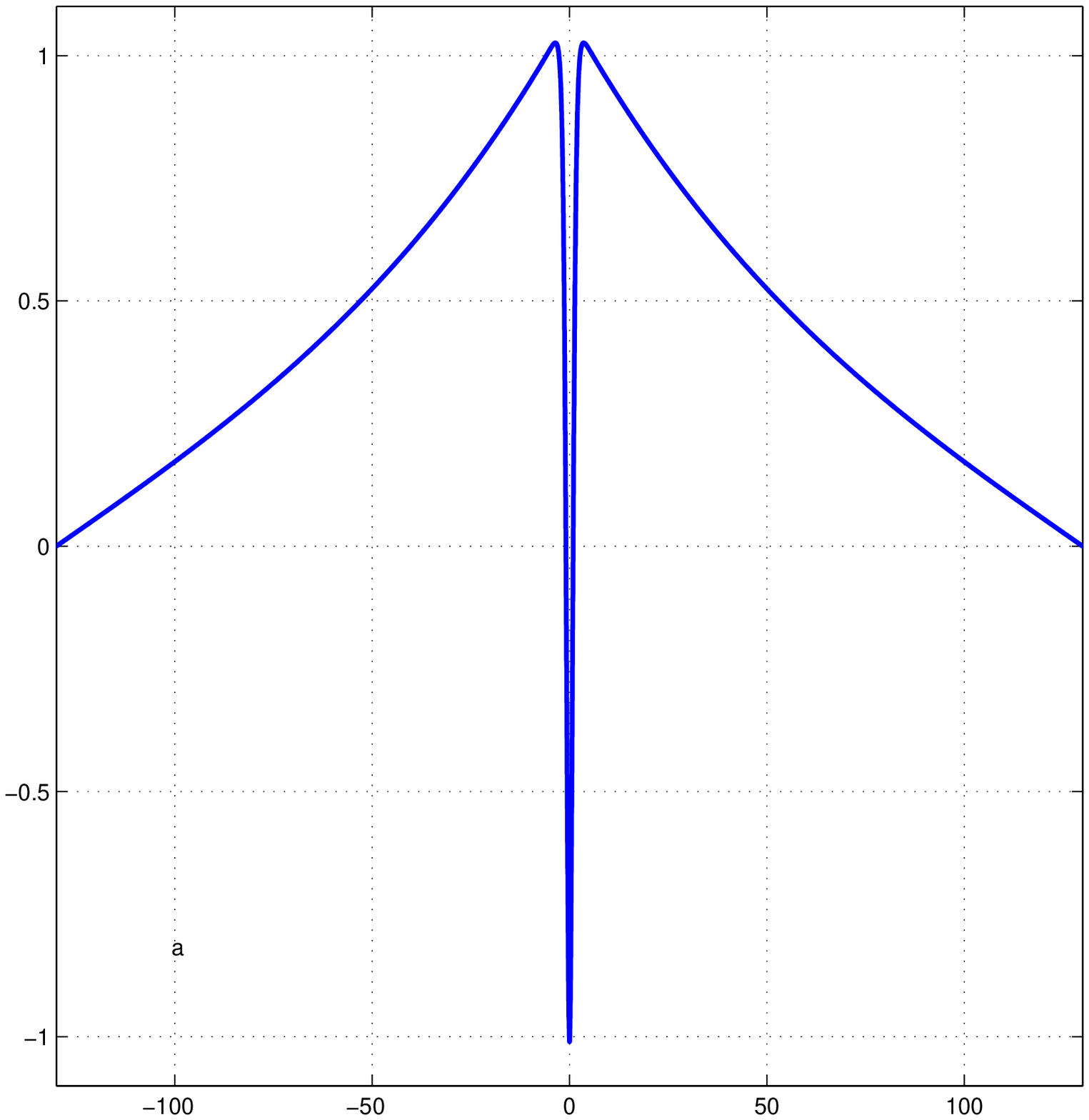,height=2.5cm,width=2.5cm}~~~~~~~
\psfrag{z}{$p=2.9$}\psfrag{x}{}\psfrag{y}{$u_{2.9}(x)$}\psfrag{a}{(b)}
\epsfig{file=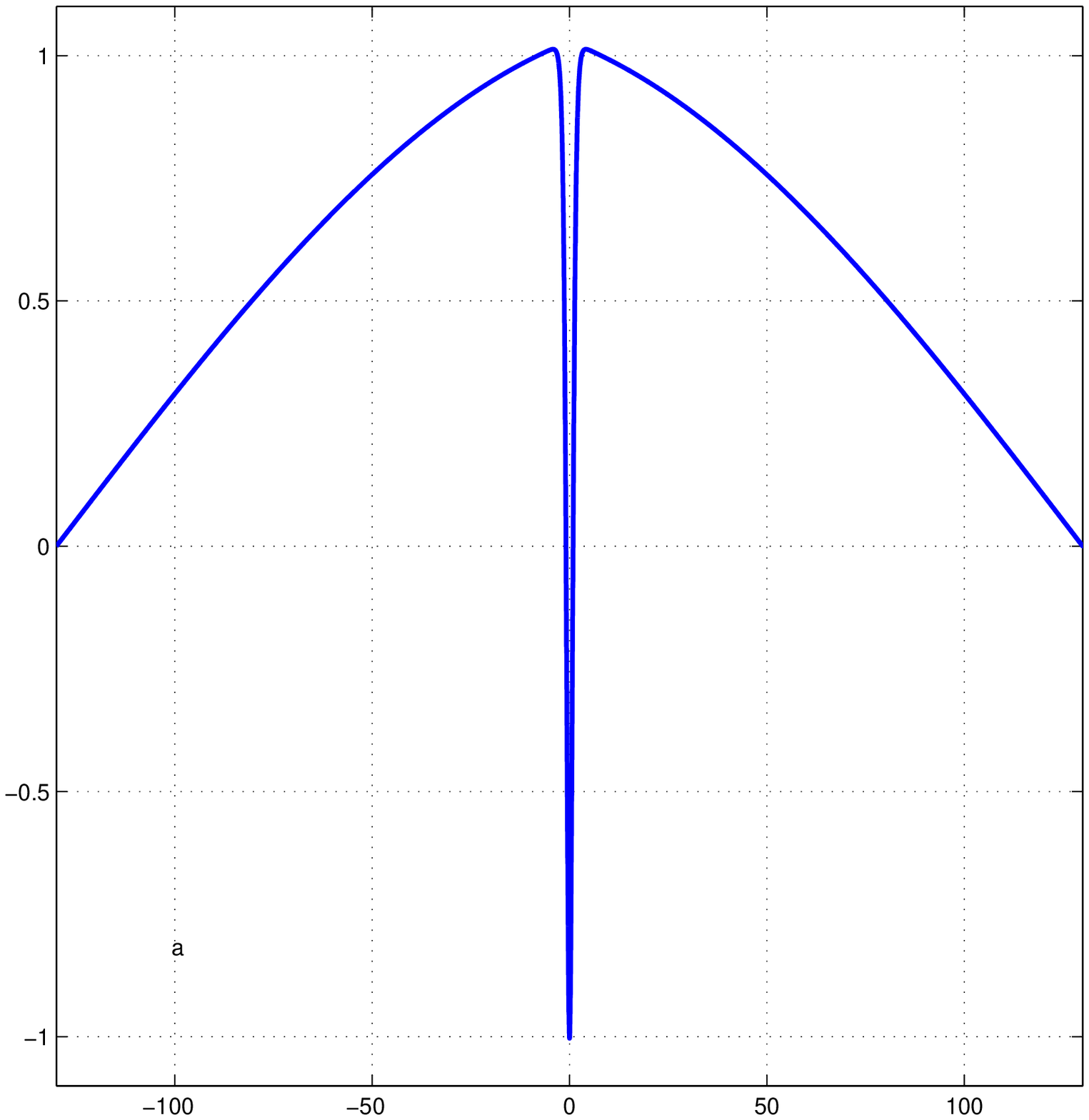,height=2.5cm,width=2.5cm}~~~~~~~
\psfrag{z}{$p=3.0$}\psfrag{x}{}\psfrag{y}{$u_{3}(x)$}\psfrag{a}{(c)}
\epsfig{file=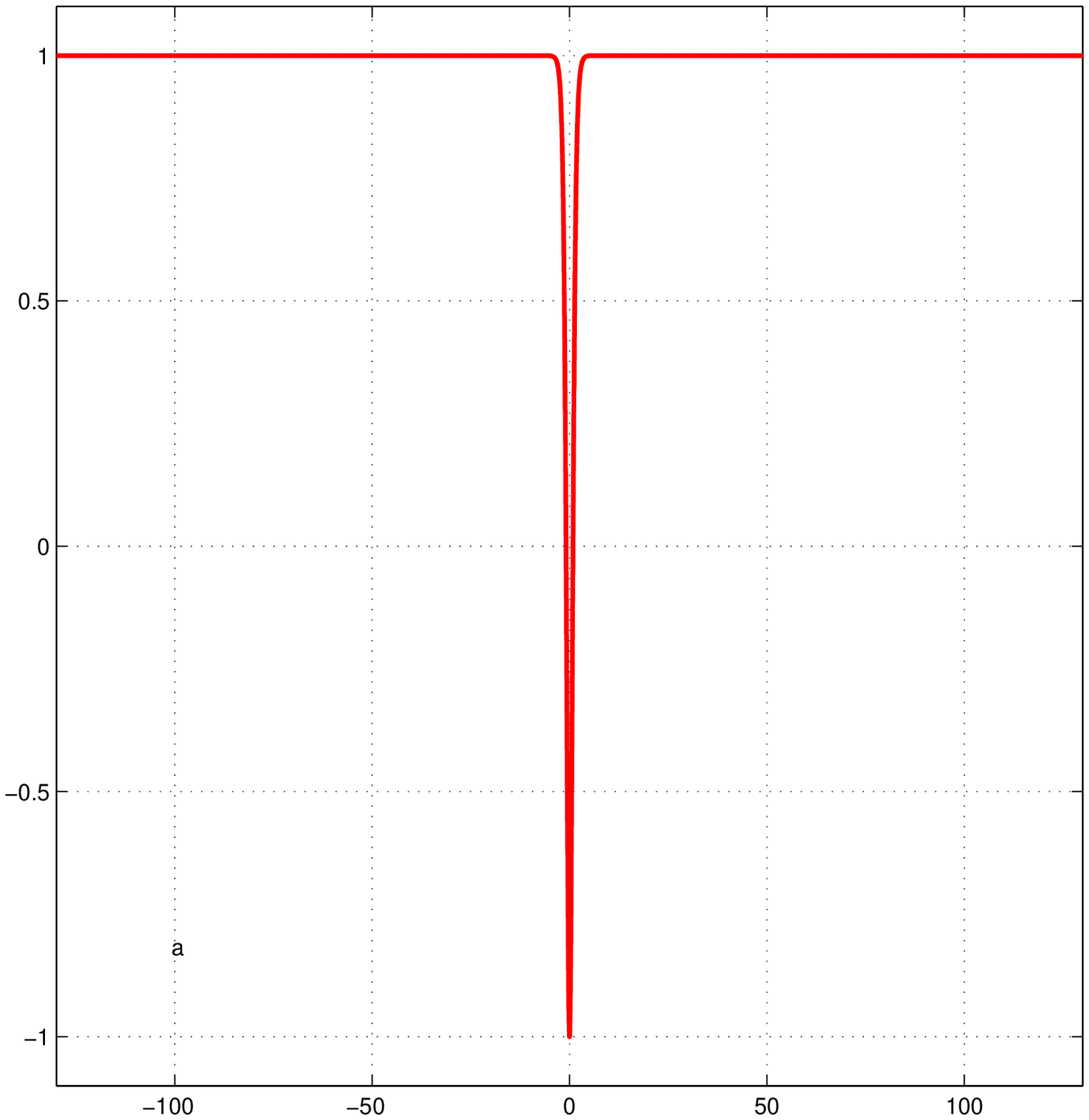,height=2.5cm,width=2.5cm}~~~~~~~
\psfrag{z}{$p=3.1$}\psfrag{x}{}\psfrag{y}{$u_{3.1}(x)$}\psfrag{a}{(d)}
\epsfig{file=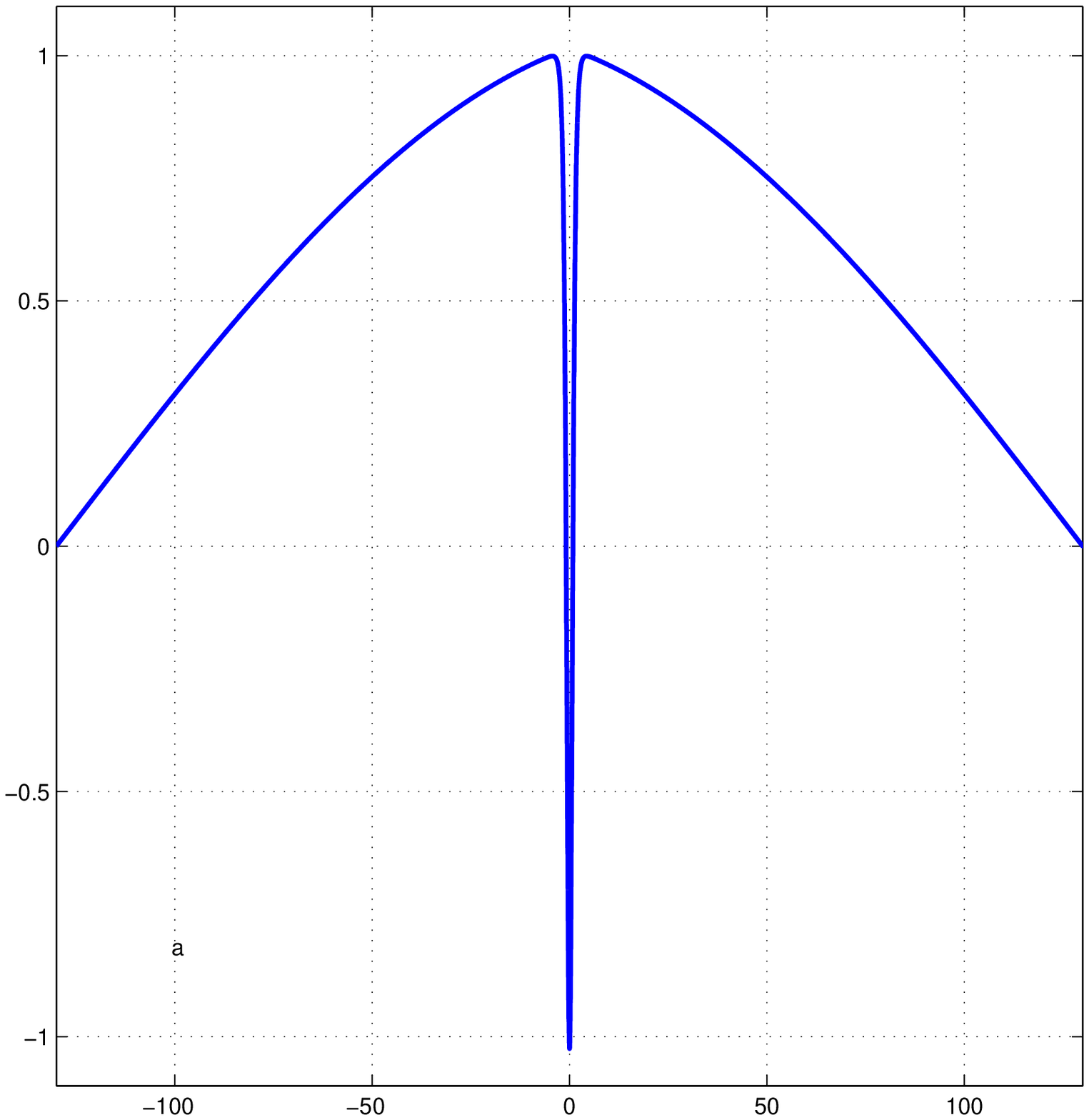,height=2.5cm,width=2.5cm}~~~~~~~
\psfrag{z}{$p=3.2$}\psfrag{x}{}\psfrag{y}{$u_{3.2}(x)$}\psfrag{a}{(e)}
\epsfig{file=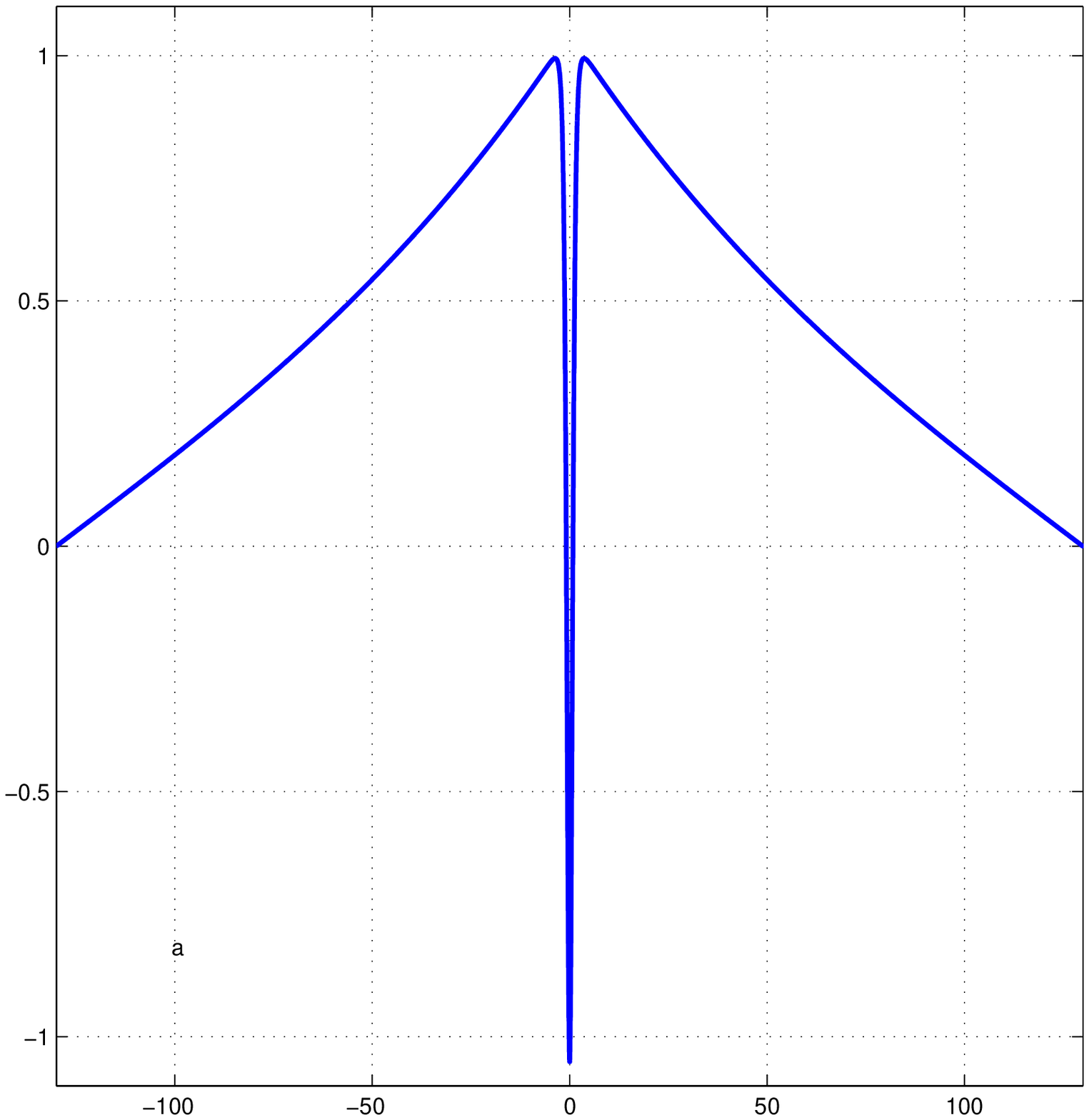,height=2.5cm,width=2.5cm}\vspace{0.2cm}\\
\psfrag{z}{}\psfrag{x}{$x$}\psfrag{y}{${\rm
w}(u_{2.8})(x)$}\psfrag{a}{(f)}
\epsfig{file=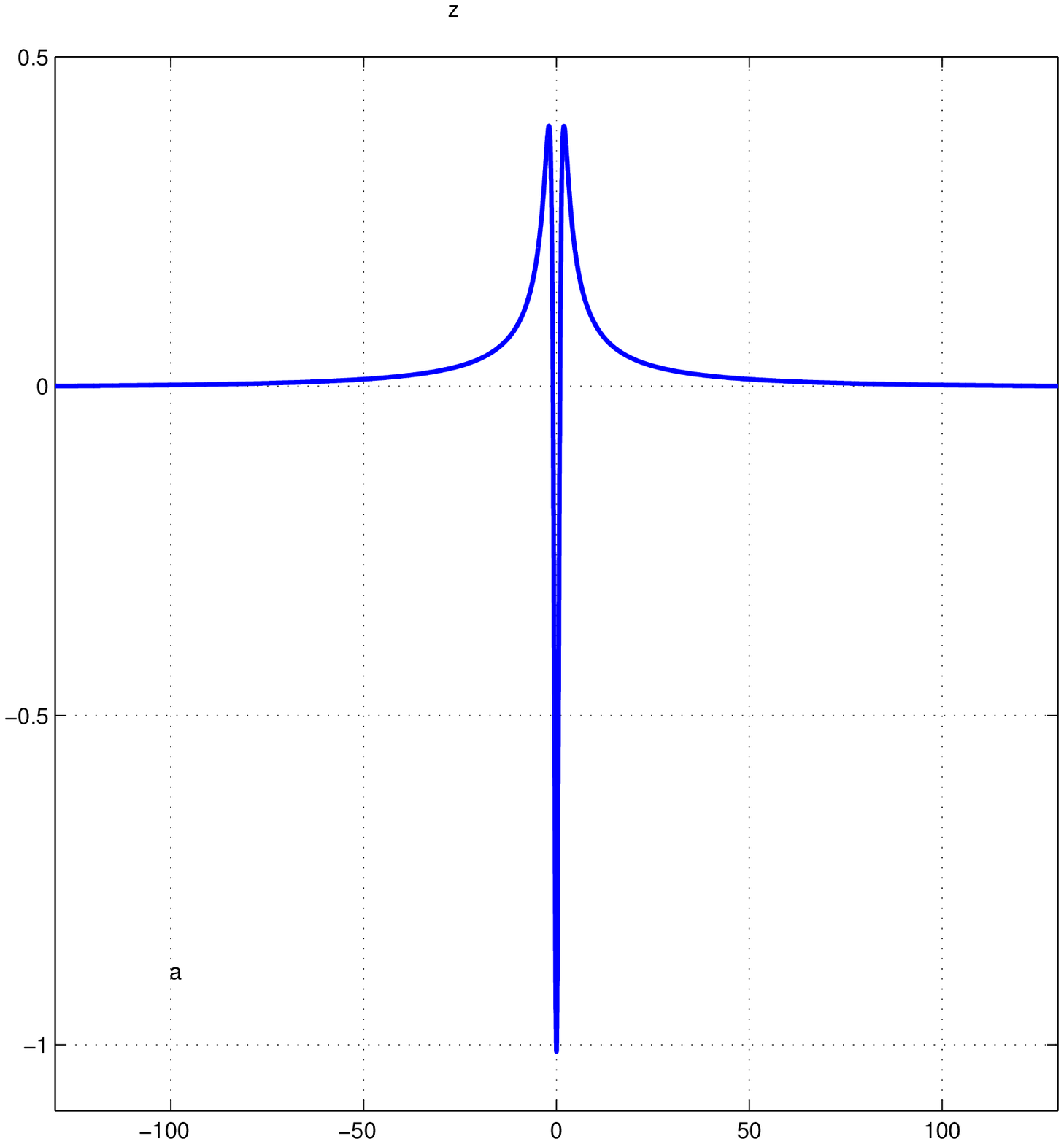,height=2.5cm,width=2.5cm}~~~~~~~
\psfrag{z}{}\psfrag{x}{$x$}\psfrag{y}{${\rm
w}(u_{2.9})(x)$}\psfrag{a}{(g)}
\epsfig{file=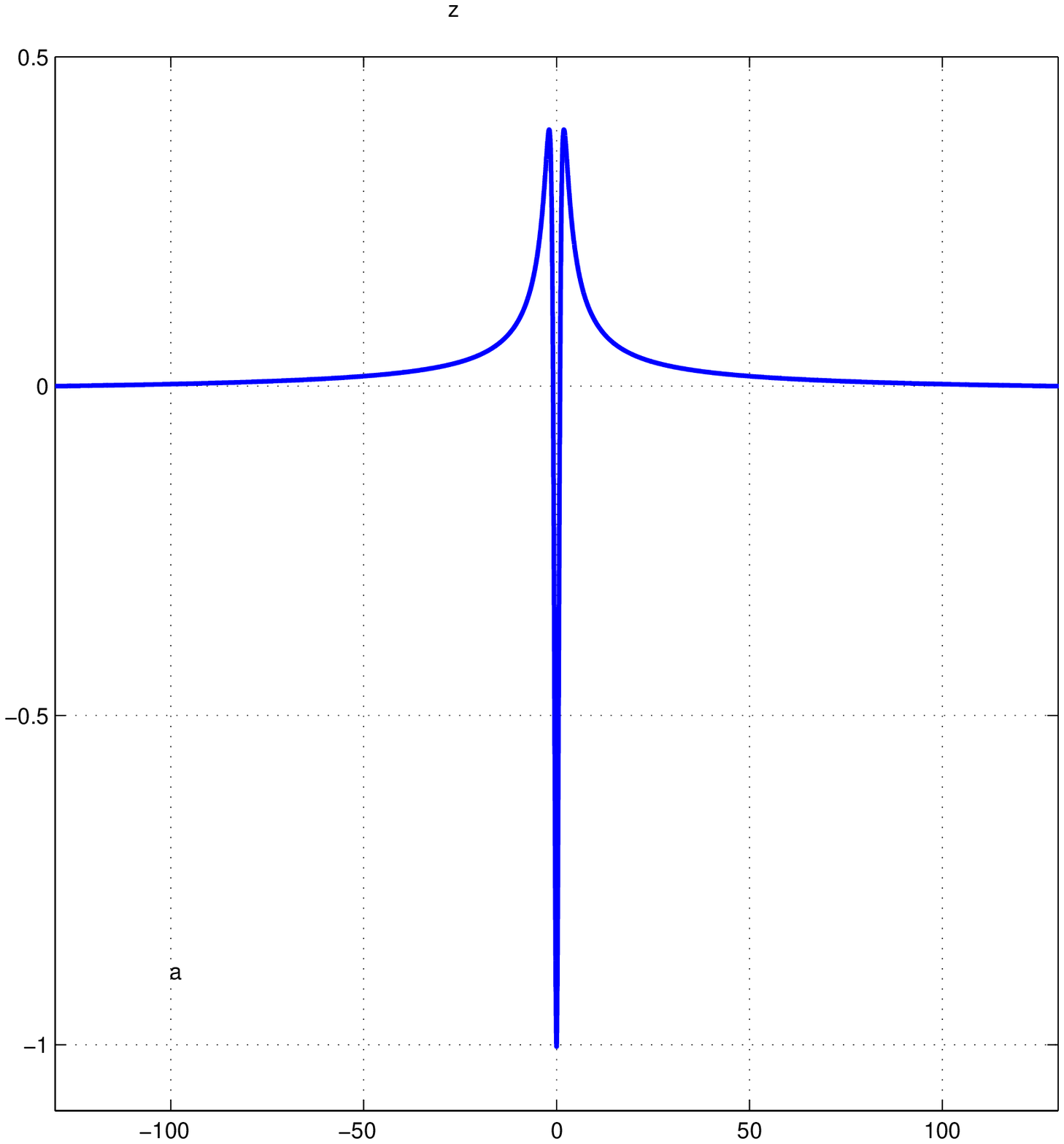,height=2.5cm,width=2.5cm}~~~~~~~
\psfrag{z}{}\psfrag{x}{$x$}\psfrag{y}{${\rm
w}(u_{3})(x)$}\psfrag{a}{(h)}
\epsfig{file=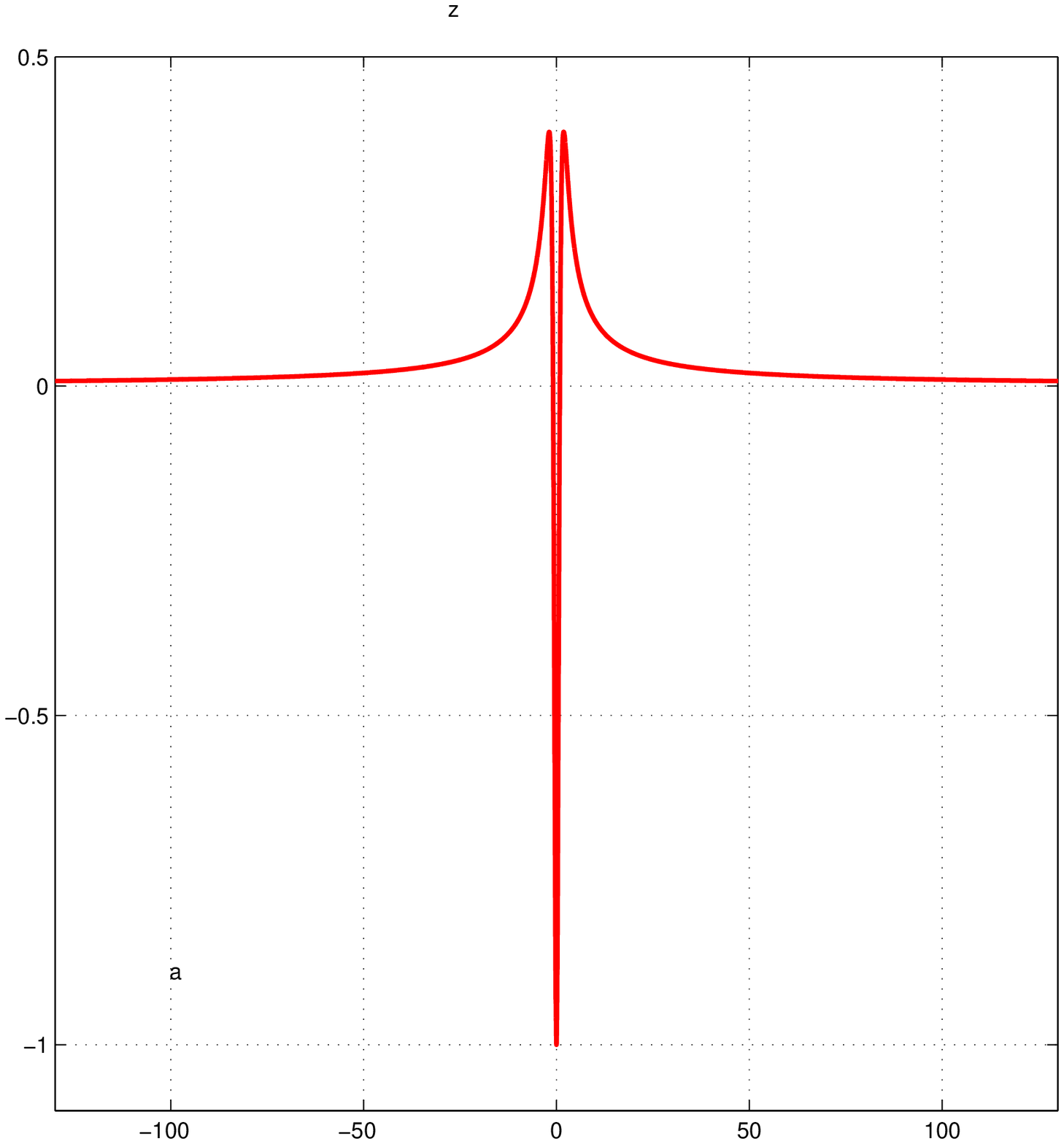,height=2.5cm,width=2.5cm}~~~~~~~
\psfrag{z}{}\psfrag{x}{$x$}\psfrag{y}{${\rm
w}(u_{3.1})(x)$}\psfrag{a}{(i)}
\epsfig{file=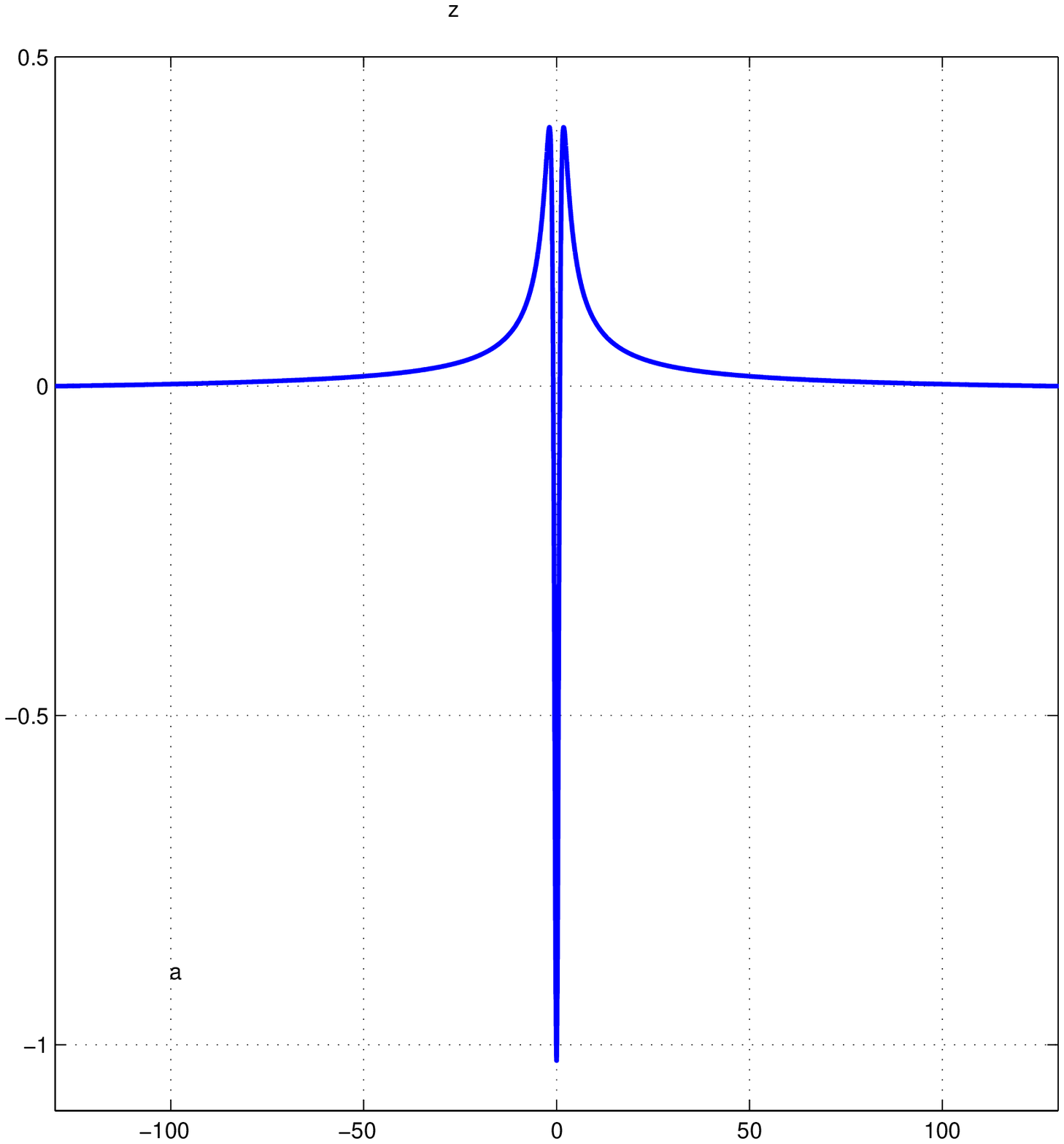,height=2.5cm,width=2.5cm}~~~~~~~
\psfrag{z}{}\psfrag{x}{$x$}\psfrag{y}{${\rm
w}(u_{3.2})(x)$}\psfrag{a}{(j)}
\epsfig{file=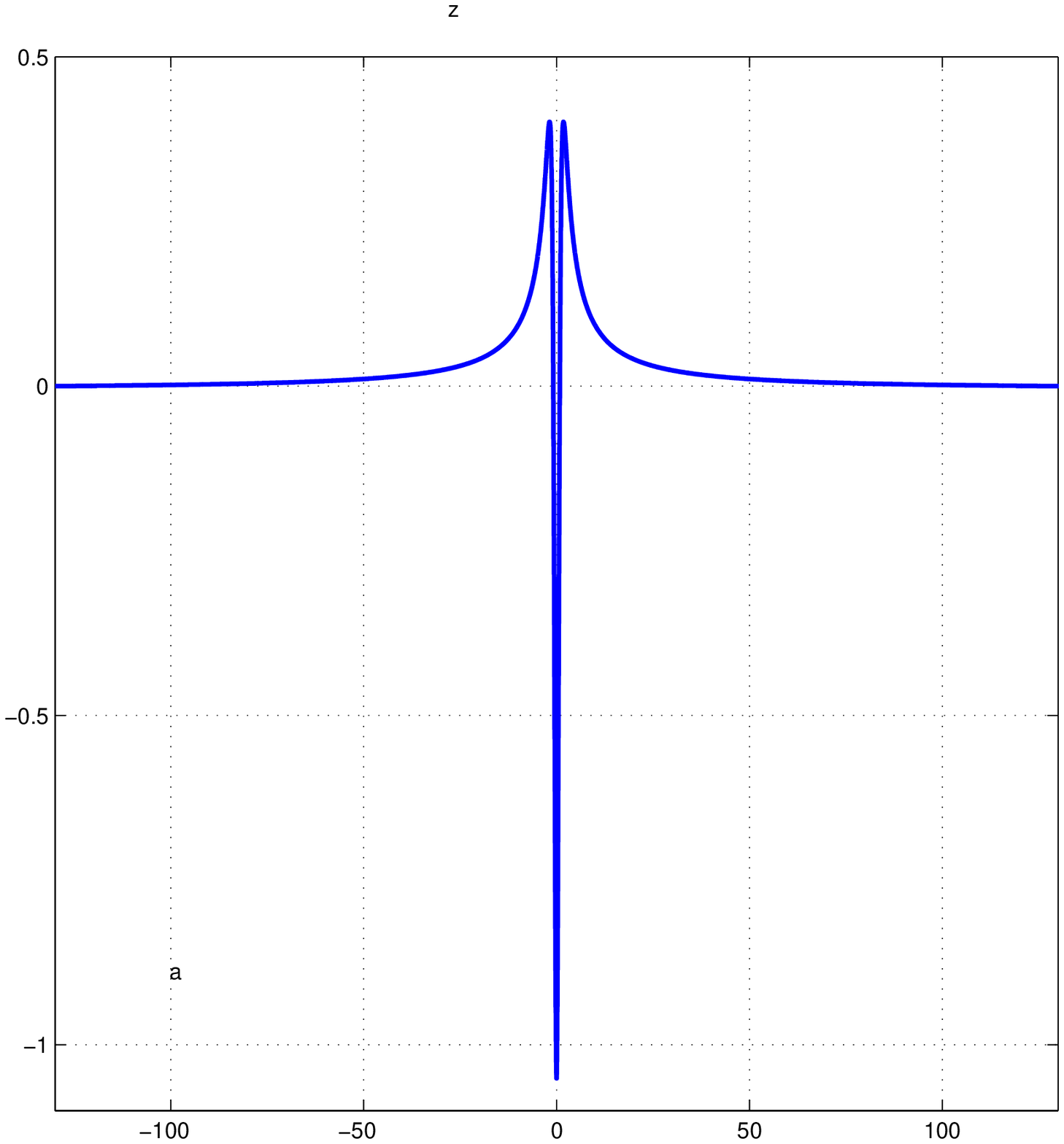,height=2.5cm,width=2.5cm}
\caption{$u_p(x)$ $\&$ ${\rm w}(u_p)(x)$ for
$p=2.8, 2.9, 3, 3.1$ and $3.2$.} \label{fig:up3}
\end{figure}

Let $u_3 := 1- Q^2$ and $\delta := \|u_3\|_{\rm w}$. In Appendix we
will propose a numerical method to solve for the eigenpair
$\{\lambda, [u_p(x), ~w_p(x)]^\top\}$ of \eqref{eq2-6} corresponding
to $\mu_2 = - \la^2$. Renormalize $u_p(x)$ for $p\ne 3$ so that it
is real-valued, $u_p(0)<0$, and $\|u_p\|_{\rm w}=\delta$. In
Figure~\ref{fig:up3}(c) we plot $u_3$ in a  large interval $|x|<130$
with $\delta = 1.3588$. According to the numerical method in
Appendix, we get $u_{2.8},u_{2.9},u_{3.1}$ and $u_{3.2}$ plotted in
Figure~\ref{fig:up3}(a), (b), (d) and (e), respectively. The
vertical range is roughly $[-1,1]$. In Figure~\ref{fig:up3}(f)--(j)
we plot ${\rm w}(u_p)$ for $p=2.8, 2.9, 3, 3.1$ and $3.2$, for
$|x|<130$ and vertical range $[-1,0.5]$.

\begin{figure}[ht]
\centering \psfrag{x}{$p$}\psfrag{y}{$\| u_p - u_3\|_{\rm w}$}
\epsfig{file=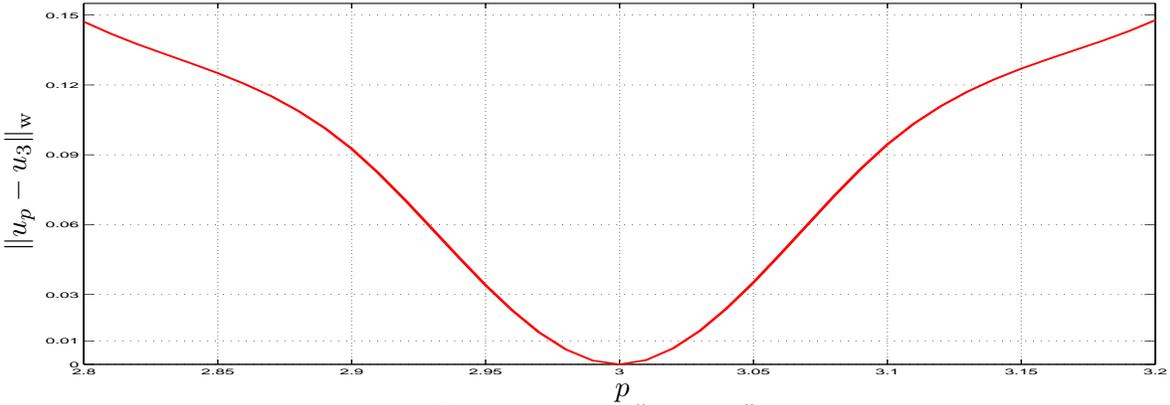,height=15cm,width=5cm, angle=270}
 \caption{$p$ vs. $\|u_p-u_3\|_{\rm w}$.} \label{fig:rate1}
\end{figure}

In Figure~\ref{fig:rate1} we plot $p$ vs.  $\|u_p-u_3\|_{\rm w}$ and
observe that $u_p(x)$ converge to $u_3(x)$ in the weighted norm
$\norm{{\cdot}}_{\rm w}$ as $p \to 3$.  In the numerical calculation
for Figure~\ref{fig:rate1}, our increment for $p$ is $0.01$.

\begin{rem}
For the operators $L_+L_-$ and $L_-L_+$, and in general 4-th order operators,
it seems difficult to exclude the possibility that $\mu=1$ is an eigenvalue.
Consider the following example.
Let $\tilde H := (\Lp)^2$ with $p=\sqrt 8 -1$. Note
$-1$ is an eigenvalue of $\Lp$ when $p=\sqrt 8 -1$. Hence $1$ is an
eigenvalue of $\tilde H $, at the endpoint of its continuous spectrum.

It would be interesting to prove the above convergence analytically
and characterize the leading order behavior near $p=3$ as we did in
Theorem~\ref{th2-2}.
\end{rem}

\section{Excited states with angular momenta}\label{sec4}

In this section we consider excited states with angular momenta in
$\R^n$, $n \ge 2$. Let $k= [n/2]$, the largest integer no larger
than $n/2$.  For $x =(x_1, \ldots,x_n) \in \R^n$, use polar
coordinates $r_j$ and $\th_j$ for each pair $x_{2j-1}$ and $x_{2j}$,
$j=1,\ldots,k$. P. L. Lions \cite{Lions} considers solutions of the
form
\[
Q(x) = \phi(r_1,r_2,\ldots r_k,x_{n}) \, e^{i (m_1 \th_1 + \cdots m_k
\th_k)}, \qquad m_j \in \mathbb{Z}.
\]
The dependence of $\phi$ in $x_n$ is dropped if $n$ is even. He proves
the existence of energy minimizing solutions in each such class.

For the simplest case $n=2$, $Q(x)=\phi(r) \, e^{i m \th}$ and, by
\eqref{Q.eq}, $\phi=\phi(r)$ satisfies
\begin{equation} \label{ESeq-1}
- \phi '' - \frac {1}r \phi' + \frac {m^2}{r^2} \phi + \phi -
|\phi|^{p-1} \phi = 0, \quad (r>0).
\end{equation}
The natural boundary conditions are
\begin{equation} \label{ESeq-2}
\lim_{r \to 0} r^{-m} \phi(r) = \al, \quad
\lim _{r \to 0} r^{-m+1} \phi'(r) =m \al, \quad
\lim _{r \to \infty} \phi(r) = 0,
\end{equation}
for some $\al \ge 0$. One can choose $\phi(r)$ real-valued.  It is
shown by Iaia-Warchall \cite{IW} that \eqref{ESeq-1}--\eqref{ESeq-2}
has countably infinite many solutions, denoted by $\phi_{m,k,p}(r)$,
each has exactly $k$ positive zeros. They correspond to
``$m$-equivariant'' nonlinear bound states of the form
\begin{equation}
Q_{m,k,p} =  \phi_{m,k,p}(r) \, e ^{i m \theta},
\qquad (k=0,1,2,\ldots).
\end{equation}
Note that $Q_{m,k,p}$ are radial if and only if $m=0$, and the ground
state $Q = Q_{0,0,p}$ is considered in the previous sections.
The uniqueness question of $\phi_{m,k,p}(r)$ is not addressed in
\cite{IW}. It is proved for the case $k=0$ in \cite{Mi1}.

Mizumachi \cite{Mi1}--\cite{Mi4} considered the stability problem for
these solutions.  He showed that
\begin{enumerate}
\item Under $m$-equivariant perturbations of the form $\e(r) e^{i m
\th}$, $Q_{m,0,p}$ are stable for $1<p<3$ and unstable for $p>3$;
\item Under general perturbations, $Q_{m,k,p}$ are unstable for $p>3$
for any $k$;

\item linear (spectral) instability implies nonlinear instability; (it
can be also obtained by extending the results in \cite{Shatah-Strauss} to
higher dimensions using the method of \cite{Almeida-Guo});

\item For fixed $p>1$, if $m > M(p)$ is sufficiently large,
$Q_{m,k,p}$ are linearly unstable and its linearized operator has a
positive eigenvalue.
\end{enumerate}

We are most interested in the last result.  Intuitively, for
$1<p<\infty$, $Q_{m,k,p}$ should be unstable for all $(m,k) \not =
(0,0)$ since they are excited states. Can this be observed
numerically? It turns out to be true for $p$ away from $1$, but false
for $p$ near 1.

In the following, we first describe our numerical methods for $k=0$
and next discuss their relations.  We will discuss our numerical
results in the end. We only compute $m=1,2$ but the same methods work
for other $m$.

\begin{rem}
Our numerical methods do not apply when $k>0$. Indeed, for $m \ge 0$
and $k>0$, the radial functions $\phi_{m,k,p}(r)$ are sign-changing
and cannot be numerically calculated using the method described in the
Appendix. In fact, it is an open question whether they are unique.
Assuming the uniqueness, one needs to develop a new algorithm to
compute them before one can compute the spectra of $\L$ for
$Q_{m,k,p}$.
\end{rem}

\subsection{Numerical algorithms}
There are two steps in our numerical method: First, compute
$\phi_{m,0,p}(r)$. Second, compute the spectra of the discretized
linearized operator around $Q_{m,0,p}$. The second step is more
involved and we will present three algorithms.

\medskip

{\bf Step 1.}  Compute $\phi(r) = \phi_{m,0,p}(r)$.  It is energy
minimizing among all solutions of \eqref{ESeq-1}--\eqref{ESeq-2} for
fixed $m,p$, and it is positive for $r>0$. Since our algorithm in the
Appendix is applicable to all positive (ground state) solutions, we
can use it to calculate the discretized vector of $\phi(r)$ with a
small change of the code.

\medskip

{\bf Step 2.} Compute the spectra of the discretized linearized
operator.  The linearized operator $\L$ has a slightly different form
than \eqref{eq1-7} because $Q = \phi_{m,0,p}(r)e^{i m \th}$ is no
longer real. With the same ansatz \eqref{eq1-5}--\eqref{eq1-6}, the
linearized operator $\L$ has the form
\begin{equation}
  \L h = i \bke{\Delta h - h + \frac {p+1}2 |Q|^{p-1} h + \frac {p-1} 2
    |Q|^{p-3} Q^2 \bar h }.
\end{equation}
We have developed three algorithms for computing the spectrum of $\L$.

\medskip

{\bf Algorithm~1.}  Write $Q=\phi(r)\, e^{i m \th}= \phi(r) \cos (m
\th) + i \phi(r) \sin(m \th)$. In vector form with $\L$ acting on
$[\Re h, \Im h]^\top$, we have
\begin{equation}\label{neweq4-5}
\L \sim
\begin{bmatrix} 0 & - \Delta + 1 \\ \Delta - 1 & 0 \end{bmatrix}
+
|\phi(r)|^{p-1}
\begin{bmatrix}
- (p-1) \cos \sin    & - \cos ^2   - p \sin^2  \\
p \cos ^2 + \sin^2  &  (p-1) \cos \sin
\end{bmatrix}(m \th).
\end{equation}

It is convenient to use polar coordinates to discretize the operator.
We use a two dimensional mesh,
\begin{equation}\label{mesh2d}
\text{2d mesh: }\quad r = 0: \de_r: r_{\max}, \quad \th = 0:
\de_\th: 2\pi.
\end{equation}
The discretized matrix has size $N T$ by $NT$ with $N=
r_{\max}/\de_r$ and $T= 2\pi/ \de_\th$. We use zero boundary
condition with $r_{\max}=15$, $\de_r = 0.04$, and $T=160$.

Although the matrix operator \eqref{neweq4-5} is slightly more
complicated than \eqref{eq1-7} and the mesh is 2-dimensional, the same
numerical routine can be applied to compute the spectrum of the
discretized matrix of \eqref{neweq4-5}. The only difference is that
the matrix size is much larger.

\medskip

{\bf Algorithm~2.} By restricting the problem to some invariant
subspaces of $\L$, as we did for the computation of figures
\ref{fig:cL2d-m1}--\ref{fig:cL2d-m2}, we will reduce the problem to
1-dimension.

Observe that functions of the form $a(r) e^{i j \th}$ with a fixed
integer $j$ are not preserved by $\L$ unless $j=m$, but the following
$L^2$-subspaces are invariant under $\L$:
\[
X_k = X_k^{(m)}= \bket{h(r,\th): h = a(r) e^{i (m+k) \th} + b(r)
e^{i(m-k)\th}}, \quad 0 \le j\in \mathbb{N}.
\]
If $k=0$, we drop $b(r)$ and $X_0 = \bket{h(r,\th) :h= a(r) e^{i m\th}
}$.  We will compute the spectra of $\L$ limited to each subspace $X_k$.
Define
\[
V = \frac {p-1}2 \phi^{p-1} , \quad H_k = - \Delta_r +1 +
\frac{(m+k)^2}{r^2} - \frac {p+1}2 \phi^{p-1}.
\]

For $k=0$, with $a=a_1+ia_2$ and $a_1,a_2\in \R$, we have
\[
\L [a(r)e^{i m \th}] = -i\bkt{H_0 a - V \bar a}e^{i m \th}
= [H_0 (a_2-ia_1)  + V(a_2+ia_1)]e^{i m \th}.
\]
Thus, acting on $[a_1, a_2]^\top$, $\L|_{X_0}$ has the matrix form
\[
L_{X_0} = \mat{ 0 & H_0 + V \\ - H_0 + V & 0}.
\]

For $k>0$, with $a=a_1+ia_2$, $b=b_1+ib_2$ and $a_1,a_2,b_1,b_2 \in
\R$, we have
\begin{align*}
&\L [a(r)e^{i (m +k)\th}+ b(r) e^{i(m-k)\th}] \\
&=  [H_{k} (a_2-ia_1) + V(b_2+ib_1)]e^{i (m +k)\th}
+[H_{-k} (b_2-ib_1) + V(a_2+ia_1)]e^{i (m-k) \th}.
\end{align*}
Thus, acting on $[a_1, a_2, b_1, b_2]^\top$, $\L|_{X_k}$ has the
matrix form
\[
L_{X_k}= \mat{
0    & H_{k} & 0 & V \\
-H_{k} & 0   & V & 0 \\
0    & V   & 0 & H_{-k} \\
V    & 0   & - H_{-k} & 0}.
\]
To discretize the operator, we use the one-dimensional mesh
\begin{equation}\label{mesh1d}
\text{1d mesh: } \quad r = 0: \de_r: r_{\max}, \quad N=
r_{\max}/\de_r.
\end{equation}
The matrix corresponding to $X_0$ has size $2N$ by $2N$. The matrix
for $X_k$ with $k>0$ has size $4N$ by $4N$.  We use zero boundary
condition with $r_{\max}=30$ and $\de_r = 0.01$.

Counting multiplicity, the eigenvalues of $\L$ is the union of
eigenvalues of $\L |_{X_k}$ with $k=0,1,2,\ldots$.

\medskip

{\bf Algorithm~3.}  Instead of the form \eqref{eq1-5}, include the
phase $e^{im\th}$ in the linearization: $\psi = (\phi + h)e^{i m
\th+it}$. Then the linearized operator acting on $[\Re h, \Im
h]^\top$ is
\[
\L' = \mat{
-2m/r^2 \pd_\th   &   -\Delta + 1 + m^2/r^2 - \phi^{p-1}  \\
-(-\Delta + 1 + m^2/r^2 - p \phi^{p-1})   &   -2m/r^2 \pd_\th}
\]
which is invariant on subspaces $Z_k =\bket{ [a_1(r), a_2(r)]^\top
e^{i k \th}}$ with integers $k$.  We have
\[
\L' \mat{a_1(r)\\ a_2(r)}e^{i k \th}
= e^{i k \th} L_{m,k} \mat{a_1(r)\\ a_2(r)},
\]
where
\[
L_{m,k} := \mat{
-\frac{2imk}{r^2}  & -\Delta_r + 1 + \frac{m^2+k^2}{r^2} - \phi^{p-1}  \\
-(-\Delta_r + 1 + \frac{m^2+k^2}{r^2} - p \phi^{p-1}) & -\frac{2imk}{r^2} }
\]
acting on radial functions.  We use the same one-dimensional mesh
\eqref{mesh1d} as in Algorithm~2.  For every $k$, the matrix size is
$2N $ by $2N$. We then compute the spectra of $L_{m,k}$ for each
$k$.

Counting multiplicity, the eigenvalues of $\L$ is the union of
eigenvalues of $L_{m,k}$ with $k=0,\pm 1,\pm 2,\ldots$.

\subsection{Properties of these algorithms}
We now address the relation between these algorithms.  First note that
$X_k$ is essentially the sum of $Z_k$ and $Z_{-k}$.  Let us make it
more precise and suppose $k>0$. The case $k=0$ is easier.  A function
$h=(a_1+ia_2)(r)e^{i (m+k)\th} + (b_1+ib_2)(r)e^{i(m-k)\th}$ in $X_k
\subset L^2(\R^2)$ can be identified with $[a_1,a_2,b_1,b_2] \in
\tilde X_k= L^2_{rad}(\R^2;\R^4)$. The space $\tilde X_k$ is a
subspace of $L^2_{rad}(\R^2;\C^4)$ on which we compute the spectrum.
The function $h$ can be also identified with
\[
\mat{a_1(r)\\ a_2(r)} e^{i k \th}+ \mat{b_1(r)\\ b_2(r)} e^{-i k \th},
\]
the collection of which form a subspace of $ Z_k \oplus Z_{-k}$ with
real components.

\medskip

{\bf Nullspace of $\L$}. \quad The nullspace of $\L$ gives a good test
of the correctness of our numerical results.  For $k=0$, the
$0$-eigenfunction $iQ$ of $\L$ corresponds to $[0 , \phi]^\top e^{im
\th}$ in $X_0$ and $[0 , \phi]^\top$ in $Z_0$. The generalized
eigenfunction $Q_1 = \frac{2}{p-1}Q + x\cdot \nabla Q$ corresponds to
$[Q_1,0]^\top e^{i m\th}$ in $X_0$ and $[\frac
2{p-1}\phi+r\phi',0]^\top$ in $Z_0$. Since $X_0 \subset L^2(\R^2,\C)$,
they also provide two (generalized) eigenvectors for Algorithm~1.

For $k=\pm 1$, the $0$-eigenfunctions
\[
2Q_{x_1} = 2(\phi' \cos \th - i \psi \sin \th)e^{im\th} = (\phi'-\psi)e^{i
(m+1) \th} + (\phi'+\psi)e^{i (m-1) \th},
\]
\[
2Q_{x_2} = 2(\phi' \sin \th + i \psi \cos \th)e^{im\th}= i(-\phi'+\psi)e^{i
(m+1) \th} + i (\phi'+\psi)e^{i (m-1) \th},
\]
where $\psi = m\phi/r$, belong to $X_1$, and correspond to
$0$-eigenvectors $[\phi'-\psi, 0, \phi'+\psi, 0]^\top$ and
$[0,-\phi'+\psi, 0,\phi'+\psi]^\top$ of $L_{X_1}$.  For Algorithm~3,
they correspond to the following vectors in $Z_1 \oplus Z_{-1}$,
\[
2\mat{ \phi' \cos \th \\ -\psi \sin \th} = W_+ e^{ i \th}+ W_- e^{ -i
\th}, \quad
2\mat{ \phi' \sin \th \\  \psi \cos \th} =  -i W_+ e^{ i \th}+i W_- e^{ -i
\th}
\]
where
\[
W_\pm = \mat{ \phi' \\ \pm i\psi}, \quad L_{m,\pm 1} W_\pm = \mat{0 \\ 0}.
\]
Thus $W_+e^{ i \th}$ is a $0$-eigenvector of $\L'$ in $Z_{1}$, and
$W_-e^{- i \th}$ is a $0$-eigenvector of $\L'$ in $Z_{-1}$.

The generalized eigenfunctions
\[
ix_1Q = i r \phi \cos \th e^{i m \th} = i r \phi e^{i (m+1)\th} + i r
\phi e^{i (m-1)\th}
\]
\[
ix_2Q = i r \phi \sin \th e^{i m \th} = r \phi e^{i (m+1)\th} - r \phi
e^{i (m-1)\th}
\]
also lie in $X_1$ and correspond to generalized $0$-eigenvectors
$[0,r\phi,0,r\phi]^\top$ and $[r\phi,0,-r\phi,0]^\top$ of $L_{X_1}$.
For Algorithm~3, they correspond to $[0, r\phi \cos \th]^\top$ and
$[0, r\phi \sin \th]^\top$ in $Z_1 \oplus Z_{-1}$. By the same
consideration as for $Q_{x_1}$ and $Q_{x_2}$, their span over $\C$
is the same as the span of $[0, r\phi]^\top e^{i \th} \in Z_1$ and
$[0, r\phi]^\top e^{-i \th} \in Z_{-1}$.  One can check that
\begin{equation}\label{eq4-test}
L_{m,\pm 1} \mat{0 \\ r\phi} = -2 \mat{ \phi' \\\pm i\psi}.
\end{equation}

Thus, the multiplicity of $0$-eigenvalue in each of $X_0$, $Z_{-1}$,
$Z_0$ and $Z_1$ is at least 2. The multiplicity of $0$-eigenvalue on
$X_1$ is at least 4.

\medskip

{\bf Symmetry of spectra}. \quad
If
\[
L_{m,k} \mat{A \\ B} = \la \mat{A \\ B}
\]
then
\[
L_{m,-k} \mat{\bar A \\ \bar B} = \bar \la \mat{\bar A \\ \bar B},
\quad
L_{m,-k} \mat{ A \\ - B} = - \la \mat{ A \\ - B},
\quad
L_{m,k} \mat{\bar A \\ -\bar B} = -\bar \la \mat{\bar A \\ -\bar B}.
\]

In particular, if $\la \in \si(L_{m,k})$, then $-\bar \la\in
\si(L_{m,k})$, and $\bar \la, -\la \in \si(L_{m,-k})$. Thus
$\si(L_{m,k})$ itself is symmetric w.r.t.~the imaginary axis, and
$\si(L_{m,k})$ and $\si(L_{m,-k})$ are symmetric w.r.t.~the real axis.

Similarly, one can show that the spectra of $L_{X_k}$ are symmetric
with respect to both real and imaginary axes.

\medskip

{\bf Equivalence of Algorithms~2 and 3.} \quad In Algorithm~2, for
$k>0$, we can write
\[
L_{X_k} = \mat{ H_kJ & VU \\ VU & H_{-k}J}
\]
where
\[
J= \mat{0 & 1 \\ -1 & 0}, \quad
U= \mat{0 & 1 \\ 1 & 0}, \quad
I=\mat{ 1 & 0 \\  0 & 1}.
\]
Let
\[
M=\mat{ I &
 -J \\ I & J}, \quad M^{-1} = \frac 12 \mat{ I & I \\ J & - J}, \quad
P=\mat{
1 & 0 & 0 & 0 \\
0 & 0 & 1 & 0 \\
0 & 1 & 0 & 0 \\
0 & 0 & 0 & 1},
\quad P^{-1}=P.
\]
Noting $JU=-UJ$, we have
\[
M^{-1} L_{X_k} M  = \mat{ \al J + cU & -\beta \\ -\beta & \al J + cU}
=\mat{
0 & \al + c & \beta & 0 \\
-\al+c & 0 & 0 & \beta  \\
-\beta & 0 & 0 & \al+c \\
0 & -\beta & -\al+c & 0}
\]
where
\[
\al = \frac 12 (H_k+H_{-k})= H_0 + \frac {k^2}{r^2}, \quad \beta = \frac
12(H_k-H_{-k})=\frac {2mk}{r^2}, \quad c=V.
\]
Let
\[
L':=P^{-1}M^{-1} L_{X_k} M P = \mat{
0 & \beta & \al + c &0\\
- \beta & 0 & 0 & \al+c  \\
- \al+c & 0 & 0 & \beta  \\
0 & -\al+c  & -\beta & 0 }.
\]

In Algorithm~3, $L_{m,k}$ acts on $[A(r),B(r)]^\top$. If we write
the enlarged matrix of $L_{m,k}$ acting on $[\Re A, \Im A, \Re B, \Im
B]^\top$, the matrix is exactly $L'$. The matrix for $L_{m,-k}$ will
be also $L'$ if it acts on $[\Re A, -\Im A, \Re B, -\Im B]^\top$.
This amounts to a choice of assigning $J$ or $-J$ to the
complexification of $i$.

More precisely, if $L_{m,k} u = \la u$ with $u = [A,B]^\top$, then
$L_{m,k} iu = \la iu$.  Write $A=A_1+iA_2$ and $B=B_1+iB_2$ and
suppose $k>0$. These two equations are equivalent to
\[
L' \mat{ A_1 \\ A_2 \\ B_1 \\ B_2} = \mat{ \Re \la A \\ \Im \la A \\
  \Re \la B \\ \Im \la B},
\qquad
L' \mat{ -A_2 \\ A_1 \\ -B_2 \\ B_1} = \mat{ -\Im \la A \\ \Re \la A \\
  -\Im \la B \\ \Re \la B}.
\]
Adding the second equation multiplied by $-i$ to the first equation,
we get
\[
L' w = \la w, \quad w = [ A ,-iA , B ,-iB]^\top .
\]
Taking conjugation we get $L' \bar w= \bar \la \bar w$.  Thus $\la$
and $\bar \la$ are eigenvalues of $L'$, and hence of $L_{X_k}$.
Since $L_{m,k} u = \la u$ iff $L_{m,-k}\bar u = \bar \la \bar u$,
eigenvalues of $L_{m,-k}$ also correspond to eigenvalues of
$L_{X_k}$.

\medskip

{\bf Counting eigenvalues}. \quad $L_{X_k}$ acts on $L^2_{rad}(\R^2,
\C^4)$ while $L_{m,\pm k}$ act on $L^2_{rad}(\R^2, \C^2)$.  The
eigenvalues of $L_{X_k}$ is the union of eigenvalues of $L_{m,k}$ and
$L_{m,-k}$. For any ball $B_R$ on the complex plane disjoint from the
continuous spectrum $\Sigma_c = \bket{ ir: r\in \R, |r|\ge 1}$,
\[
\#( \si(L_{X_k}) \cap B_R) = \# (\si(L_{m,k}) \cap B_R) + \#
(\si(L_{m,-k})) \cap B_R
\]
which is equal to $2  \#( \si(L_{m,k}) \cap B_R)$ if the center of $B_R$
is on the real axis.

\medskip

{\bf Numerical efficiency.} \quad Algorithm~1 is 2-dimensional, and
thus more expensive to compute and less accurate.  Both Algorithms~2
and 3 are one-dimensional and more accurate.

The benefit of Algorithm~3 than Algorithm~2 is that it further
decomposes the subspace of $L^2(\R^2,\C^4)$ corresponding to $X_k$
to two subspaces. Although its matrix size is only half that of
Algorithm~2, its components are complex and hence require more
storage space. Numerically these two algorithms are not very
different.

\subsection{Numerical results}
\begin{figure}[ht]
\psfrag{p}{$p=1.06$}
\epsfig{file=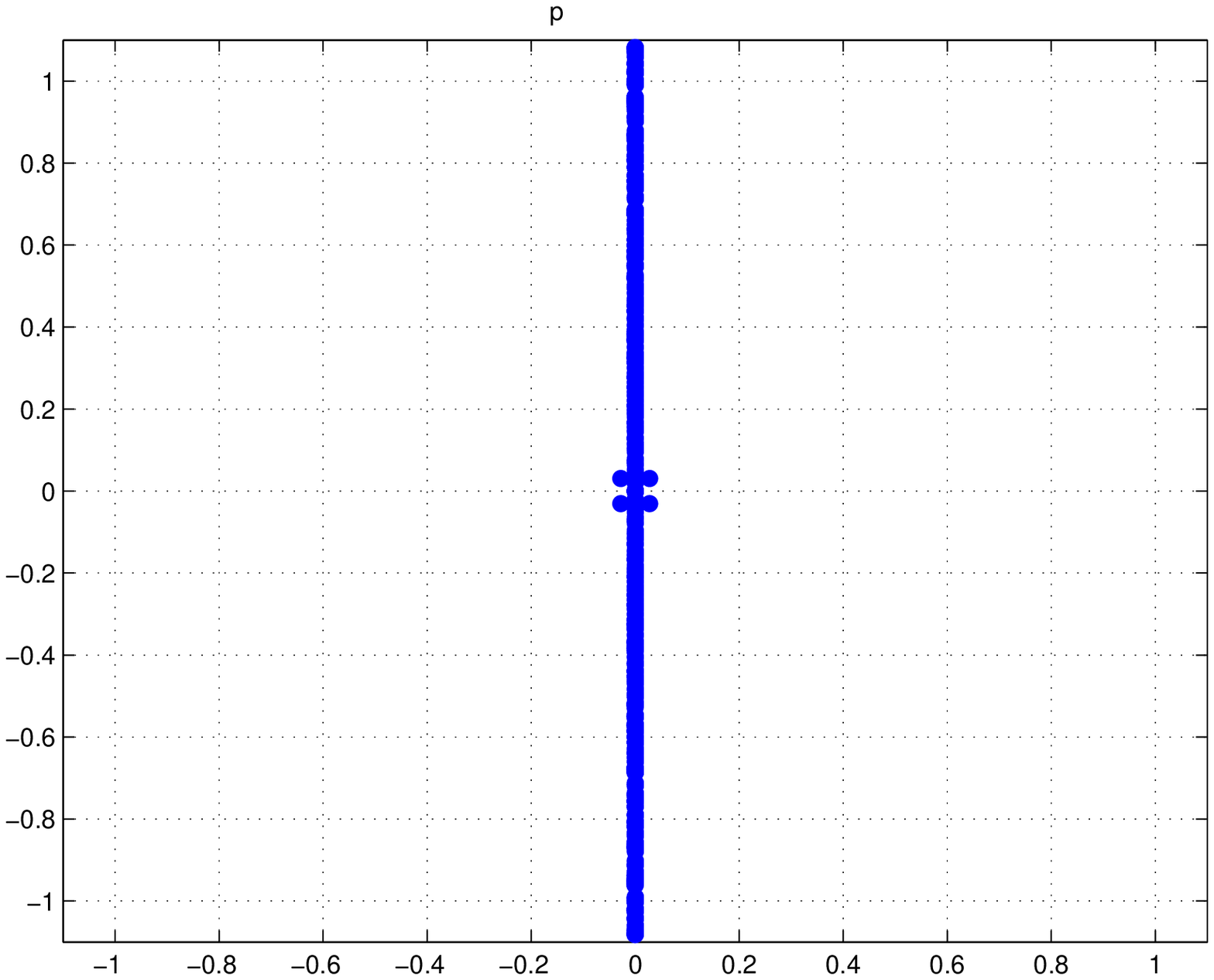,height=5cm,width=5cm}
\psfrag{p}{$p=1.1$}
\epsfig{file=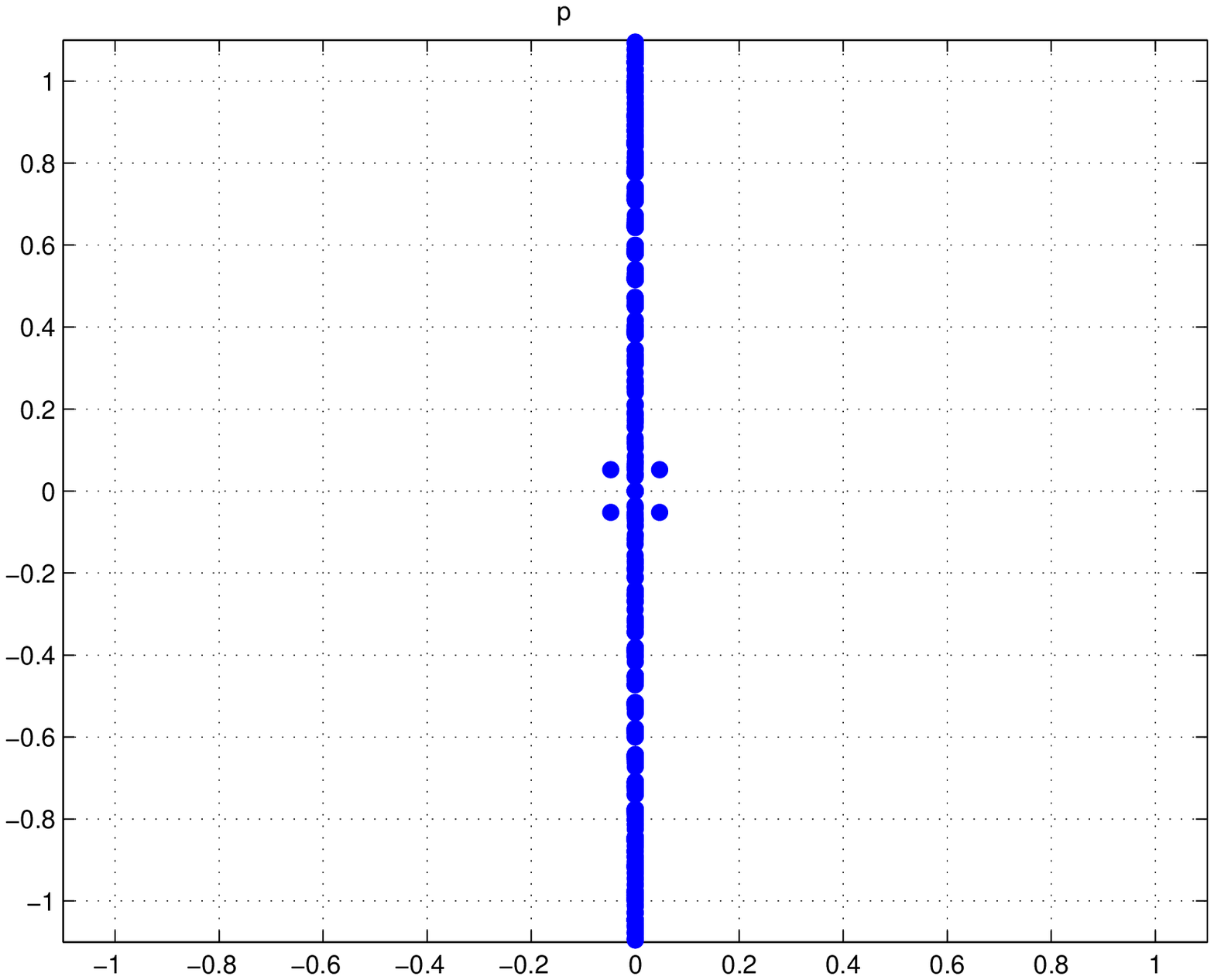,height=5cm,width=5cm}
\psfrag{p}{$p=1.4$}
\epsfig{file=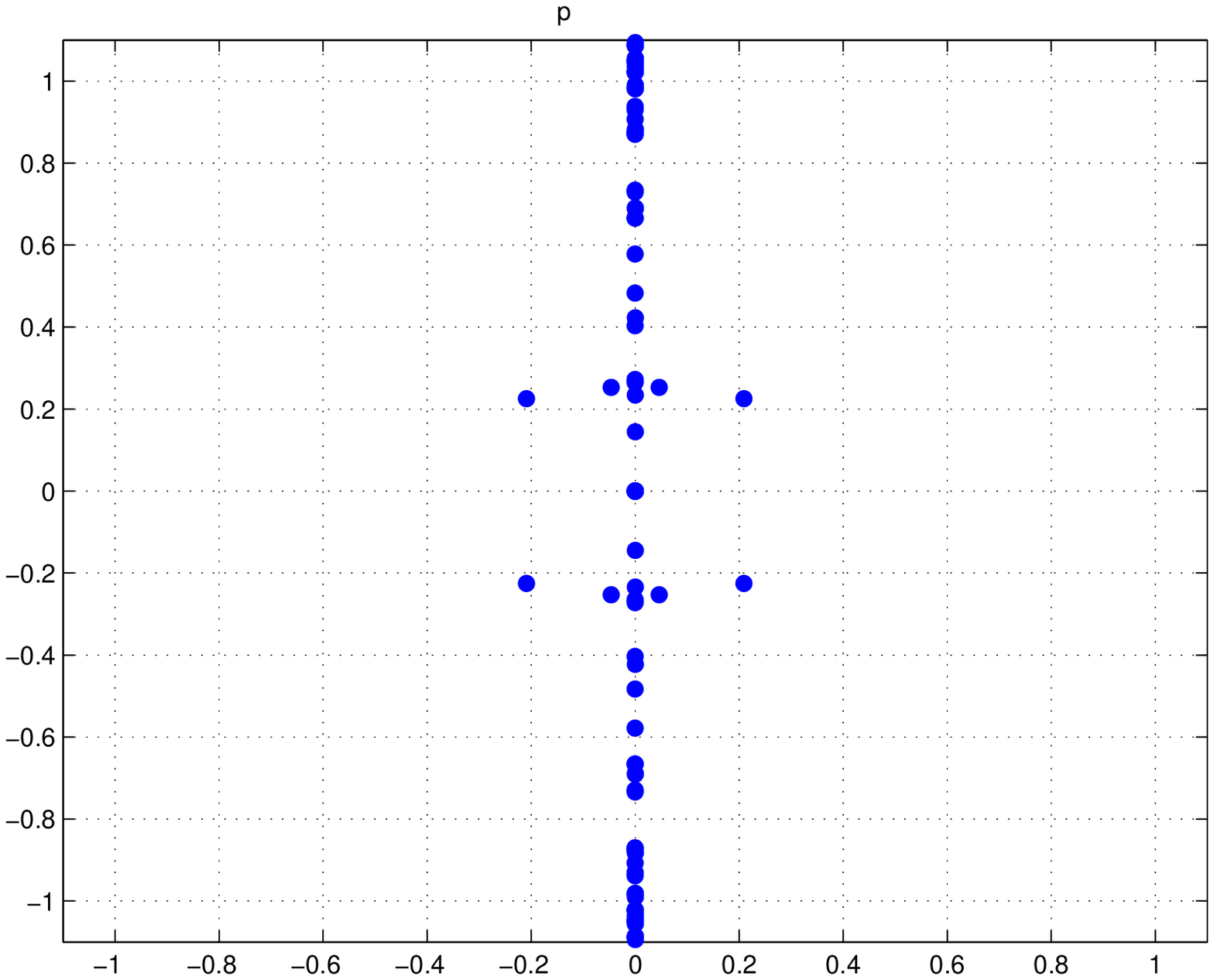,height=5cm,width=5cm} \vspace*{0.5cm}
\psfrag{p}{$p=1.55$}
\epsfig{file=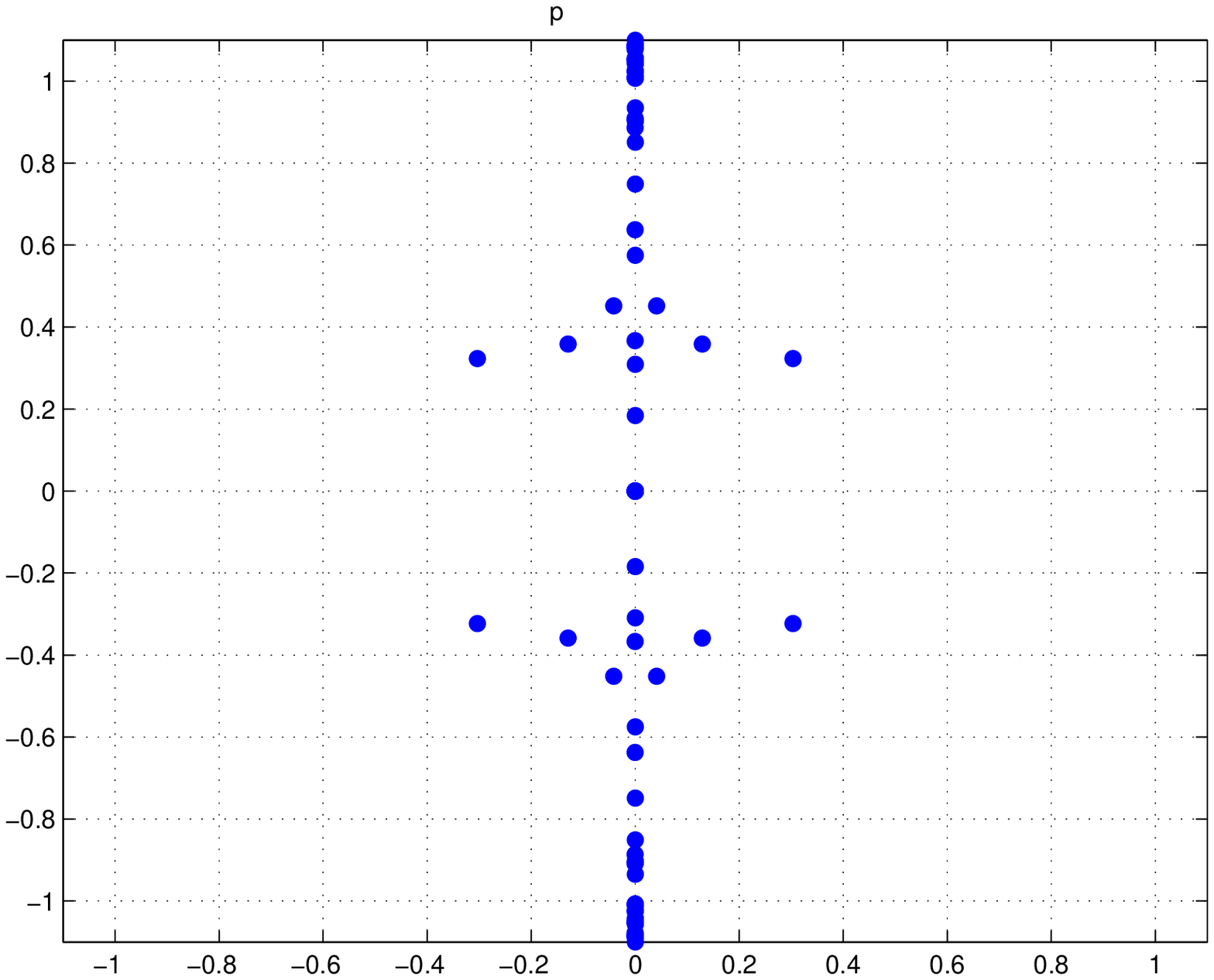,height=5cm,width=5cm}
\psfrag{p}{$p=1.8$}
\epsfig{file=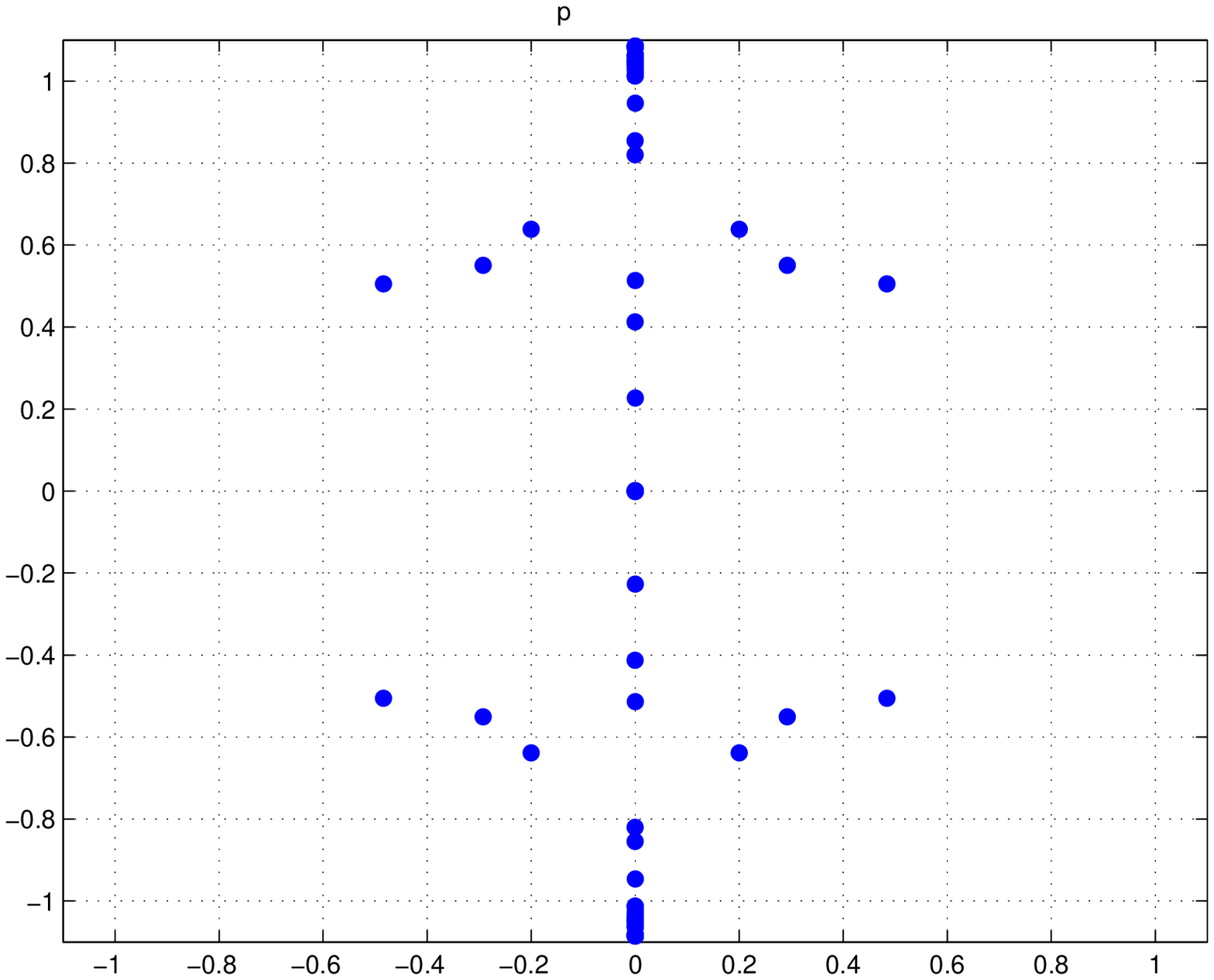,height=5cm,width=5cm}
\psfrag{p}{$p=2.3$}
\epsfig{file=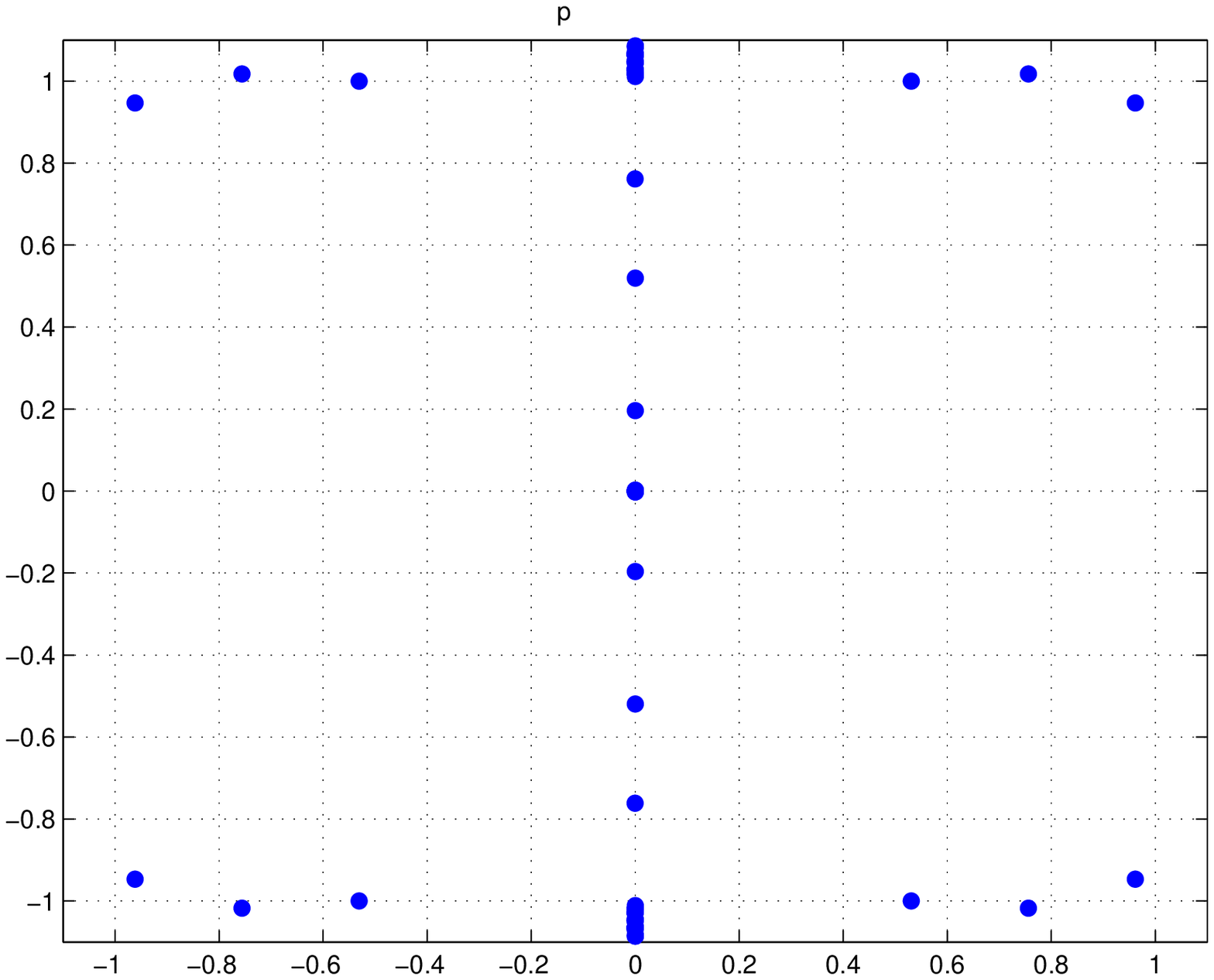,height=5cm,width=5cm} \vspace*{0.5cm}
\psfrag{p}{$p=3$} \epsfig{file=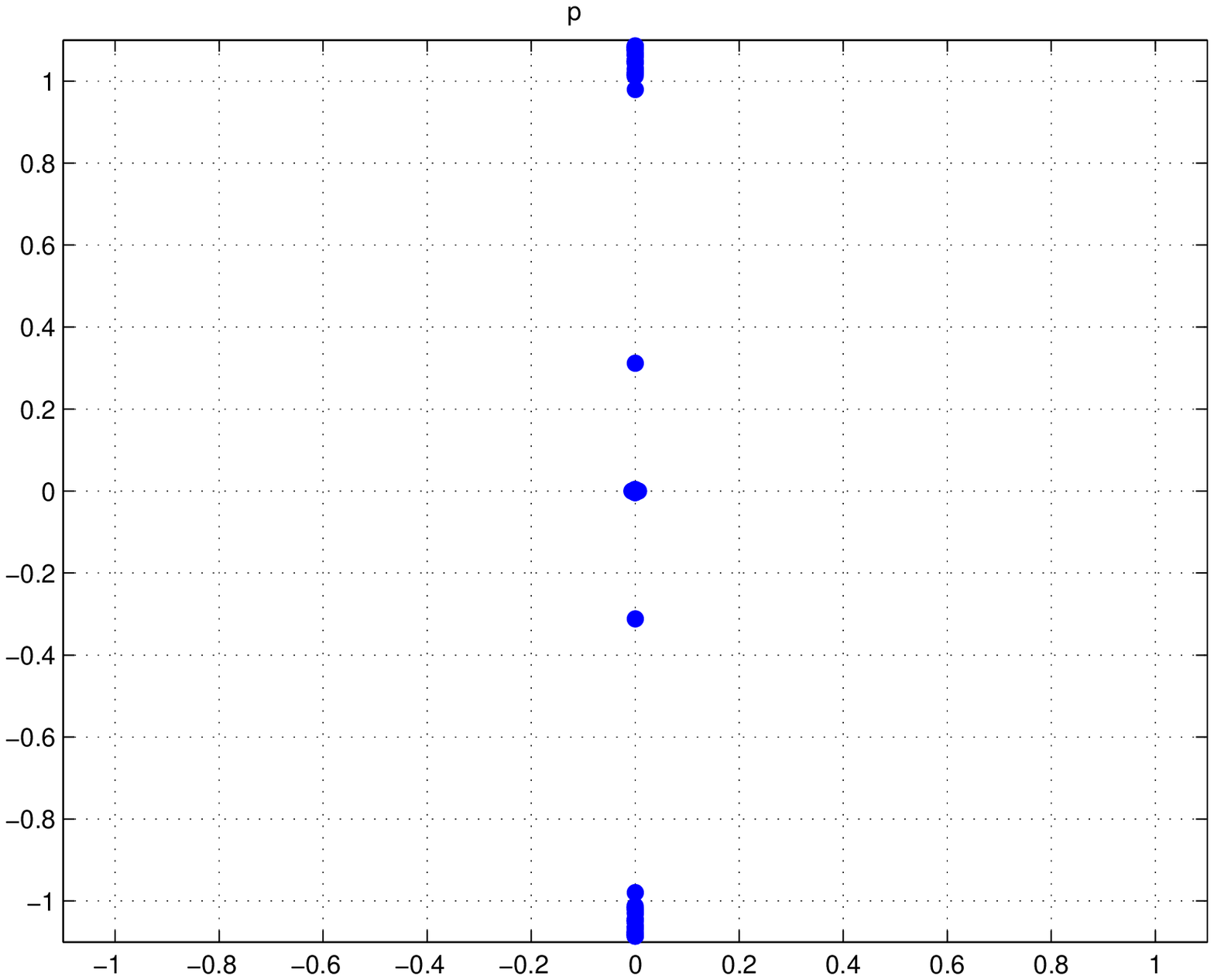,height=5cm,width=5cm}
\psfrag{p}{$p=3.1$}
\epsfig{file=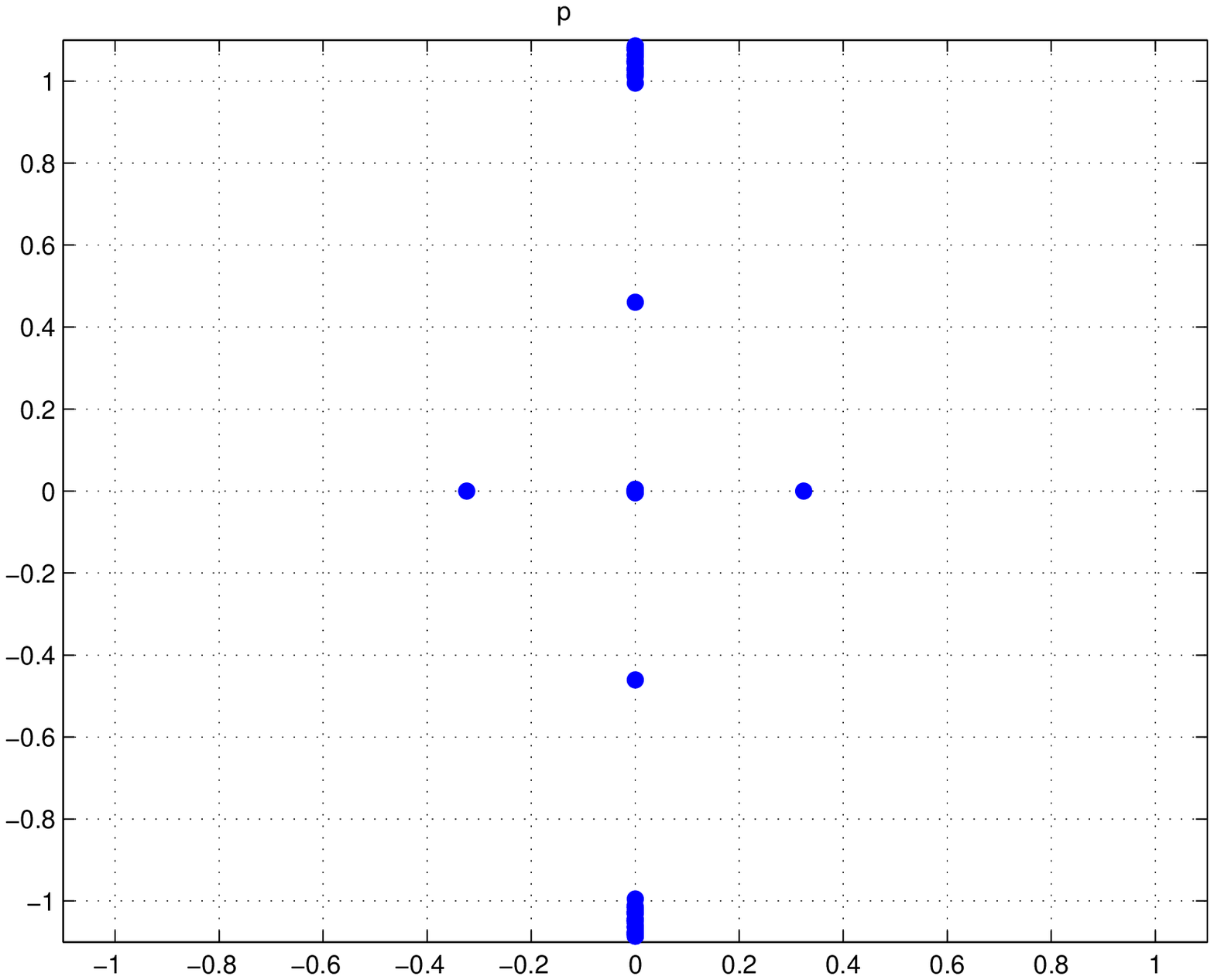,height=5cm,width=5cm}
\psfrag{p}{$p=3.5$}
\epsfig{file=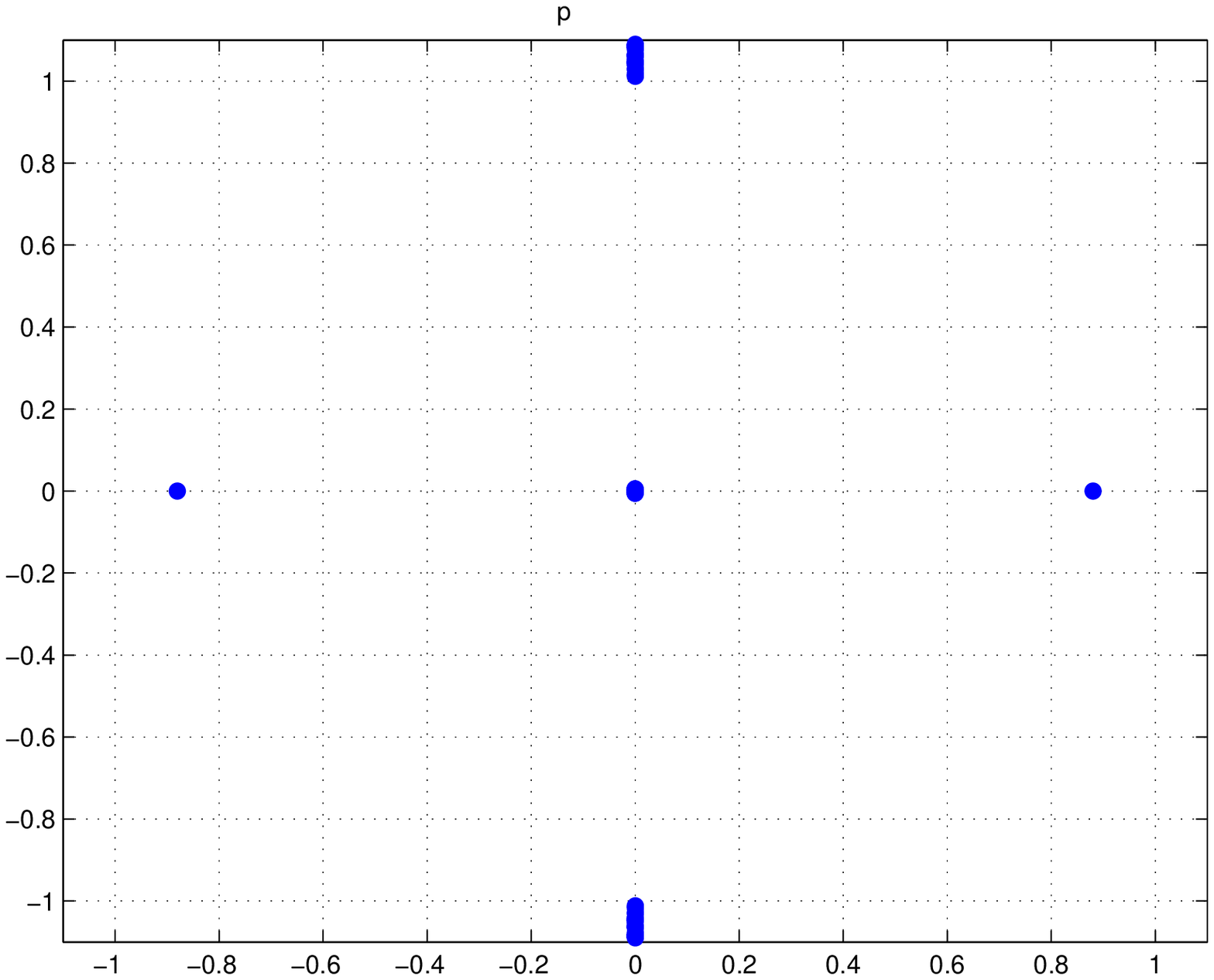,height=5cm,width=5cm}
 \caption{Spectra of $\L$ in $\R^2$ for $m=1$ and $0 \le k \le
   24$ as
$p = 1.06, 1.1, 1.4, 1.55, 1.8, 2.3, 3, 3.1, 3.5$ computed by
Algorithm~3.} \label{fig:d2m1}
\end{figure}

\begin{figure}[ht]
\psfrag{p}{$p=1.06$}
\epsfig{file=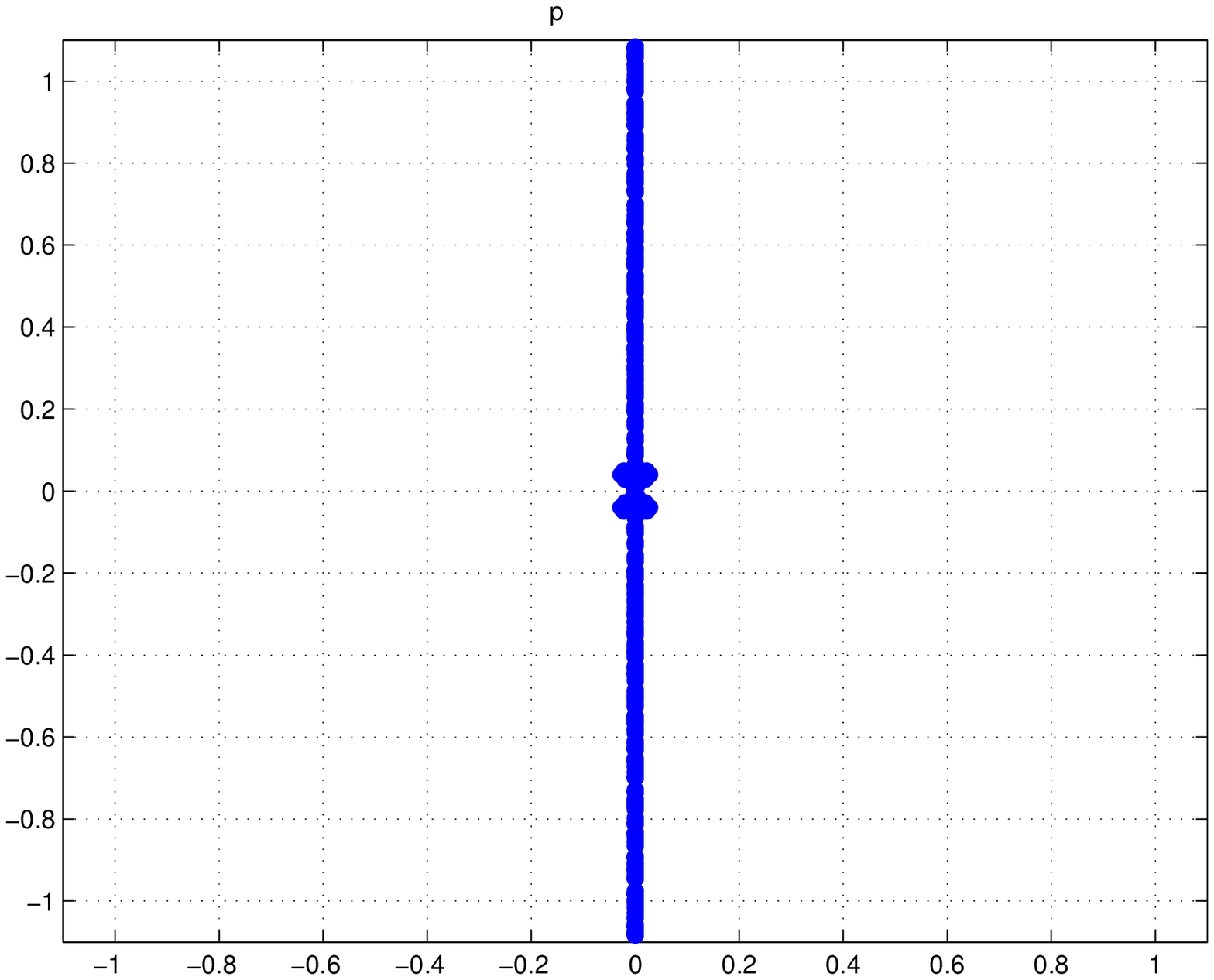,height=5cm,width=5cm}
\psfrag{p}{$p=1.1$}
\epsfig{file=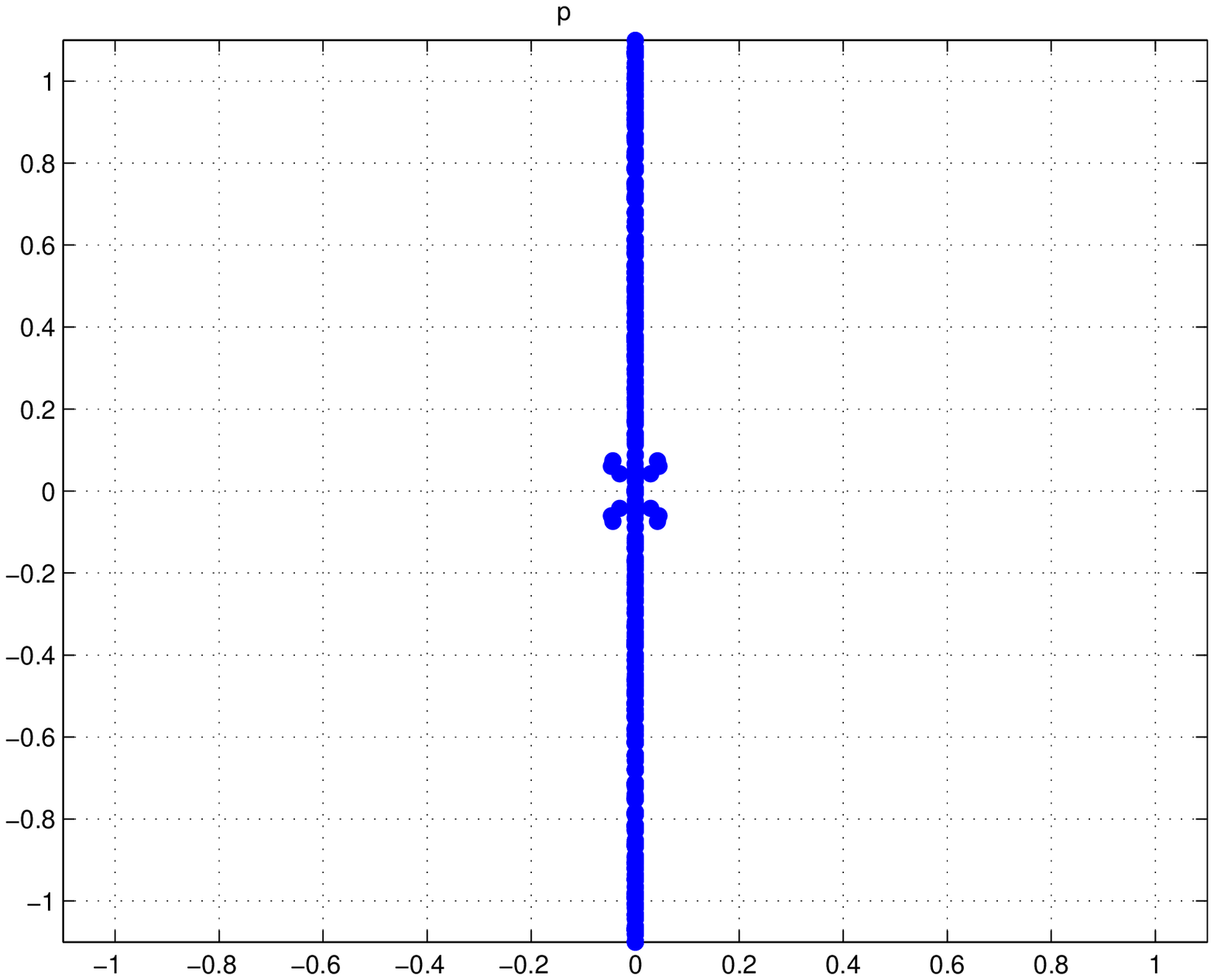,height=5cm,width=5cm}
\psfrag{p}{$p=1.3$}
\epsfig{file=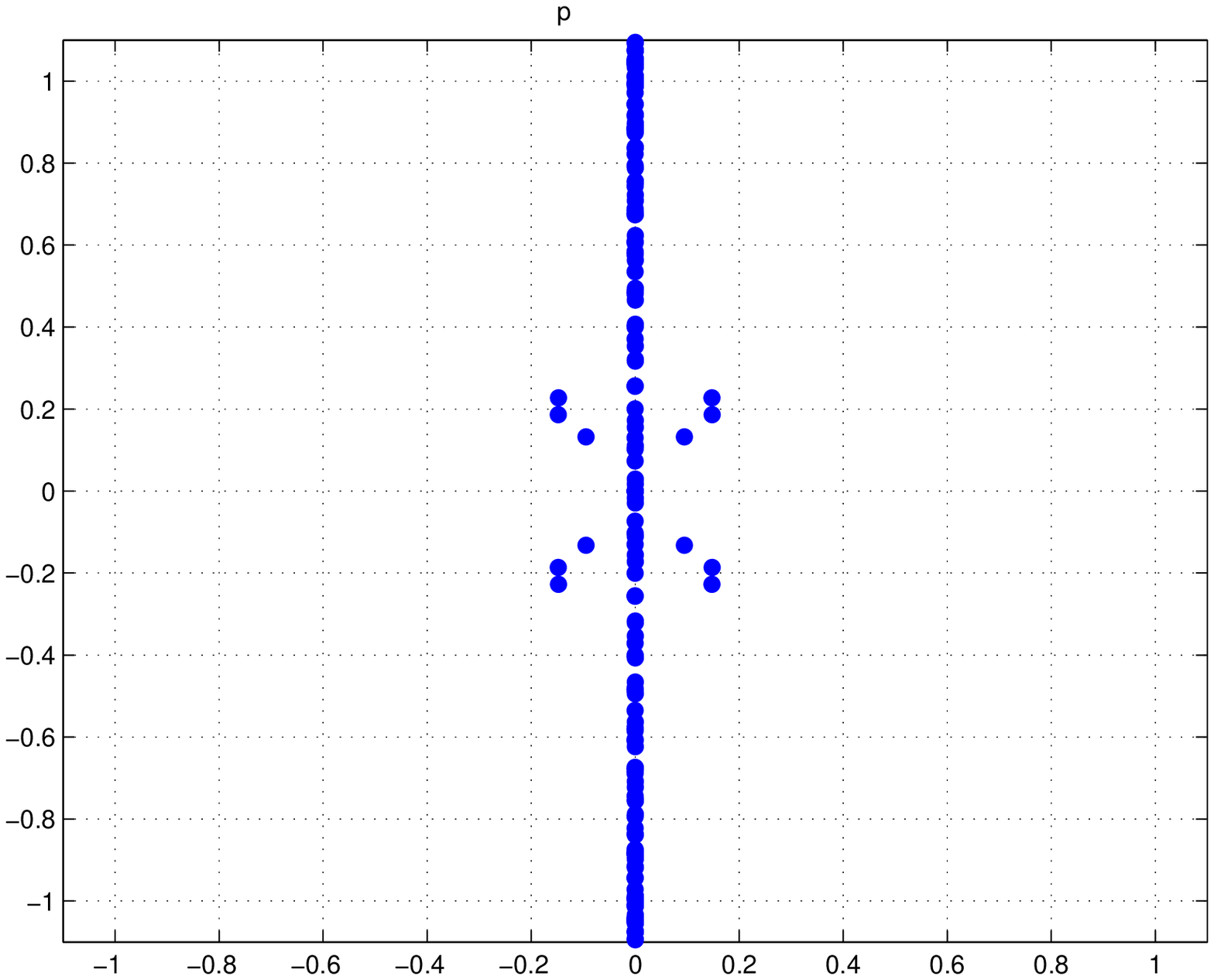,height=5cm,width=5cm} \vspace*{0.5cm}
\psfrag{p}{$p=1.4$}
\epsfig{file=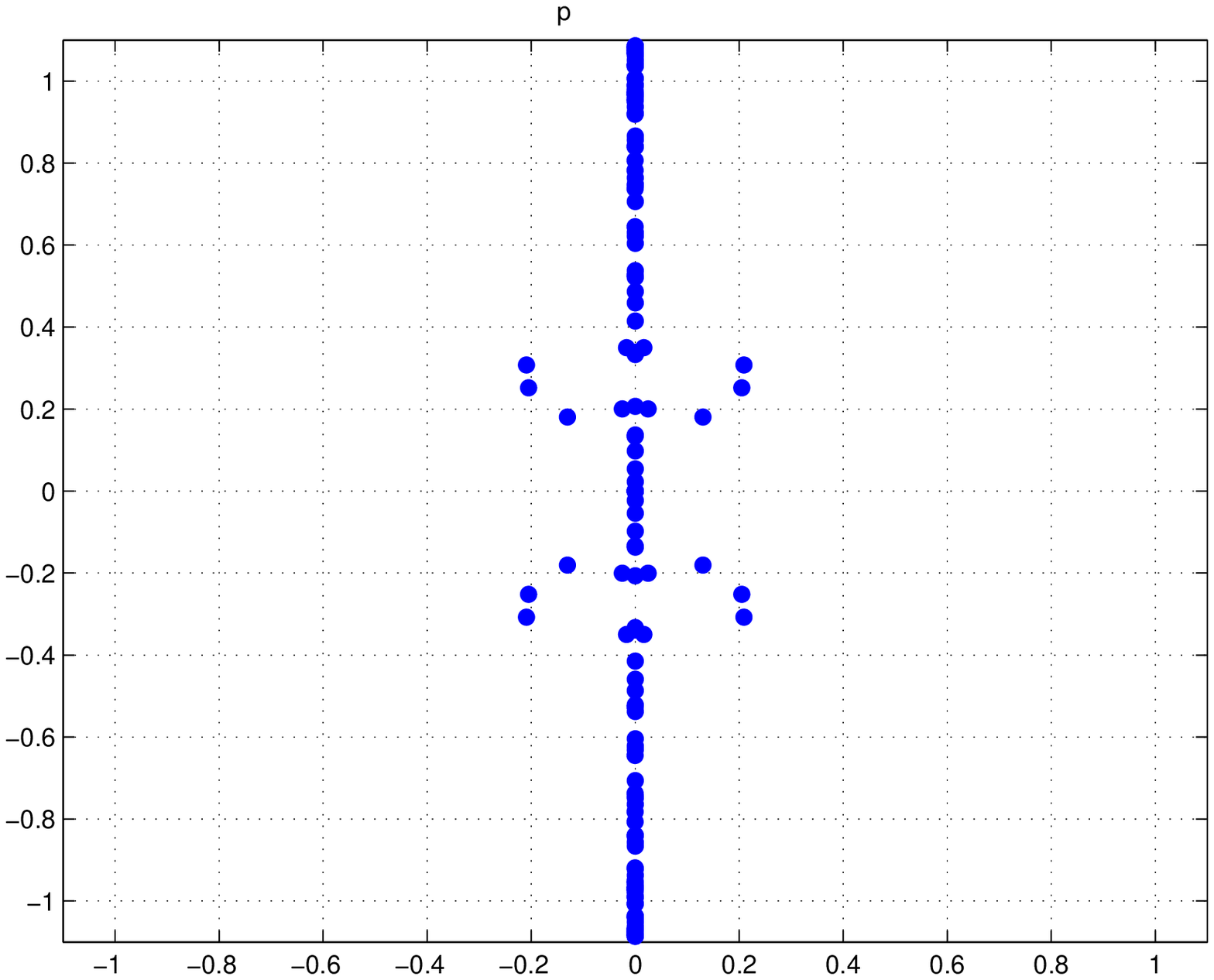,height=5cm,width=5cm}
\psfrag{p}{$p=1.8$}
\epsfig{file=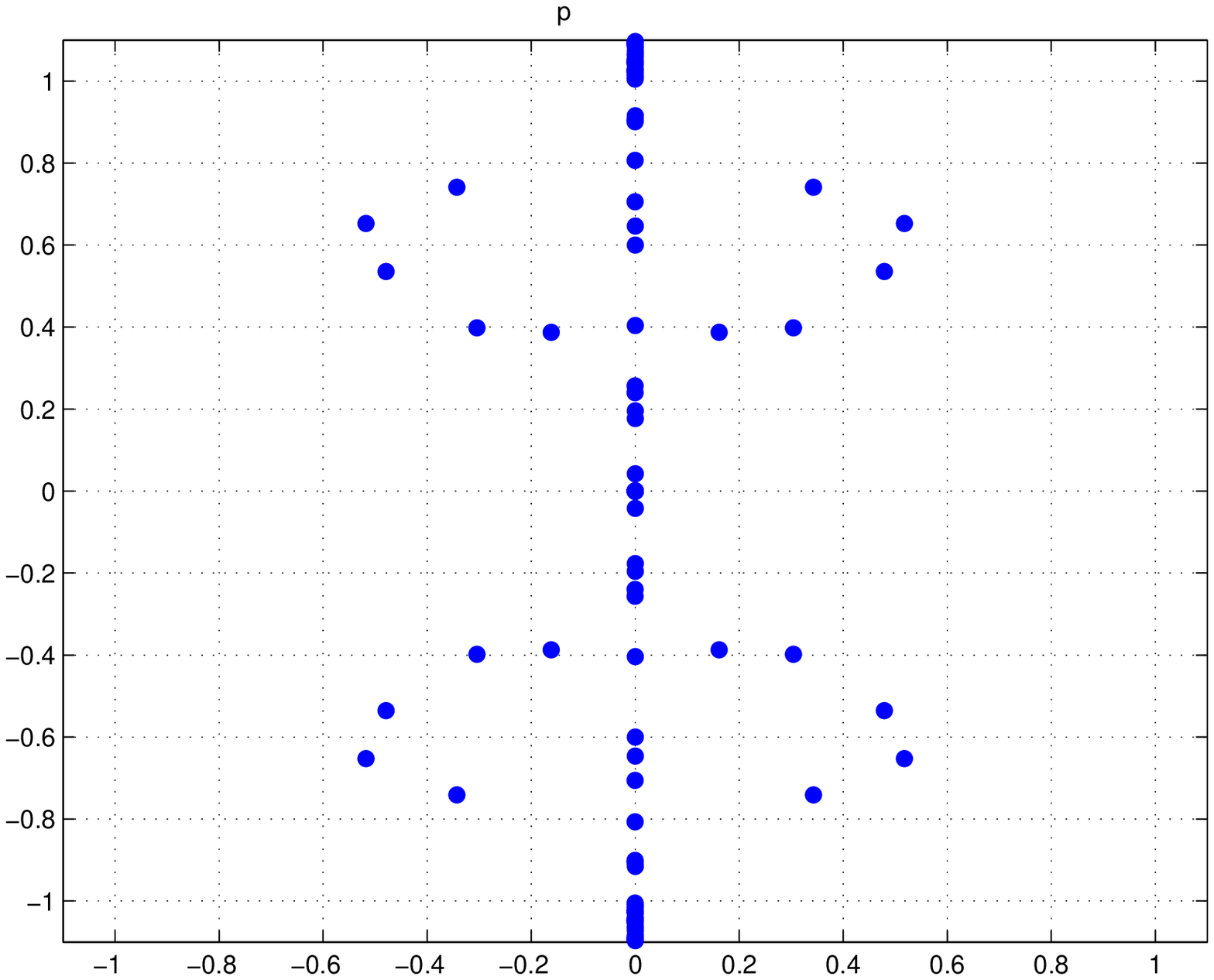,height=5cm,width=5cm}
\psfrag{p}{$p=2.2$}
\epsfig{file=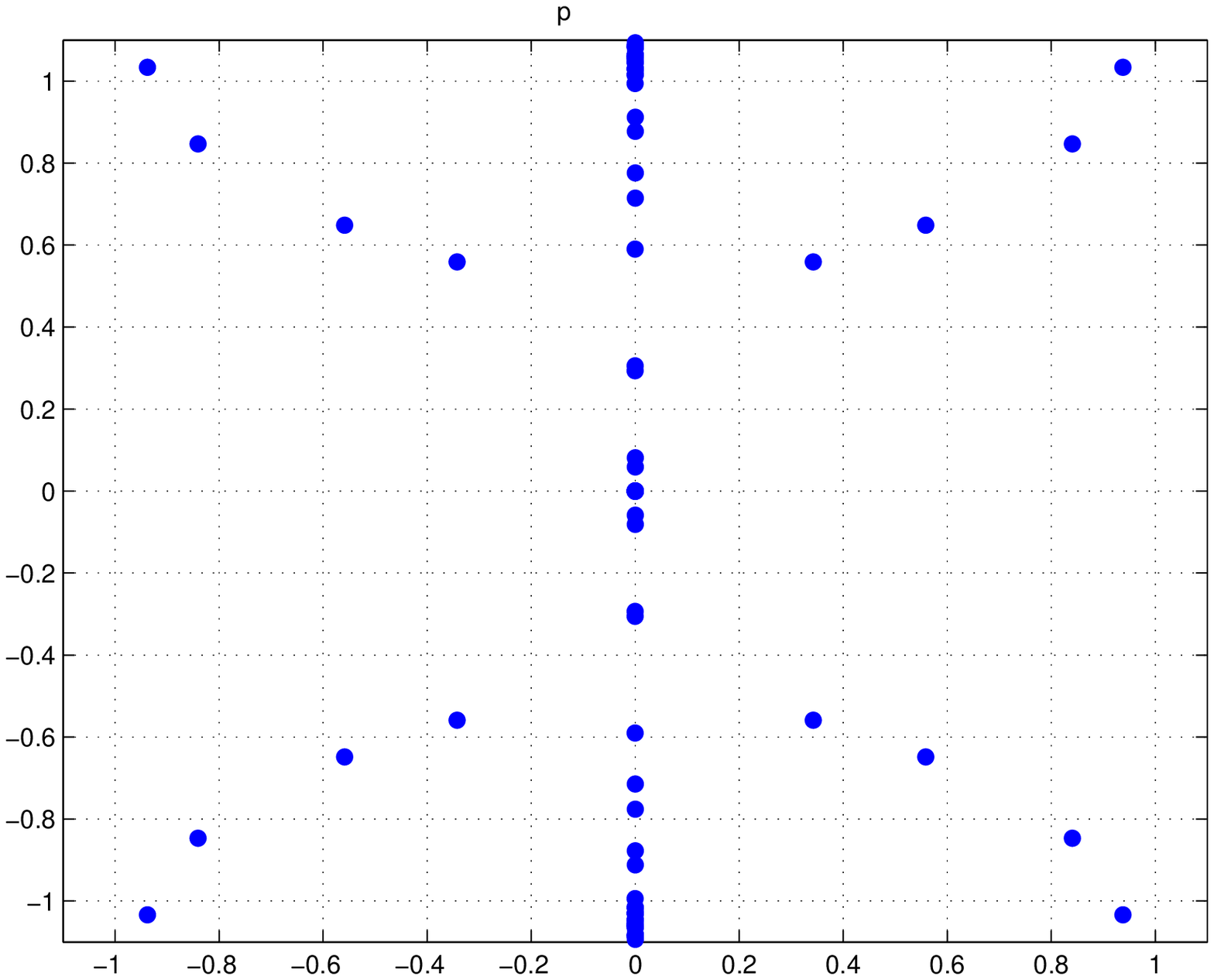,height=5cm,width=5cm} \vspace*{0.5cm}
\psfrag{p}{$p=2.6$}
\epsfig{file=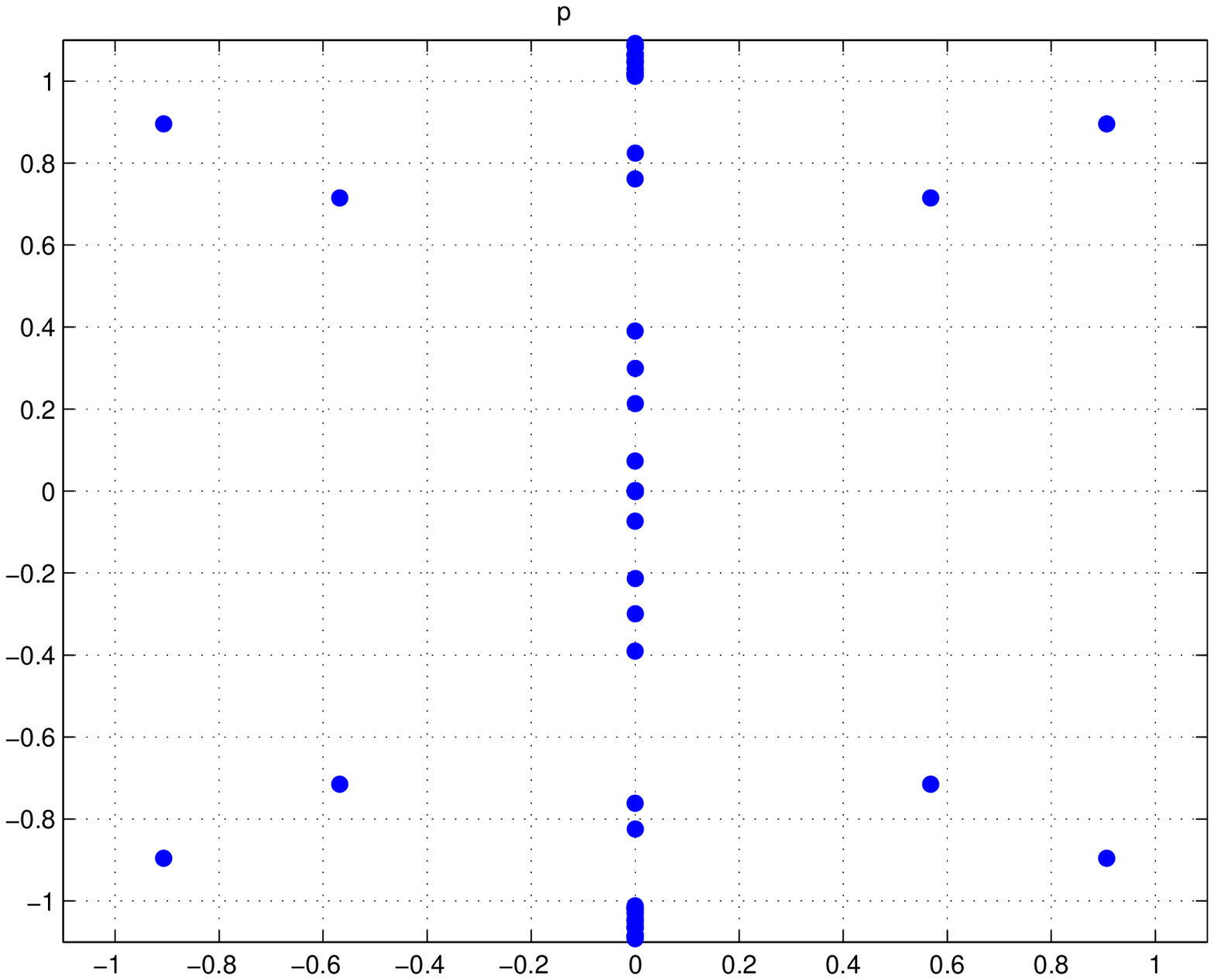,height=5cm,width=5cm} \psfrag{p}{$p=3$}
\epsfig{file=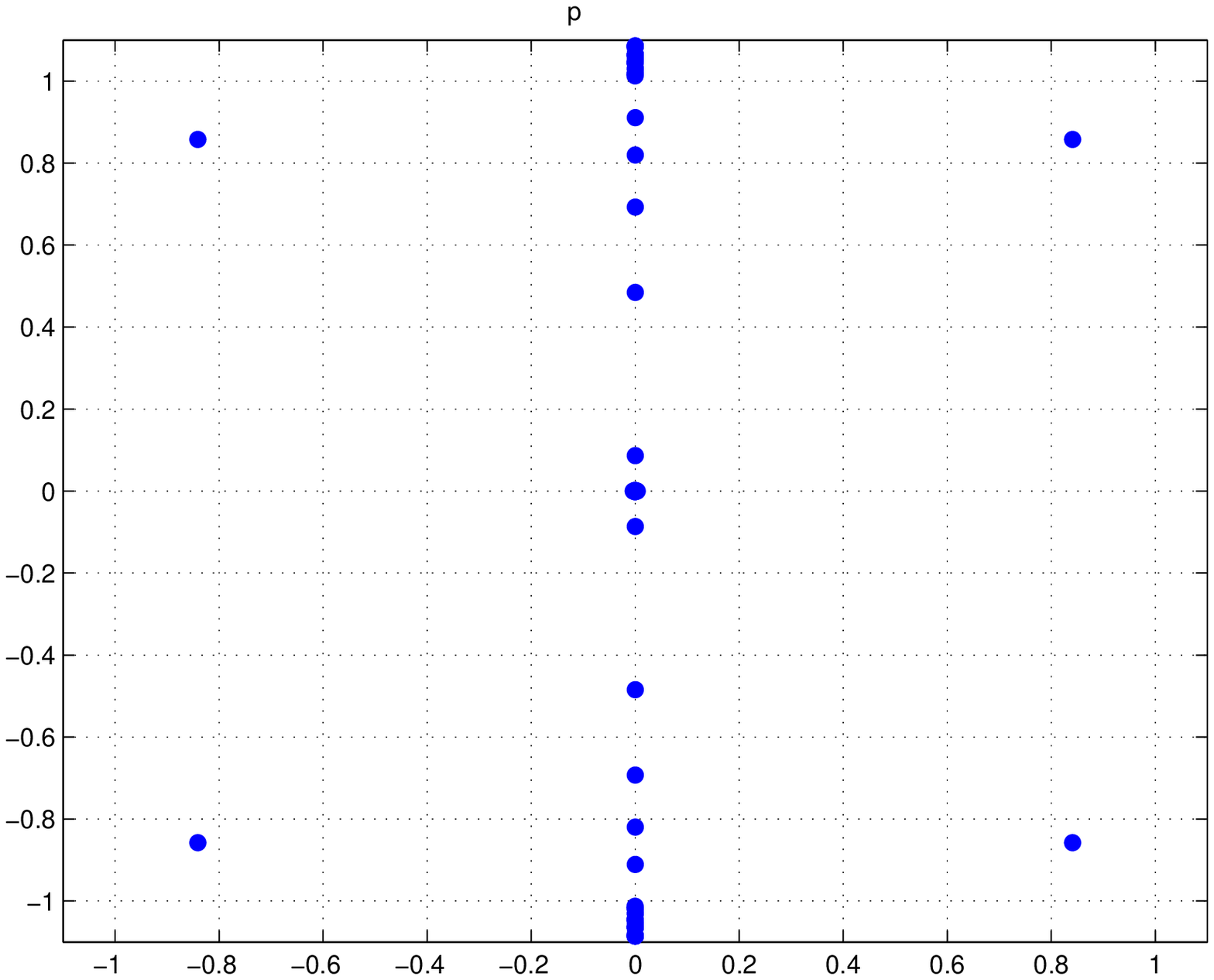,height=5cm,width=5cm}
\psfrag{p}{$p=3.5$}
\epsfig{file=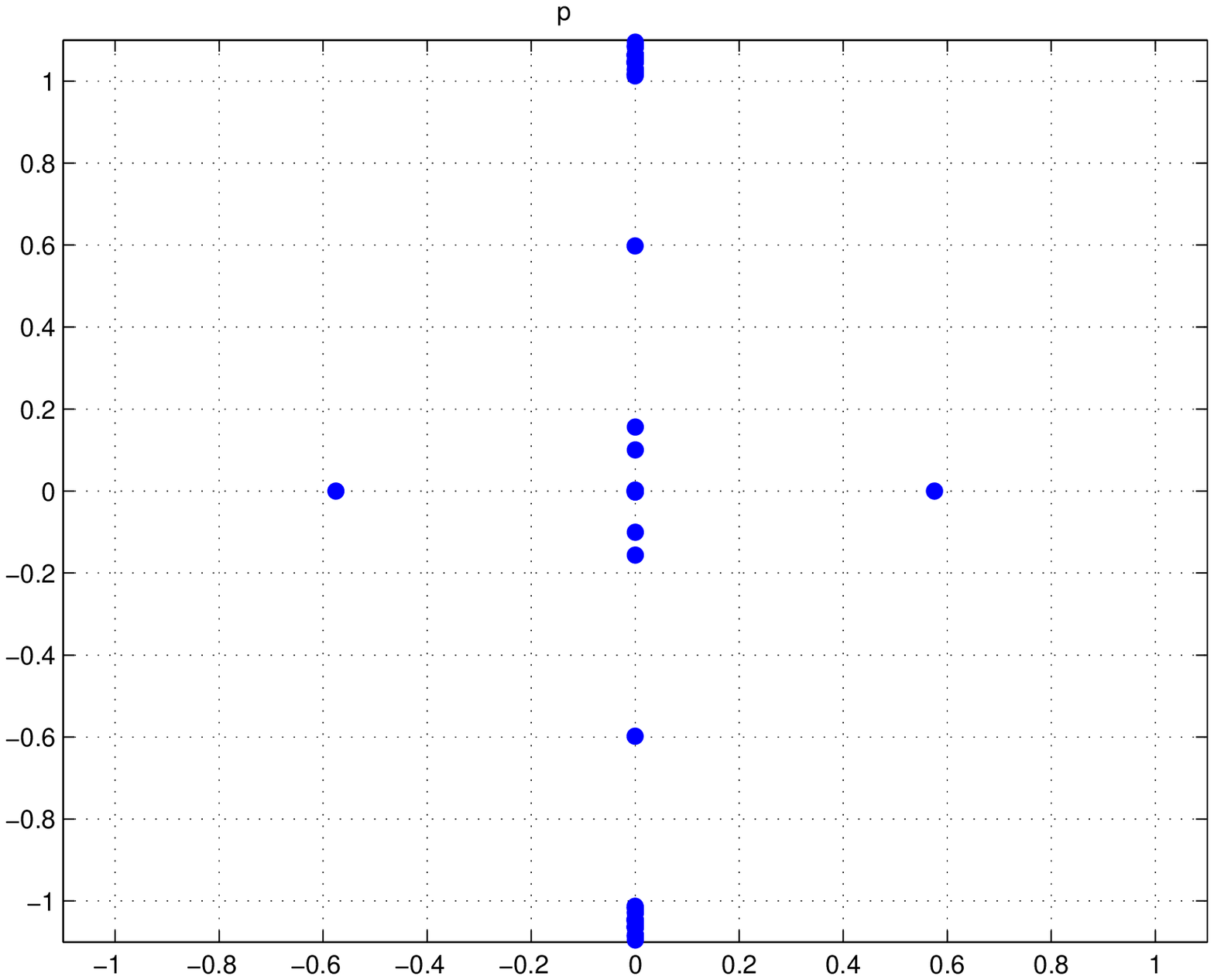,height=5cm,width=5cm}
 \caption{Spectra of $\L$ in $\R^2$ for $m=2$ and $0 \le k\le 28$ as
$p = 1.06, 1.1, 1.3, 1.4, 1.8, 2.2, 2.6, 3, 3.5$ computed by
Algorithm~3.} \label{fig:d2m2}
\end{figure}

\begin{figure}[ht]
\psfrag{p}{$p=1.6$}
\epsfig{file=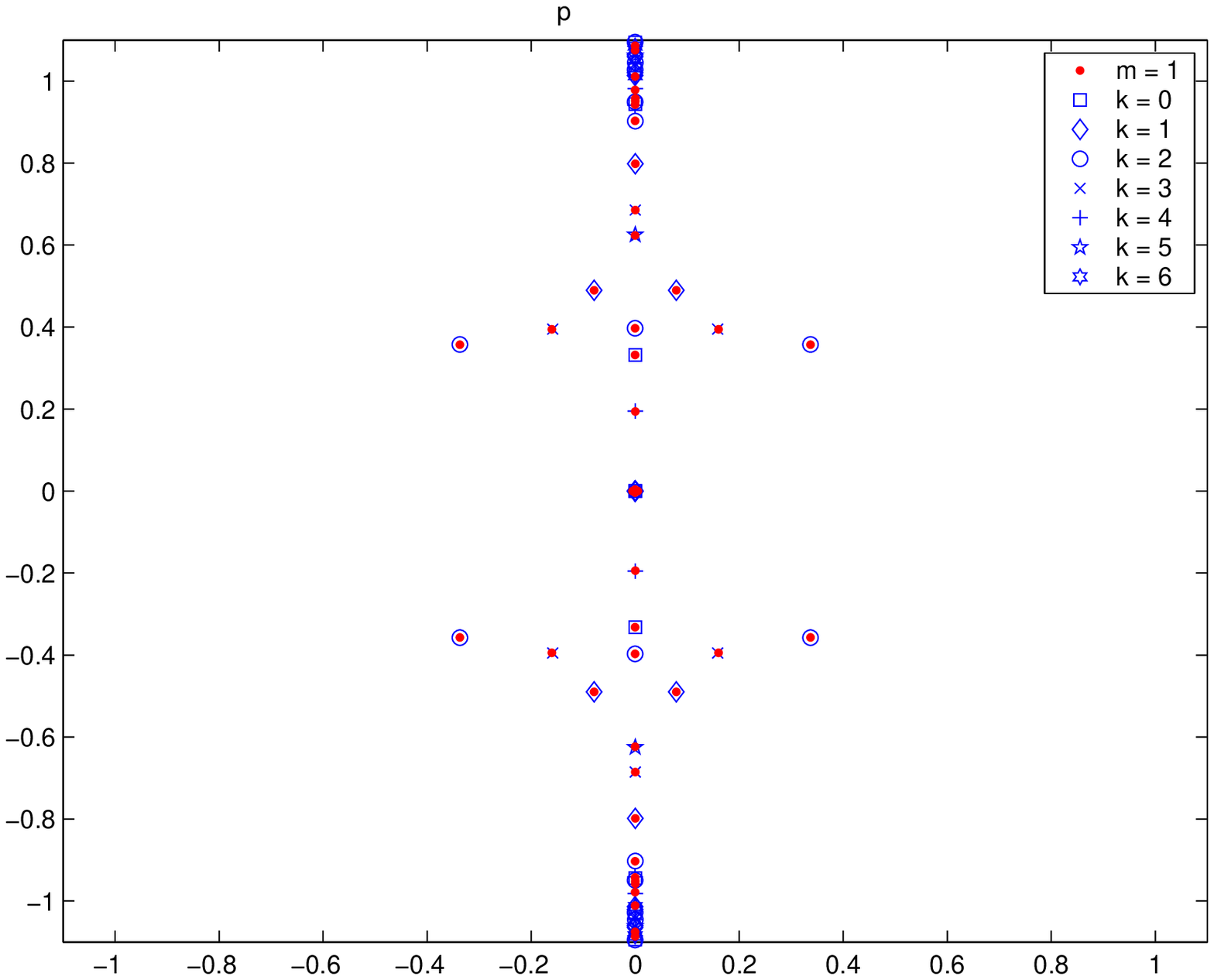,height=7.5cm,width=7.5cm}
\psfrag{p}{$p=2.1$}
\epsfig{file=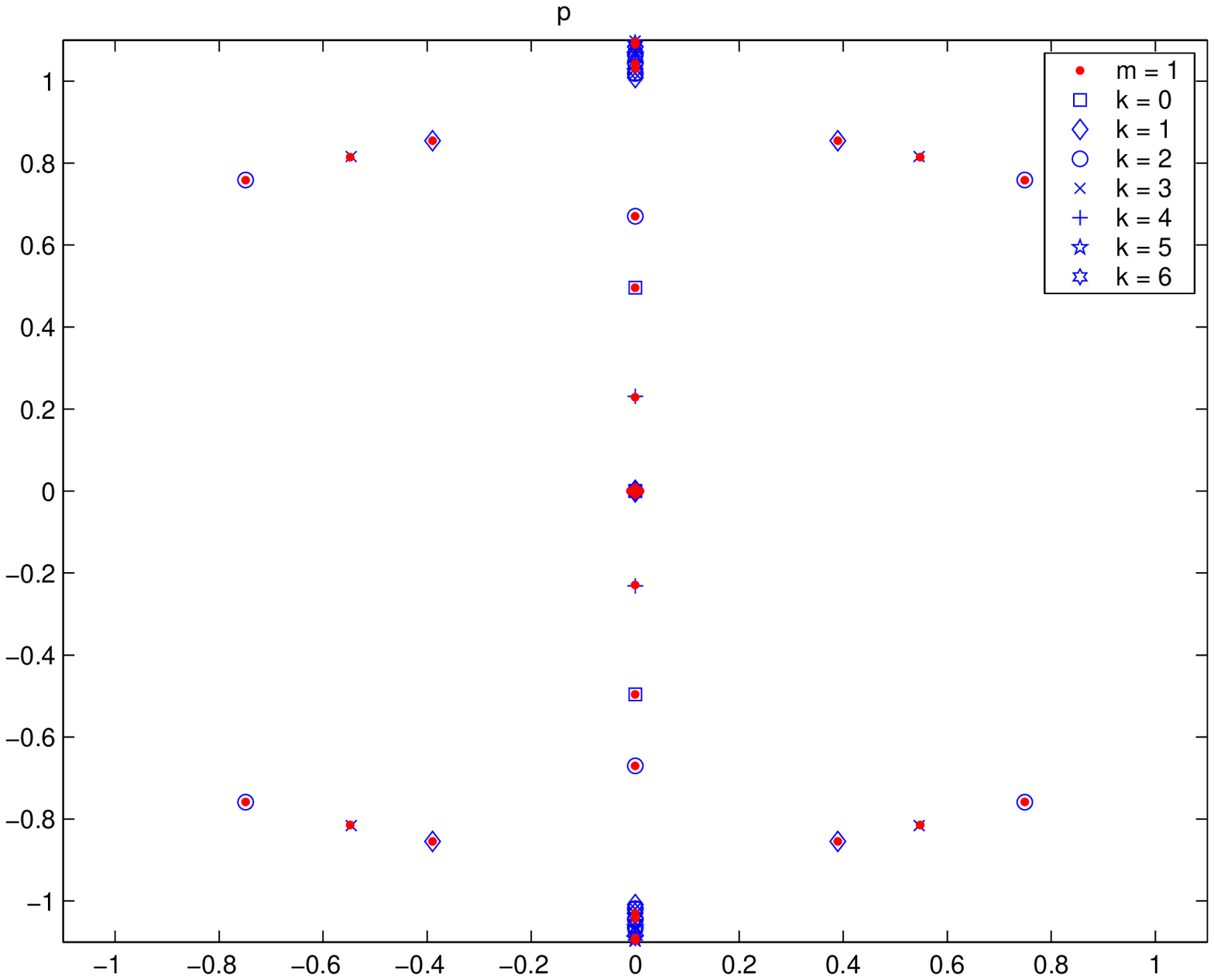,height=7.5cm,width=7.5cm}
\vspace*{0.5cm} \psfrag{p}{$p=3$}
\epsfig{file=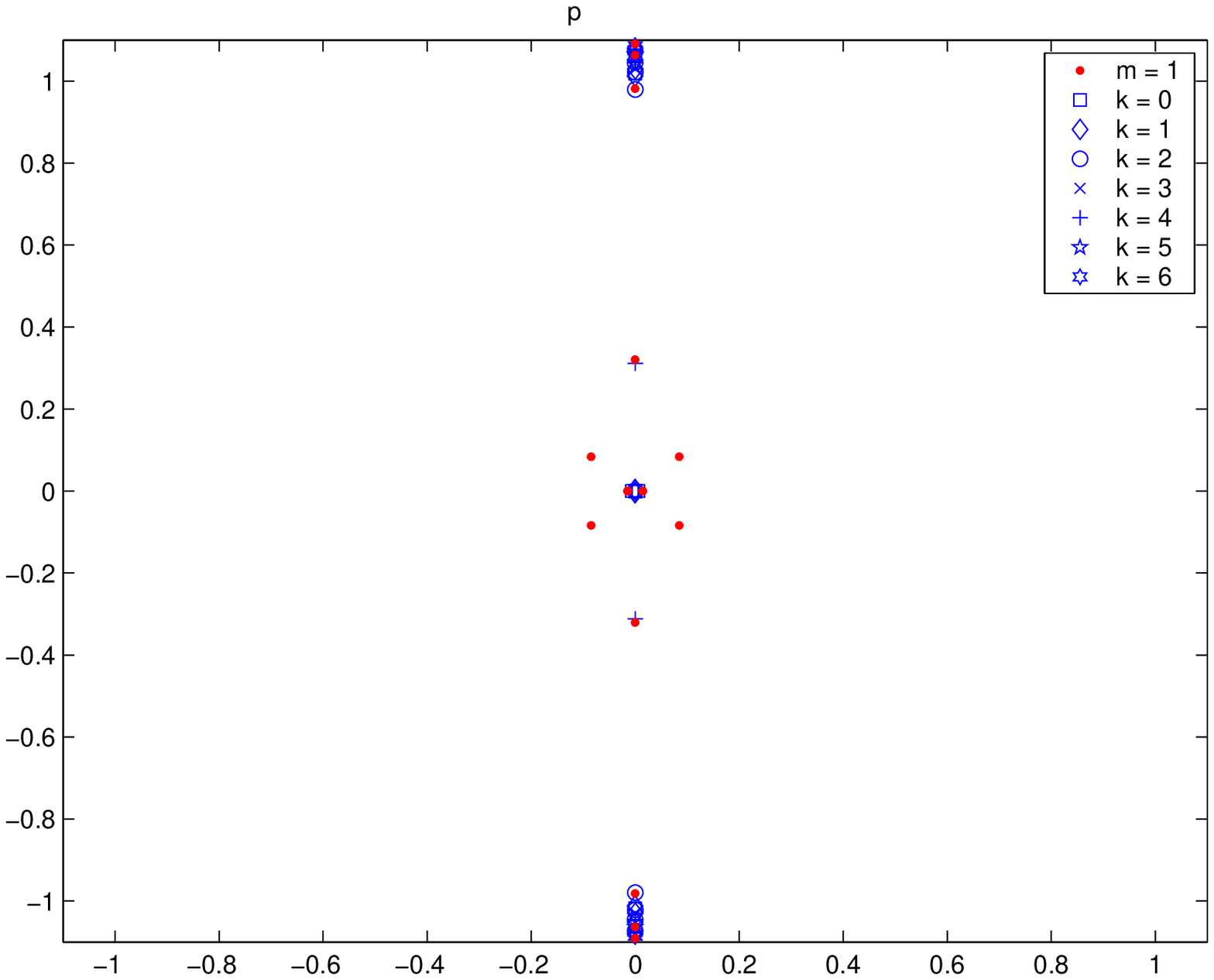,height=7.5cm,width=7.5cm}
\psfrag{p}{$p=3.2$}
\epsfig{file=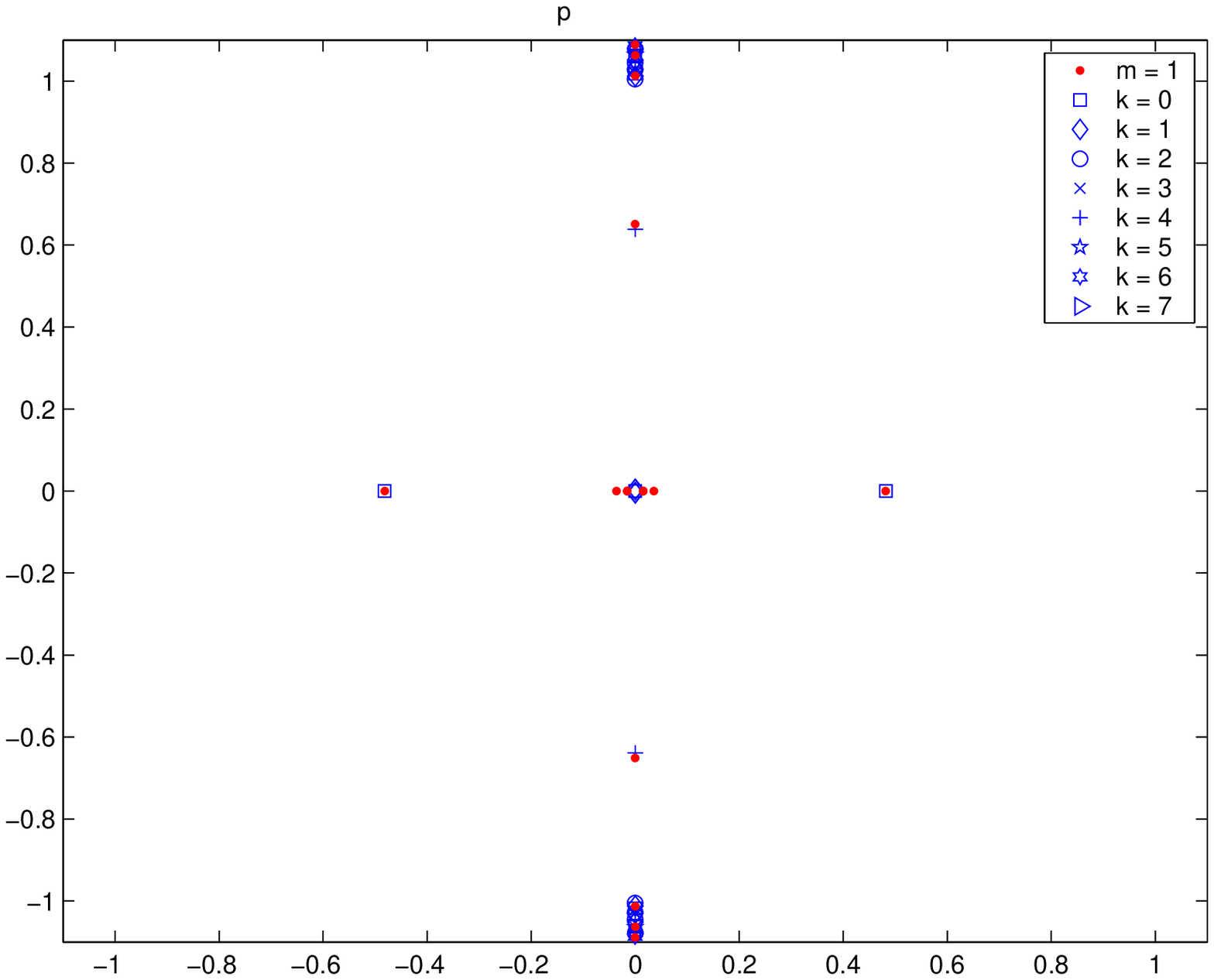,height=7.5cm,width=7.5cm}
 \caption{Spectra of $\L$ in $\R^2$ for $m=1$ and various
$p = 1.6, 2.1, 3, 3.2$. Point ``$\cdot$'' denotes the spectra
computed by Algorithm~1 and the others symbols denote the spectra
computed by Algorithms~2 and 3.} \label{fig:compm1}
\end{figure}

\begin{figure}[ht]
\psfrag{p}{$p=1.6$}
\epsfig{file=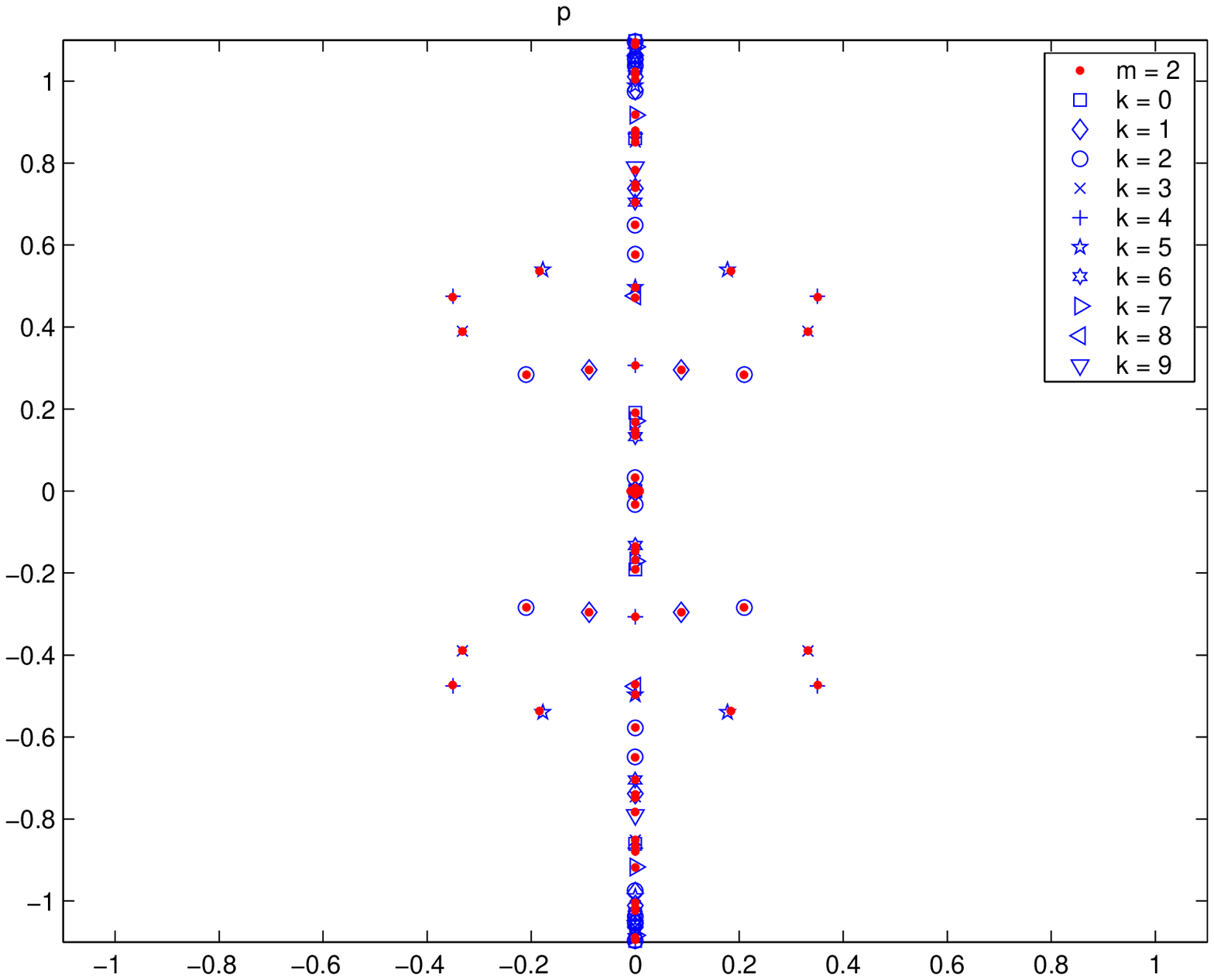,height=7.5cm,width=7.5cm}
\psfrag{p}{$p=2.1$}
\epsfig{file=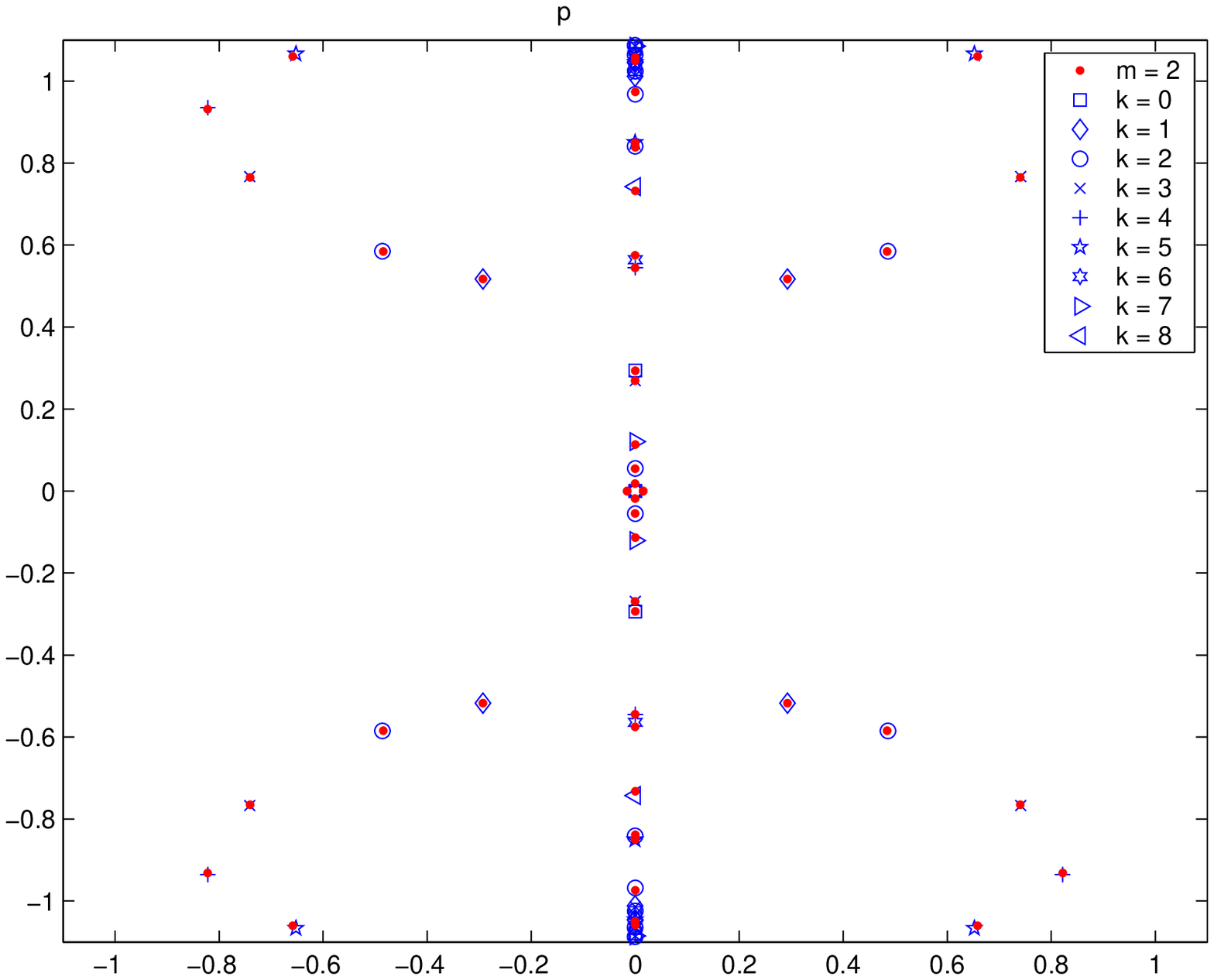,height=7.5cm,width=7.5cm}
\vspace*{0.5cm} \psfrag{p}{$p=3$}
\epsfig{file=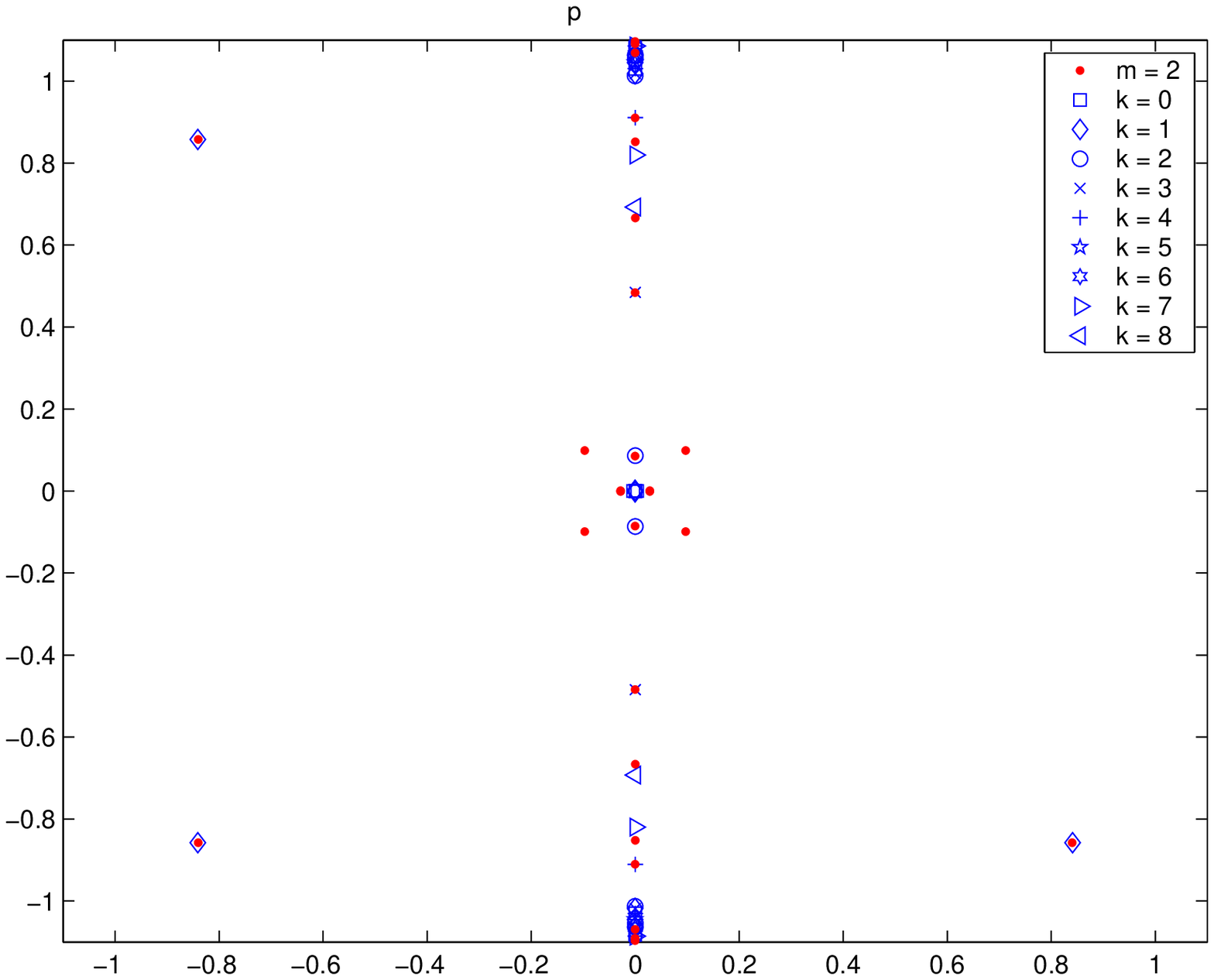,height=7.5cm,width=7.5cm}
\psfrag{p}{$p=3.2$}
\epsfig{file=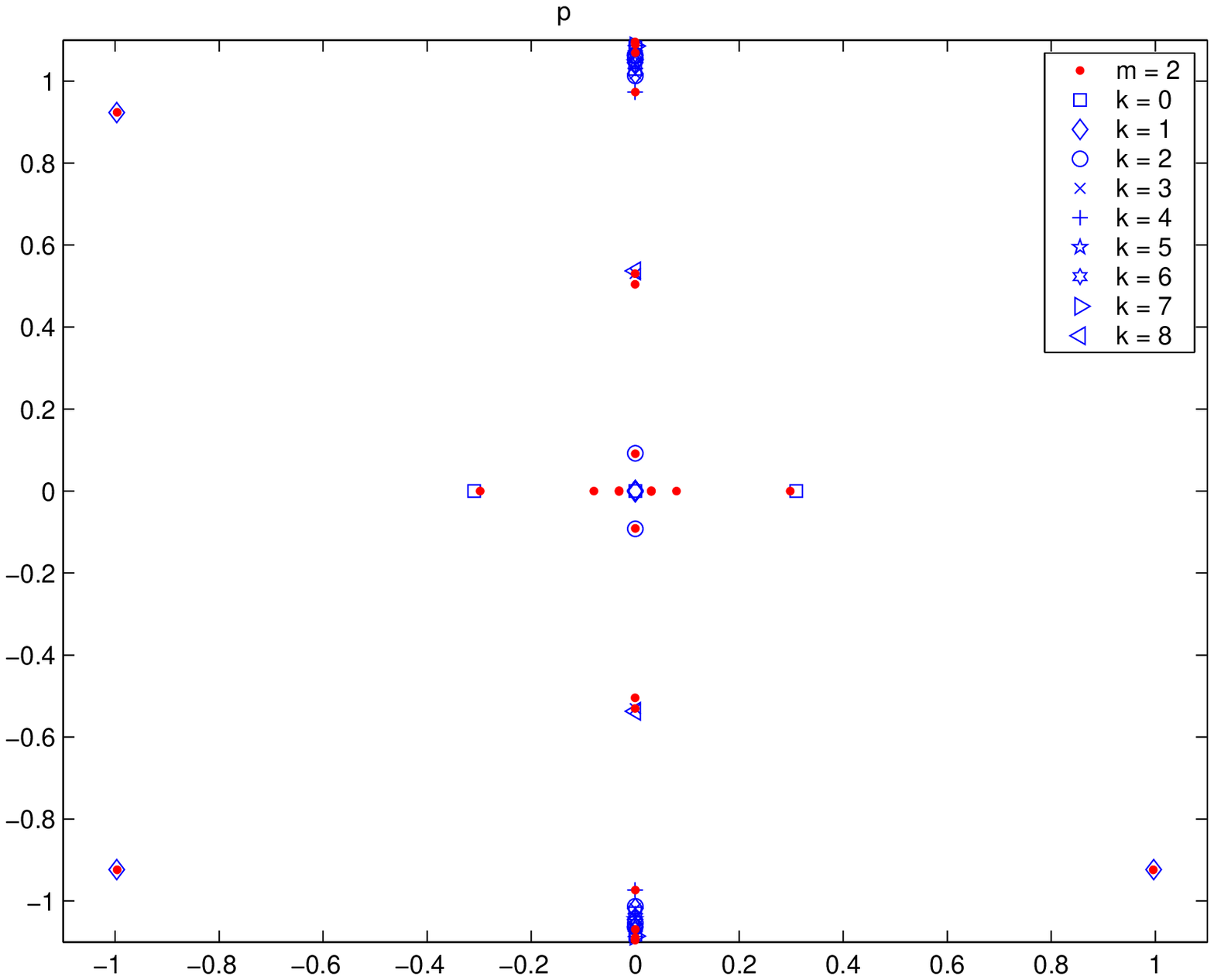,height=7.5cm,width=7.5cm}
 \caption{Spectra of $\L$ in $\R^2$ for $m=2$ and various
$p = 1.6, 2.1, 3, 3.2$. Point ``$\cdot$'' denotes the spectra
computed by Algorithm~1 and the others symbols denote the spectra
computed by Algorithms~2 and 3.} \label{fig:compm2}
\end{figure}

\begin{figure}[ht]
\psfrag{p}{$k=2$}\psfrag{a}{$p=1.01$}\psfrag{b}{$p=1.0165$}\psfrag{c}{$p=1.0
5$}
\epsfig{file=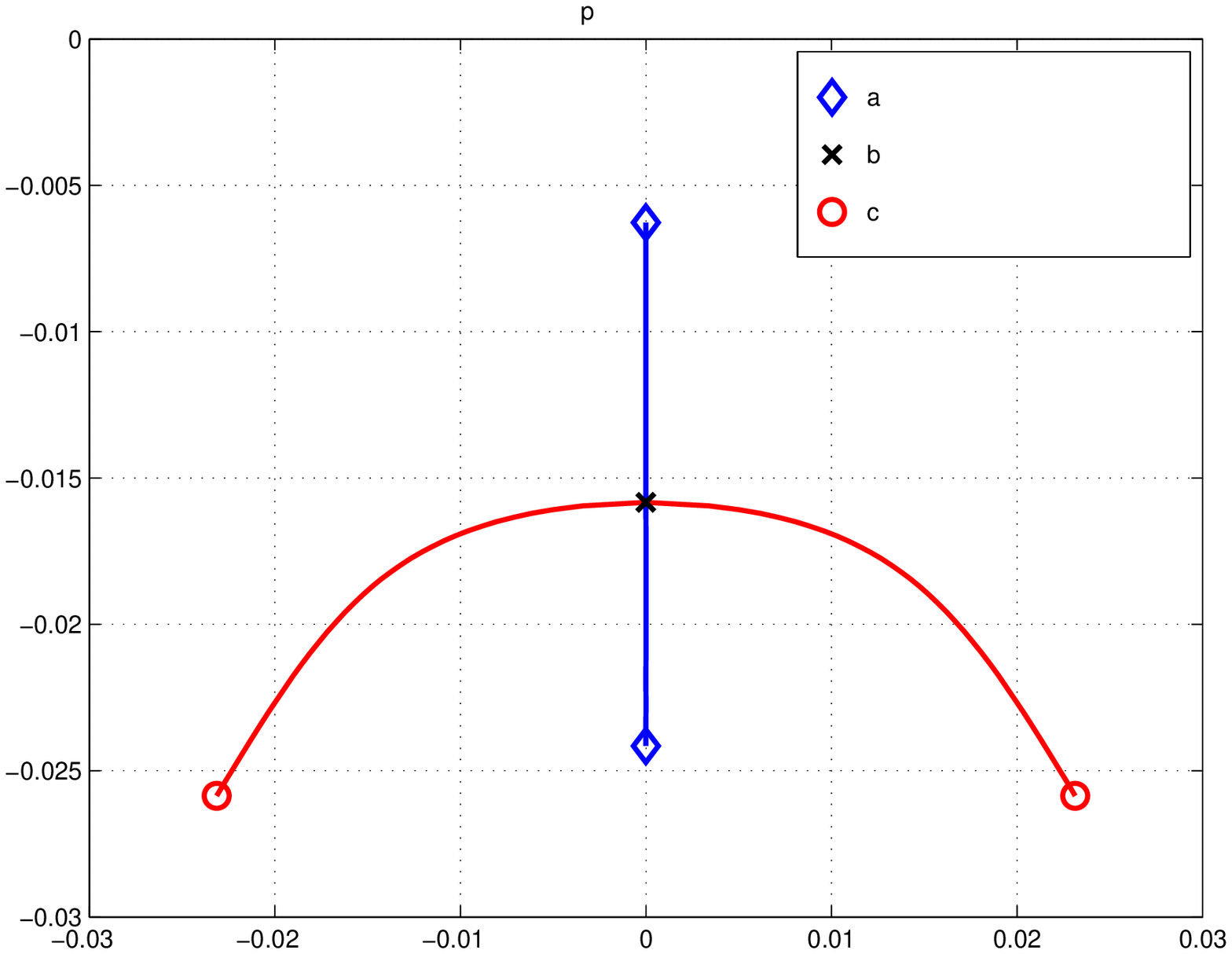,height=7.5cm,width=7.5cm}
\psfrag{p}{$k=3$}\psfrag{a}{$p=1.3$}\psfrag{b}{$p=1.3495$}\psfrag{c}{$p=1.4$
}
\epsfig{file=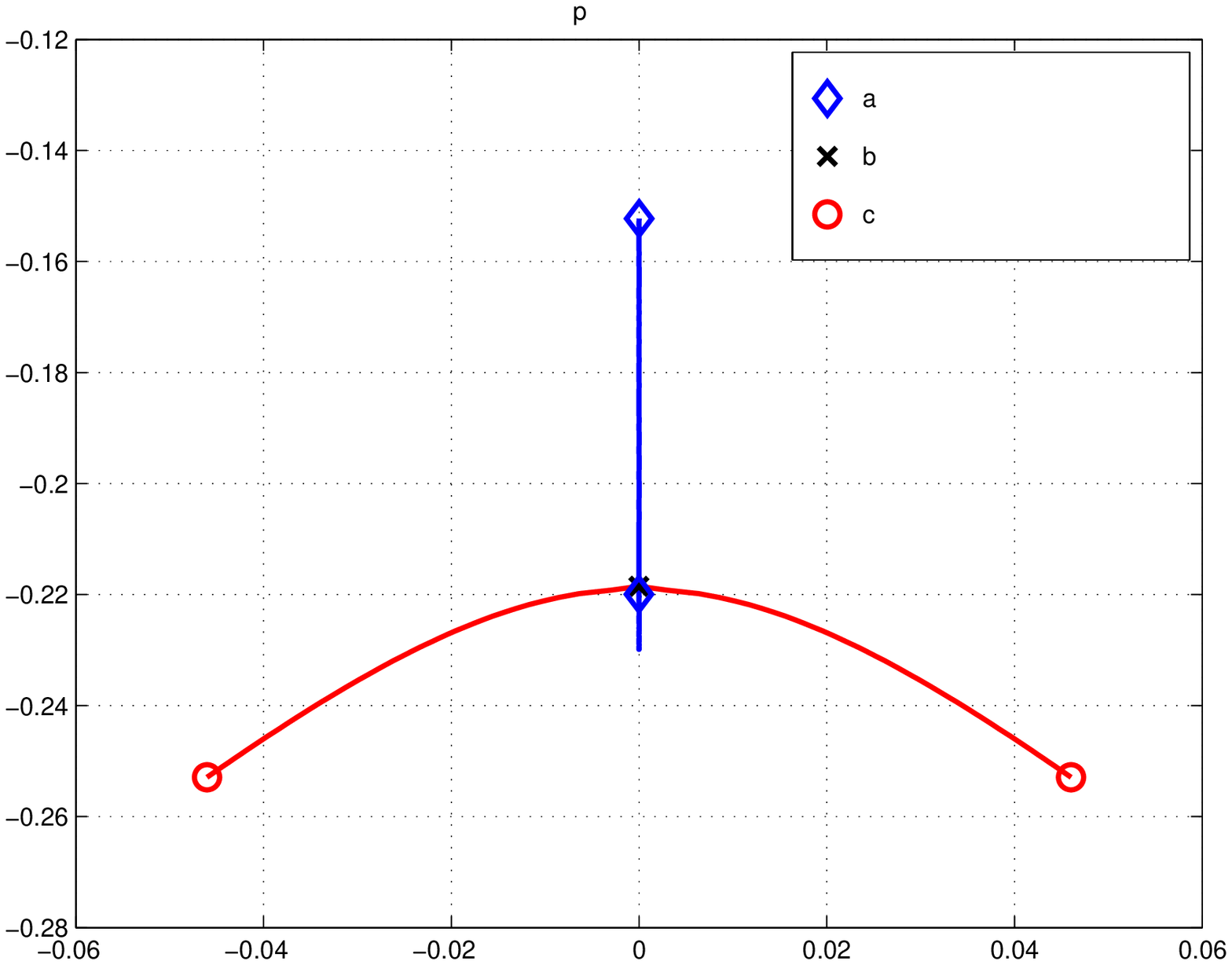,height=7.5cm,width=7.5cm}
\vspace*{0.5cm}
\psfrag{p}{$k=1$}\psfrag{a}{$p=1.48$}\psfrag{b}{$p=1.5265$}\psfrag{c}{$p=1.5
6$}
\epsfig{file=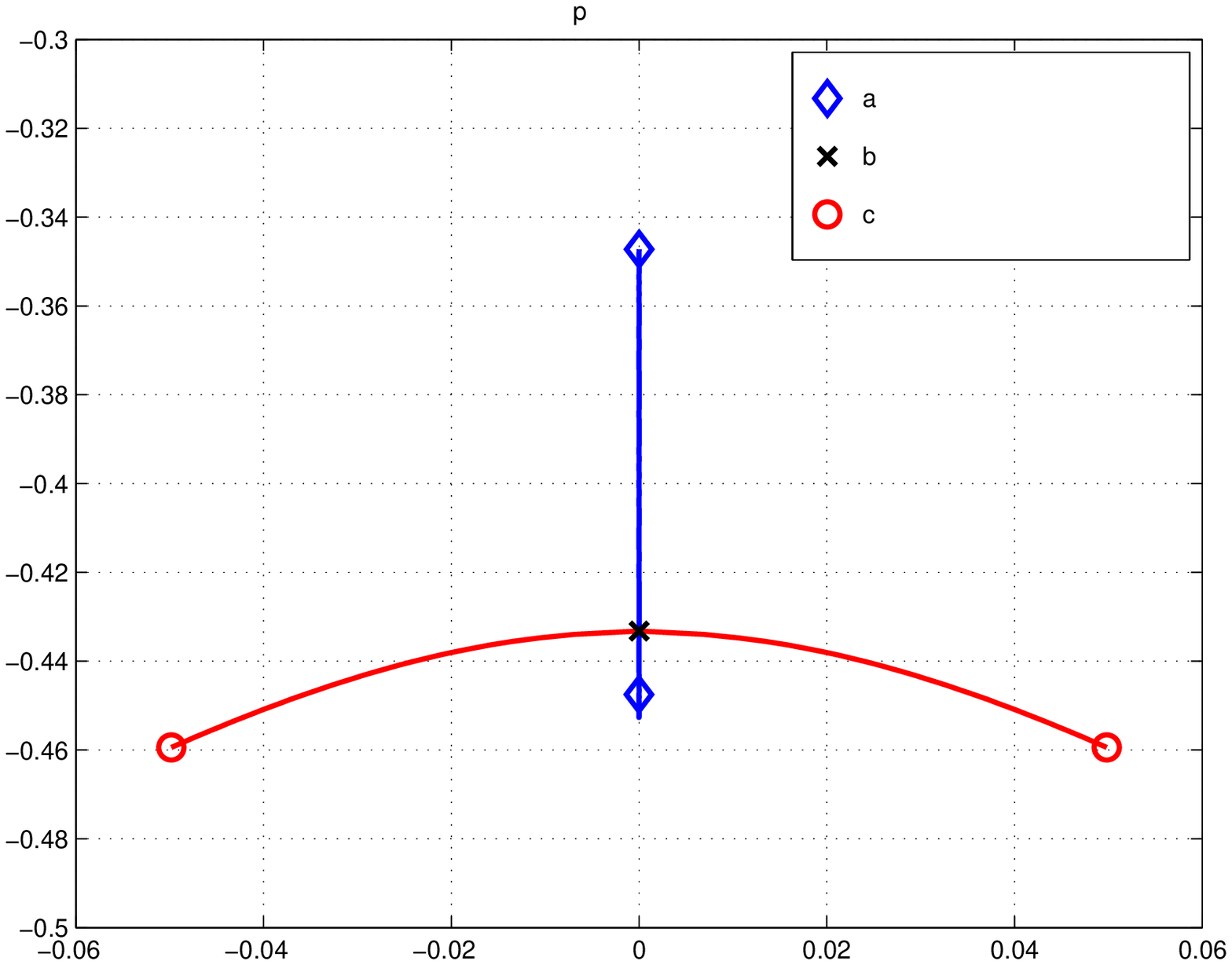,height=7.5cm,width=7.5cm}
\psfrag{p}{$k=0$}\psfrag{a}{$p=2.997$}\psfrag{b}{$p=3$}\psfrag{c}{$p=3.003$}
\epsfig{file=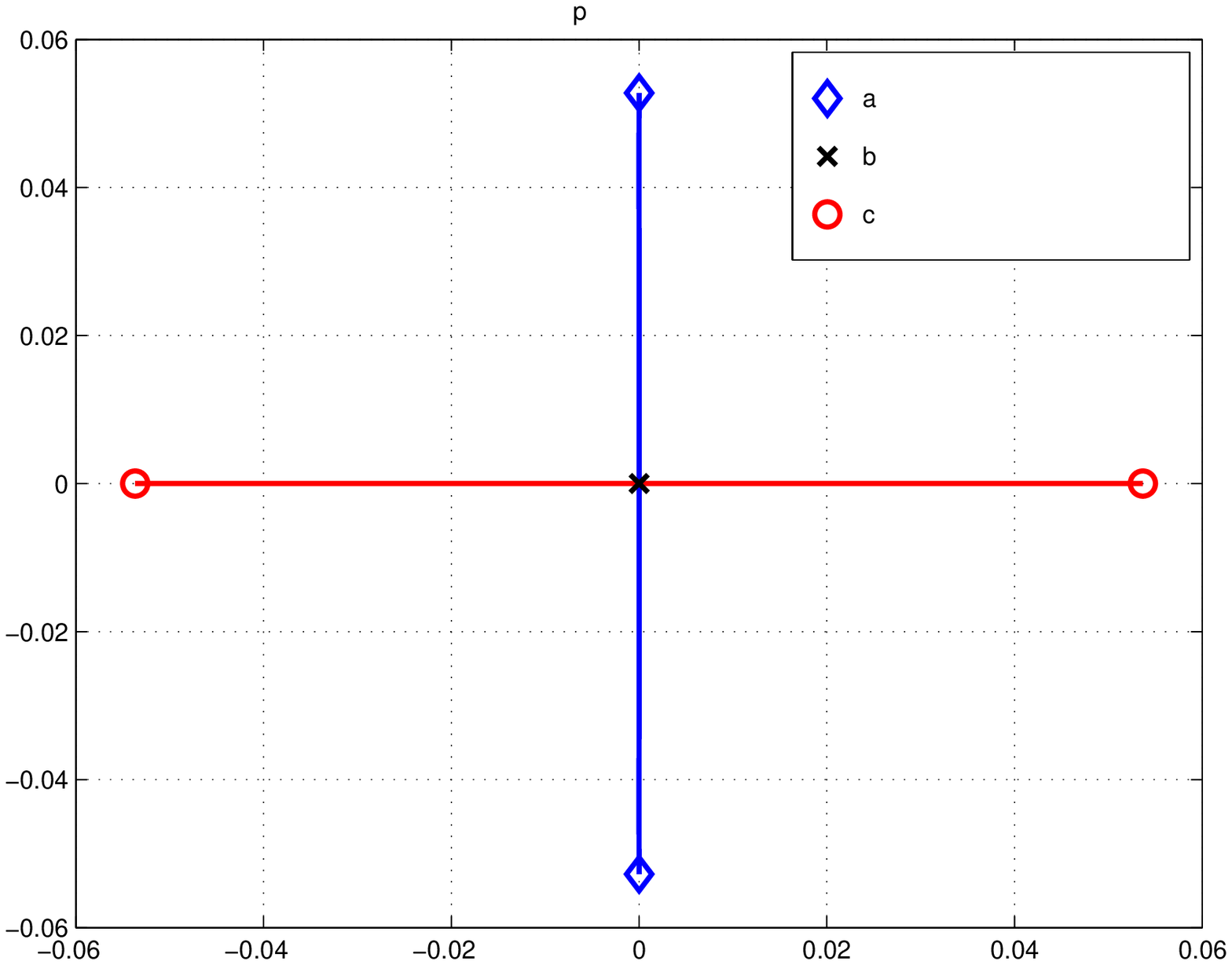,height=7.5cm,width=7.5cm}
 \caption{Bifurcation diagrams of $L_{m,k}$ for $m=1$ and $0\le k\le
 3$.} \label{fig:bifm1}
\end{figure}

\begin{figure}[ht]
\psfrag{p}{$k=2$}\psfrag{a}{$p=1.001$}\psfrag{b}{$p=1.0085$}\psfrag{c}{$p=1.
03$}
\epsfig{file=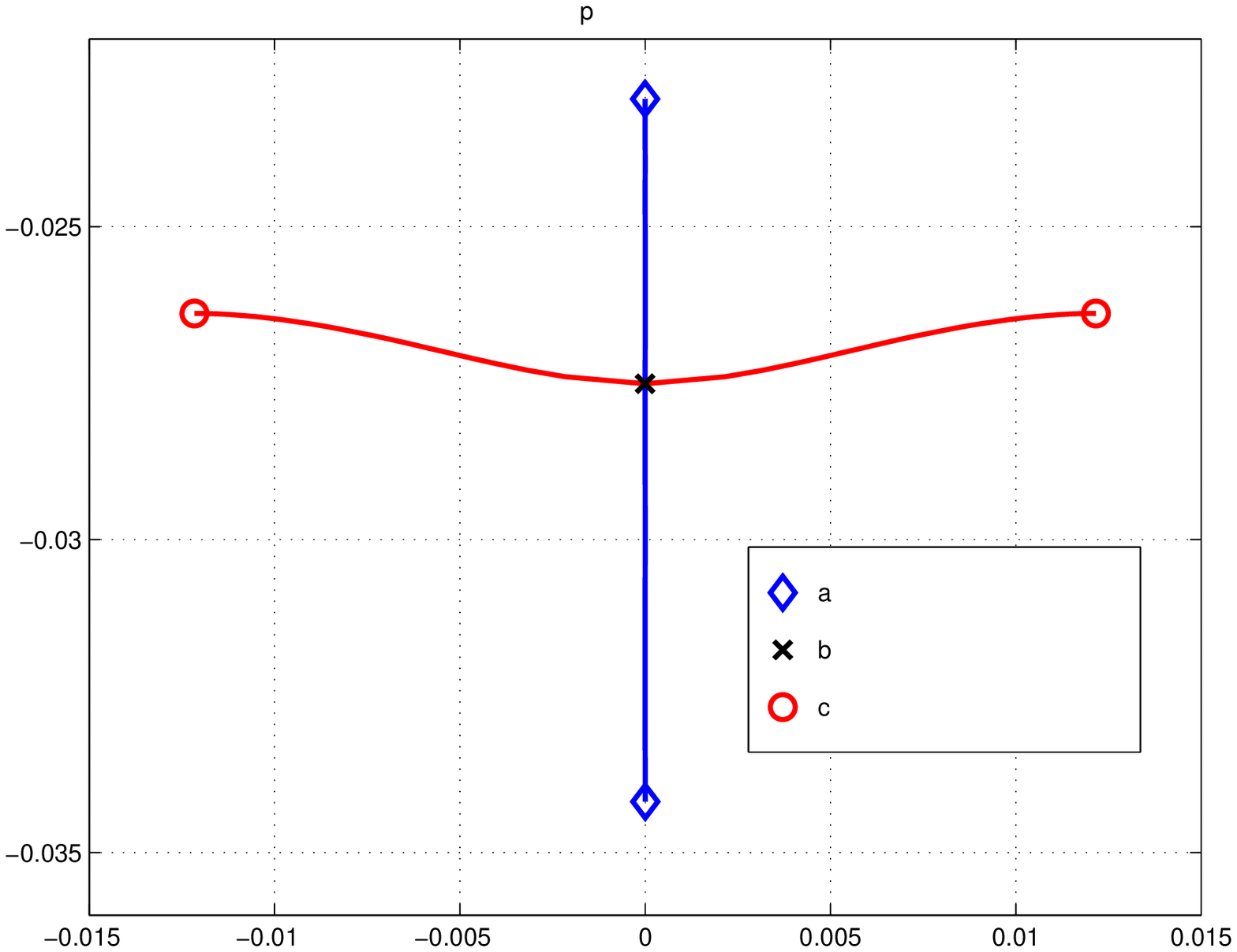,height=7.5cm,width=7.5cm}
\psfrag{p}{$k=3$}\psfrag{a}{$p=1.01$}\psfrag{b}{$p=1.0245$}\psfrag{c}{$p=1.1
$}
\epsfig{file=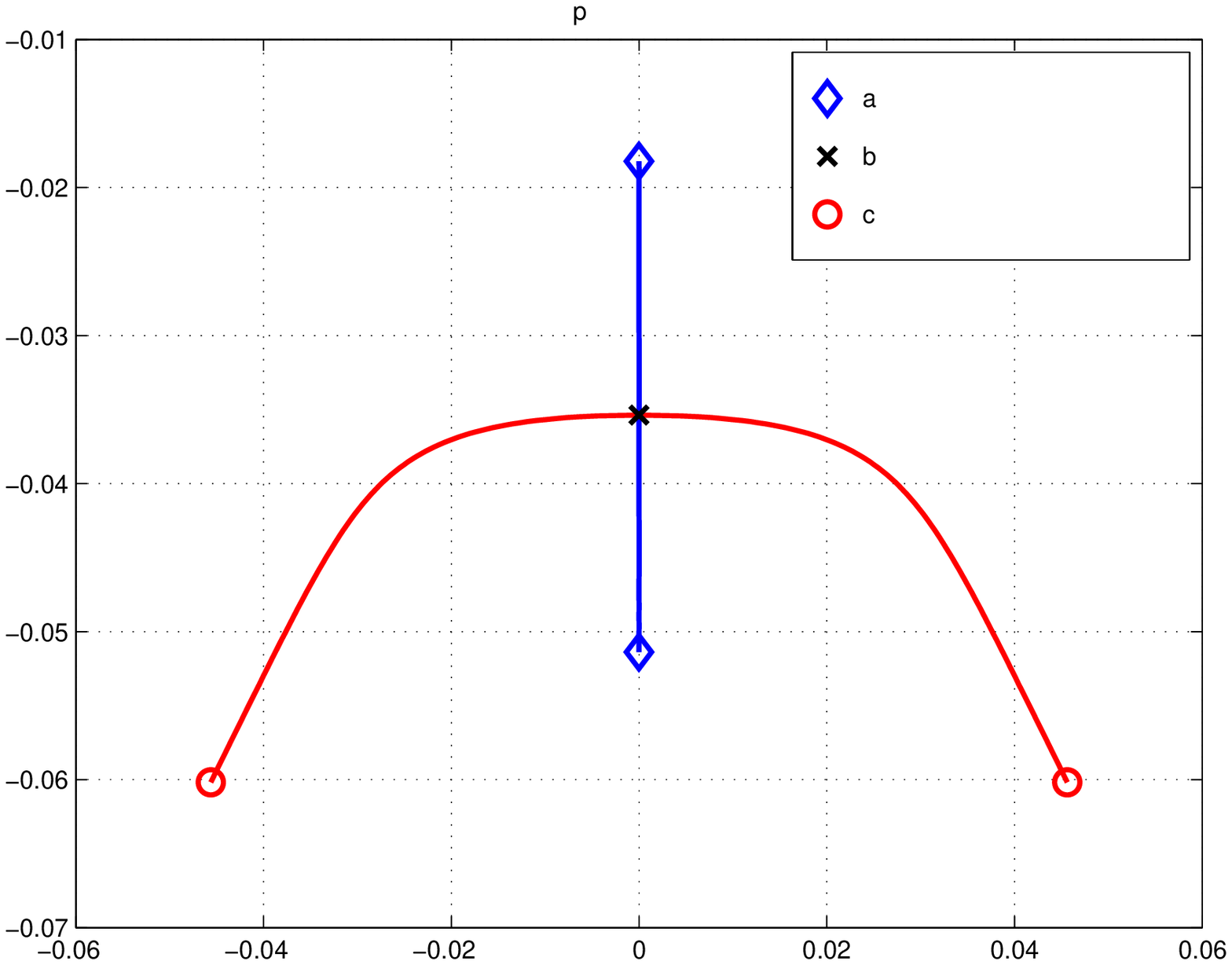,height=7.5cm,width=7.5cm}\vspace*{0.5cm
}
\psfrag{p}{$k=4$}\psfrag{a}{$p=1.01$}\psfrag{b}{$p=1.0455$}\psfrag{c}{$p=1.1
$}
\epsfig{file=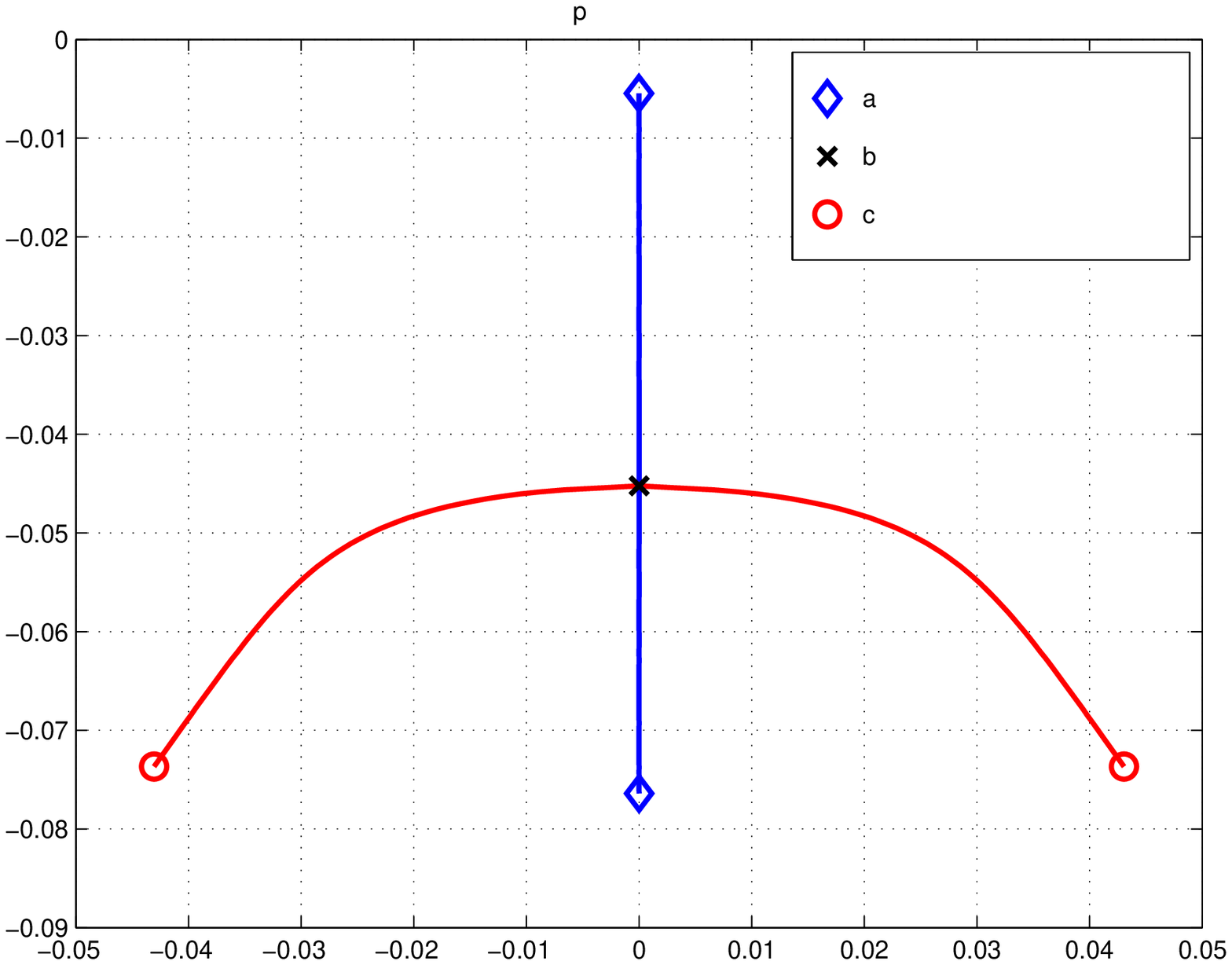,height=7.5cm,width=7.5cm}
\psfrag{p}{$k=1$}\psfrag{a}{$p=1.2$}\psfrag{b}{$p=1.3565$}\psfrag{c}{$p=1.5$
}
\epsfig{file=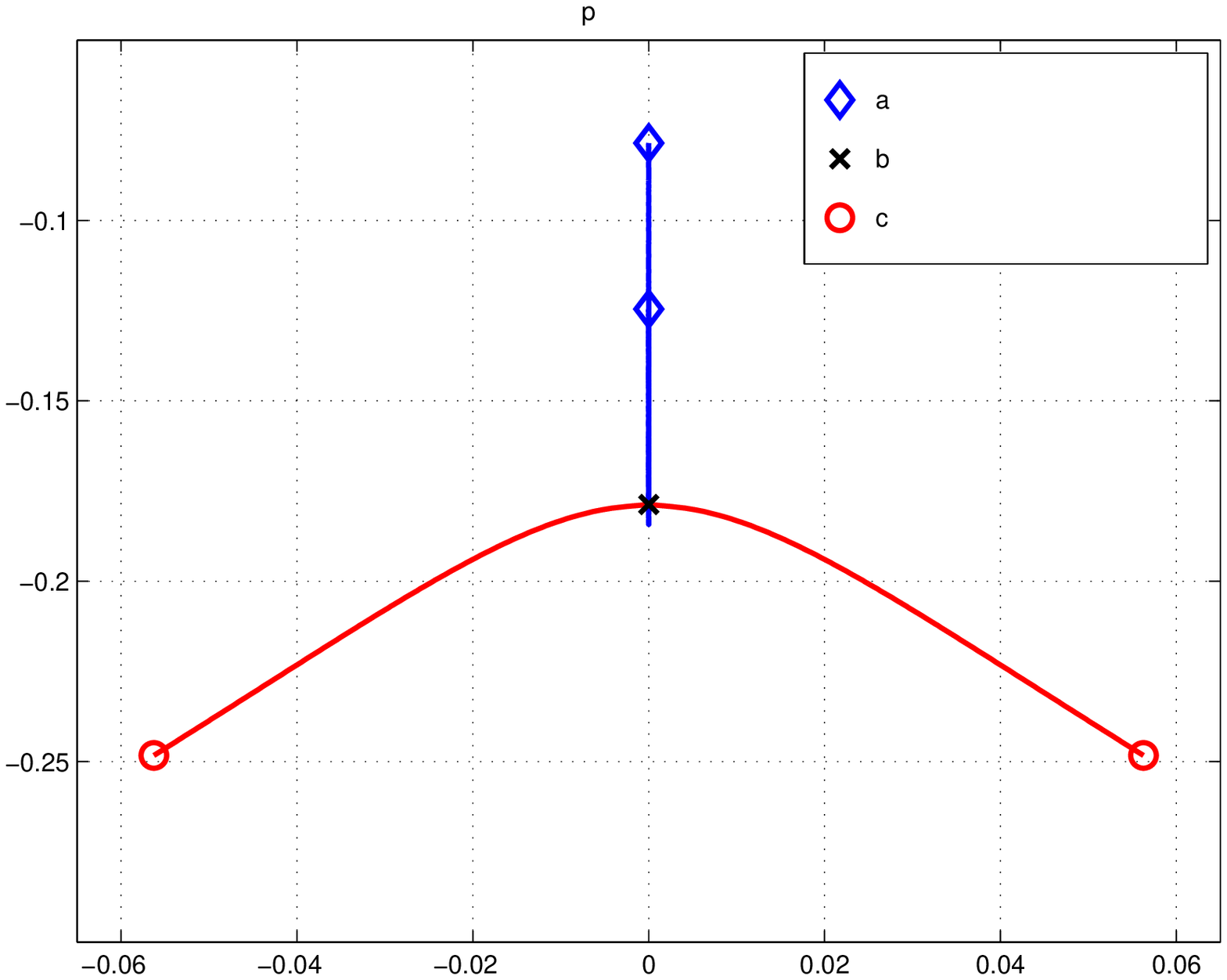,height=7.5cm,width=7.5cm}\vspace*{0.5cm
}
\psfrag{p}{$k=5$}\psfrag{a}{$p=1.355$}\psfrag{b}{$p=1.3955$}\psfrag{c}{$p=1.
435$}
\epsfig{file=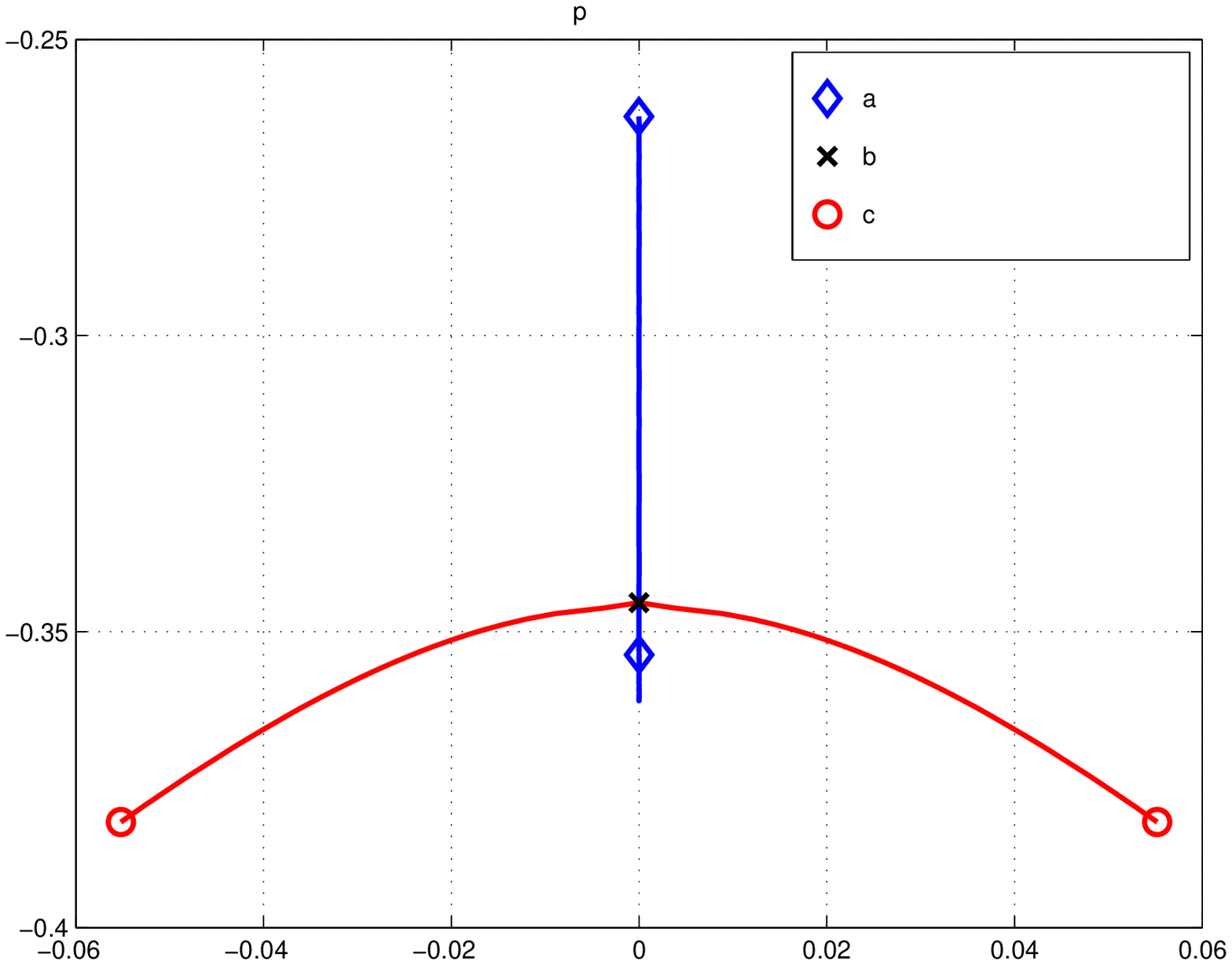,height=7.5cm,width=7.5cm}
\psfrag{p}{$k=0$}\psfrag{a}{$p=2.993$}\psfrag{b}{$p=3$}\psfrag{c}{$p=3.007$}
\epsfig{file=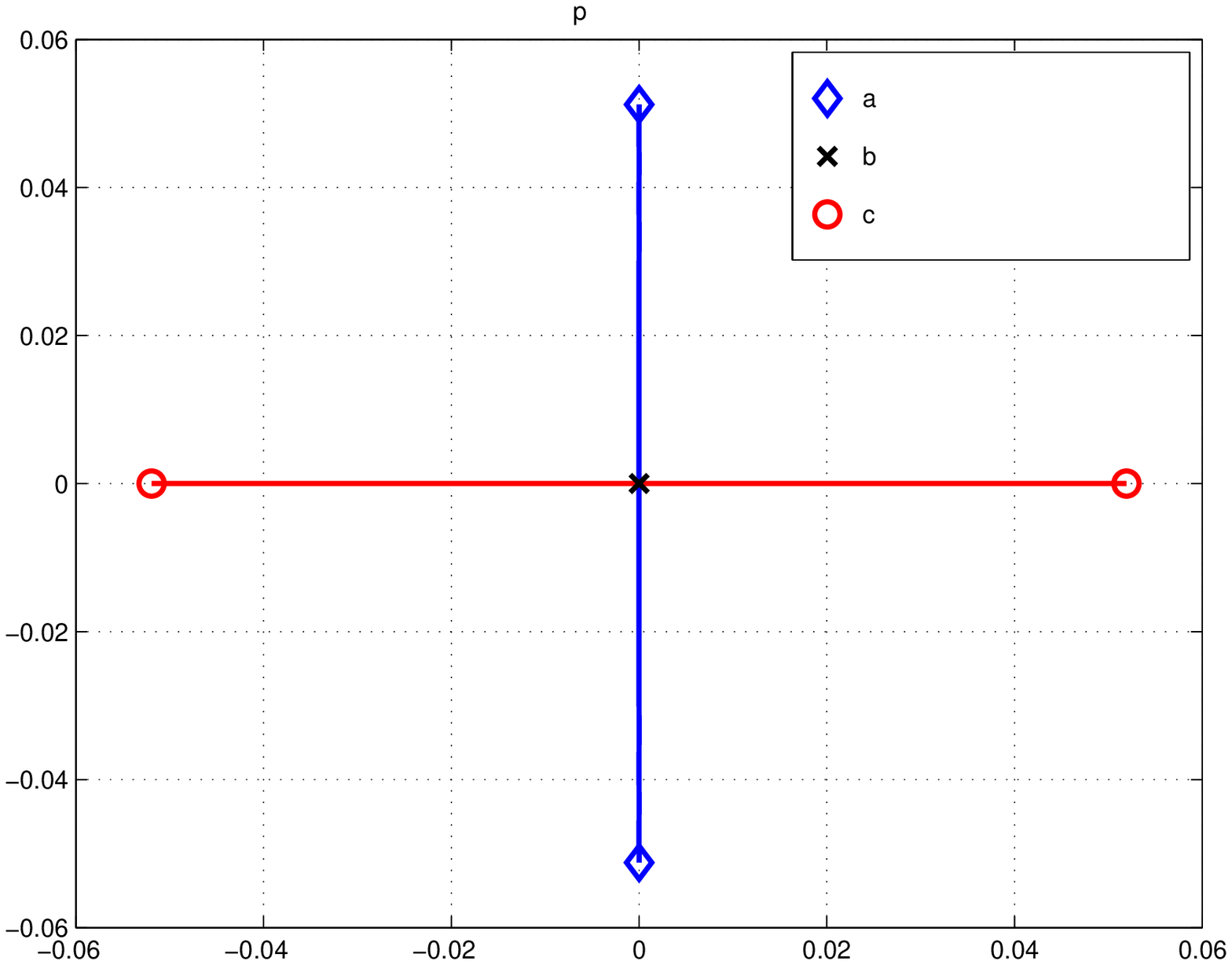,height=7.5cm,width=7.5cm}
 \caption{Bifurcation diagrams of $L_{m,k}$ for $m=2$ and $0\le k\le
 5$.} \label{fig:bifm2}
\end{figure}

The results of our numerical computations of the spectra of $\L$ for
$m=1,2$ and various $k$ and $p$ are shown in
Figures~\ref{fig:d2m1}--\ref{fig:bifm2}. As before, we focus on
eigenvalues in the square $\bket{a+bi: |a|<1, |b|<1}$. Purely
imaginary eigenvalues with modulus greater than 1 correspond to the
continuous spectrum of $\L$, and are discrete due to discretization.

Let us first describe some simple observations:
\begin{enumerate}

\item The distribution of eigenvalues, see
Figures~\ref{fig:d2m1}--\ref{fig:d2m2}, is more complicated and
interesting than Figures~\ref{fig:cL1d}--\ref{fig:cL2d-m2}.  There are
not only purely imaginary eigenvalues and real eigenvalues but also
complex eigenvalues, whose existence implies instability.

\item
In Figures~\ref{fig:compm1}--\ref{fig:compm2}, we compare the results
obtained from three algorithms. For Algorithms~2 and 3 the parameter
$k$ ranges from $0$ to $9$. Results from these algorithms have high
degree of agreement, except when $p$ near $3$ and eigenvalues near
$0$. We will discuss this exceptional case in the end.

\item For Algorithms~2 and 3, the numerical $0$-eigenvalue occurs only
when $k=0$ and $k=\pm 1$. It agrees with our discussion in the
previous subsection. Their multiplicities also match and there is no
unaccounted eigenvector. In particular, $N_g(\L)$ has dimension $6$ if
$p \not= 3$ and $8$ if $p=3$, the same as the case of ground
states. We also numerically verified the nullspace, for example, the
discrete version of \eqref{eq4-test} is correct.

\item As $p$ increases,  two pairs of purely imaginary eigenvalues may
collide away from $0$, and then split to a quadruple of complex
eigenvalues which are neither real nor purely imaginary.  For $m=1$,
this bifurcation phenomenon appears three times before $p=1.55$ and
there are 3 complex quadruples for $p>1.55$. For $m=2$, it occurs five
times before $p=1.5$ and there are 5 quadruples for $p>1.5$.  These
complex eigenvalues seem to move away from the imaginary axis as $p$
increases further.

\item As $p$ increases to $3$, (by Algorithms 2 and 3), a pair of
purely imaginary eigenvalues from $0$-th subspace collide at $0$ and
then split to a pair of real eigenvalues as $p$ increases
further. This is the same picture as in the ground state case in section
3.  Indeed, Mizumachi \cite{Mi1} proves that $Q_{m,0,p}$ are stable in
the $0$-th subspace if $p<3$ and unstable if $p>3$. Thus $p=3$ is a
bifurcation point. Also note that when $p=3$ the NLS \eqref{eq1-1} has
conformal invariance and explicit blow-up solutions can be found as in
the ground state case.

\item In Figures~\ref{fig:bifm1}--\ref{fig:bifm2} we observe the
bifurcation more closely.  For $m=1$, the bifurcation occurs when
$(k,p,\la)$ equal
\[
(0,\,3,\,0), \quad (1,\, 1.52765\, -0.436i), \quad (2,\, 1.0165,\,
-0.016i), \quad (3,\, 1.3495,\, -0.219i).
\]
For $m=2$, the bifurcation occurs when $(k,p,\la)$ equal
\[
(0,\,3,\,0), \quad (1,\, 1.357,\, -0.180i), \quad (2,\, 1.007,\,
-0.027i), \quad (3,\, 1.0245,\, -0.035i) \] and \[ (4,\, 1.0455,\,
-0.045i), \quad (5,\, 1.3955,\, -0.347i).
\]

\item Due to the existence of complex eigenvalues for $m=1,2$ and
$p\ge 1.02$, $Q_{m,0,p}$ is spectrally unstable for these
parameters. However, all these complex eigenvalues bifurcate from some
discrete eigenvalues $\pm bi$ with $|b|<1$ and $p>1.008$. 
Our computation for both $m=1,2$ and
\[
p= \ell \cdot 0.001, \quad \ell = 1,2,3,\ldots , 8, \quad \text{(up to
  15 if }m=1)
\]
does not find any complex eigenvalues.  This suggests that the two
excited states $\phi_{1,0,p}(r)e^{i \th}$ and $\phi_{2,0,p}(r)e^{i2
\th}$ are {\it linearly stable} when $p$ is sufficiently close to 1.
It is possible that the numerical error increases enormously as $p \to
1_+$ due to the artificial boundary condition, since the spectrum is
approaching to the continuous one for $p=1$. This has to be verified
analytically in the future.

\end{enumerate}

We finally discuss the exceptional case when $p$ is near $3$ for
eigenvalues near $0$.  In this case Algorithm~1 produces an quadruple
of complex eigenvalues $\pm 0.0849 \pm 0.0836i$, and the
$0$-eigenvalue has multiplicity 4.  We expect to see larger errors
from Algorithm~1 but the error in this exceptional case is much
larger.  It is related to the large size of a Jordon block for the $0$
eigenvalue. As discussed in the previous subsection, the nullspace is
at least 6 dimensional. The analysis in section \ref{sec2-3} suggests
that (we do not claim a proof), as $p$ goes to the bifurcation
exponent $p_c=3$ from below, a pair of imaginary eigenvalues merges
into the Jordon block containing the eigenfunctions $iQ$ and $Q_1$,
and the Jordon block becomes size 4.  As is well-known in matrix
analysis (see \cite[p.324]{GvL}, \cite{Wilkinson}), if a matrix
contains a Jordon block of size $\ell$, the computed eigenvalues
corresponding to that block have errors of order $\e^{1/\ell}$, where
$\e$ is the sum of the machine zero, the truncation error from
discretization, and the perturbation (from varying $p$). Since $\de_r
= 0.04$ for Algorithm~1 and the truncation error of a central
difference scheme for $\Delta_r$ has order $ O(\de_r^2)$, the error
for the zero eigenvalue near $p=3$ could be
\[
(\de^2_r)^{1/4} \approx 0.2.
\]
In contrast, for other bifurcation points on the imaginary axis, the
Jordon block at the bifurcation exponent is of size 2 and the error is
of order $(\de^2_r)^{1/2}=0.04$. In practice, the error is smaller due
to cancellation and the numerical results by Algorithm~1 do not differ
too much from those by Algorithms~2 and 3.  Also note that numerically
the $0$ eigenspace has dimension~4, accounting for $ Q_{x_j}$ and
$ix_jQ$. The complex quadruple correspond to $iQ$, $Q_1$ and the
joining pair of nonzero eigenvalues.


\section*{Appendix: Numerical method}

In this section we describe a numerical method to compute the
spectrum of the linear operator ${\cal L}$ defined by (\ref{eq1-9})
for $p>1$ and space dimension $n\ge 1$. There are two main steps in
this method. First, we will solve the nonlinear problem \eqref{Q.eq}
for $Q$: we will discretize it into a nonlinear algebraic equation,
and then solve it by an iterative method. Second, we will compute
the spectrum of $\L$: we will discretize the operator $\L$ into a
large-scale linear algebraic eigenvalue problem and then use
implicitly restarted Arnoldi methods to deal with this problem.

Hereafter, we use the bold face letters or symbols to denote a
matrix or a vector.  For ${\mathbf A}\in{\mathbb R}^{M\times N}$,
${\mathbf q}=(q_1,\ldots,q_N)^\top\in{\mathbb R}^N$, ${\mathbf
q}^\pcirc={\mathbf q}\circ\cdots\circ{\mathbf q}$ denotes the
$p$-time Hadamard product of ${\mathbf q}$, and $\lb{\mathbf
q}\rb:=\diag({\mathbf q})$ the diagonal matrix of ${\mathbf q}$.

{\bf Step I.} \ We first discretize equation \eqref{Q.eq} into a
nonlinear algebraic equation and consider it on an $n$-dimensional
ball $\Omega=\{ x\in {\mathbb R}^n:| x|\le R, R\in {\mathbb R}\}$.
We rewrite the Laplace operator $-\Delta$ in the polar coordinate
system with a Dirichlet boundary condition. Based on the recently
proposed discretization scheme \cite{Lai:2001}, the {\it standard
central finite difference method}, we discretize $-\Delta{\mathbf
q}({x})$ into
\begin{equation}\label{eq5-1}\tag{A.1}
{{\mathbf A}}{\mathbf q}={{\mathbf
A}}[q_{1},\ldots,q_{N}]^\top,~{{\mathbf A}}\in \mathbb{R}^{N\times
N},
\end{equation}
where ${\mathbf q}$ is an approximation of the function $Q(x)$. The
matrix ${{\mathbf A}}$ is irreducible and diagonally-dominant with
positive diagonal entries.  The discretization of the nonlinear
equation (\ref{Q.eq}) can now be formulated as the following
nonlinear algebraic equation,
\begin{equation}\label{eq5-2}\tag{A.2}
{{\mathbf A}}{\mathbf q} + {\mathbf q} -{\mathbf q}^\pcirc = 0.
\end{equation}
We introduce an iterative algorithm \cite{HW} to solve
(\ref{eq5-2}):
\begin{equation}\label{eq5-3}\tag{A.3}
{{\mathbf A}}\widetilde{\mathbf q}_{j+1} + \widetilde{\mathbf
q}_{j+1} ={\mathbf q}_{j}^\pcirc,
\end{equation}
where $\widetilde{\mathbf q}_{j+1}$ and ${\mathbf q}_{j}$ are the
unknown and known discrete values of the function $Q({\mathbf x})$,
respectively. The iterative algorithm is shown below.\vspace*{0.5cm}

\noindent\hspace*{-0.3cm}\textbf{Iterative Algorithm for Solving
$Q({\mathbf x})$}.
\begin{enumerate}\itemsep=0pt
\item[Step 0]   Let $j=0$.\\
         Choose an initial solution $\widetilde{{\mathbf q}}_{0}>0$ and let $%
         {\mathbf q}_{0}=\frac{\widetilde{{\mathbf q}}_{0}}{\parallel \widetilde{{\mathbf q}}_{0}\parallel_2 }.$
\item[Step 1]    Solve the equation (\ref{eq5-3}), then obtain $\widetilde{\mathbf q}_{j+1}$.

\item[Step 2]    Let $\alpha _{j+1}=\frac{1}{\parallel \widetilde{{\mathbf
q}}_{j+1}\parallel_{2}}$ and normalize $\widetilde{{\mathbf
q}}_{j+1}$ to obtain ${\mathbf q}_{j+1}=\alpha _{j+1}
\widetilde{{\mathbf q}}_{j+1}$.

\item[Step 3]    If (convergent) then\\
           \hspace*{0.5cm} Output the scaled solution ($\alpha _{j+1}$)$^{\frac{1}{p-1}}{\mathbf q}_{j+1}.$
           Stop.\\
         else\\
           \hspace*{0.5cm} Let $j:=j+1.$\\
          \hspace*{0.5cm} Goto Step 1.\\
         end 
\end{enumerate}


If the components of ${\mathbf q}_0$ are nonnegative, this property
is preserved by each iteration ${\mathbf q}_{j}$, and hence also by
the limit vector if it exists (see \cite[Theorem~3.1]{HW}).  The
convergence of a subsequence of this iteration method to a nonzero
vector is proved in \cite[Theorem~2.1]{HW}.  Although the
convergence of the entire sequence is not proved, it is observed
numerically to be very robust. See Chen-Zhou-Ni \cite{CZN} for a
survey on numerically solving nonlinear elliptic equations.

\medskip

{\bf Step II.} \ Next we discretize the operator ${\cal L}$ of
(\ref{eq1-11}) into a linear algebraic eigenvalue problem:
\begin{equation}\label{eq5-4}\tag{A.4}
{\mathbf L}\left[\begin{array}{c}
                         {\mathbf u}\\
                        {\mathbf w}
                    \end{array}
              \right]=\lambda \left[\begin{array}{c}
                         {\mathbf u}\\
                        {\mathbf w}
                    \end{array}\right],
\end{equation}
where
\begin{equation*}\label{eq5-5}
{\mathbf L} = \left[\begin{array}{cc}
                        0 & {\mathbf A}+{\rm I}-\lb{\mathbf q}^\ppcirc\rb\\
                        -{\mathbf A}-{\rm I}+\lb p~\!{\mathbf q}^\ppcirc\rb  &0
                    \end{array}
              \right],
\end{equation*}
$\gamma = p-1$, ${\mathbf u}=(u_1,\ldots,u_N)^\top\in{\mathbb R}^N$,
${\mathbf w}=(w_1,\ldots,w_N)^\top\in{\mathbb R}^N$, and ${\mathbf
q}$ is the output of the previous step, and satisfies the equation
in (\ref{eq5-2}). We use ARPACK \cite{LSY98} in MATLAB version~6.5
to deal with the linear algebraic eigenvalue problem (\ref{eq5-4})
and obtain eigenvalues $\lambda$ of ${\mathbf L}$ near the origin
for $p>1$ and space dimension $n\ge 1$. Furthermore, the
eigenvectors of ${\mathbf L}$ can be also produced.

The Step II above can in principle be used to compute all
eigenfunctions in $L^2(\R^n)$. However, in producing
Figures~\ref{fig:cL2d}--\ref{fig:cL2d-m2}, we look for
eigenfunctions of the form $\phi(r) e^{im\th}$. These problems can
be reformulated as 1-D eigenvalue problems for $\phi(r)$, which can
be computed using the same algorithm and MATLAB code. This
dimensional reduction saves a lot of computation time and memory.
Even with this dimensional reduction, and applying an algorithm for
sparse matrices, the computation is still very heavy, and we cannot
compute all eigenvalues in one step. We can only compute a portion
of them each time.

\section*{Acknowledgments}
Special thanks go to Wen-Wei Lin, whose constant advice and interest
were indispensable for this project.  We thank V.S.~Buslaev,
Y. Martel, F.~Merle, T.~Mizumachi, C.~Sulem, and W.C.~Wang for
fruitful discussions and providing many references to us.  We also
thank the referees for many valuable suggestions.  Part of the work
was done when Nakanishi was visiting the University of British
Columbia, Vancouver, and when Tsai was visiting the National Center
for Theoretical Sciences, Hsinchu \& Taipei. The hospitality of these
institutions are gratefully acknowledged.  The research of Chang is
partly supported by National Science Council in Taiwan.  The research
of Gustafson and Tsai is partly supported by NSERC grants. The
research of Nakanishi is partly supported by the JSPS grant
no.~15740086.


%
\newpage
\newpage

\clearpage
\newpage
\thispagestyle{empty}
\topmargin -0.5in

\end{document}